\documentclass[12pt,a4paper]{article}

\usepackage[margin=25.4mm]{geometry}

\usepackage[utf8]{inputenc}
\usepackage{xcolor,epsfig,amsmath,mathrsfs,amsfonts,bef_alex,%
            amssymb,datetime,bm}
\renewcommand{\FG}{\bm}
\renewcommand{\ti}{{\times}} 

\usepackage[final,color,notref,notcite]{showkeys}

\usepackage[utf8]{inputenc}
\usepackage{enumerate}

\usepackage{todonotes}

\allowdisplaybreaks

\def\sfx{\mathsf x}  \def\sfy{\mathsf y}  
  
\def\sfC{\mathsf C}  \def\sfE{\mathsf E}   
\def\sfK{\mathsf K}  \def\sfM{\mathsf M} 
\def\sfW{\mathsf W}  \def\sfX{\mathsf X}
\def\sfY{\mathsf Y}  \def\sfZ{\mathsf Z} 
\def\mfu{\mathfrak u}  \def\mfw{\mathfrak w}
  \def\mfB{\mathfrak B}  
  \def\mfE{\mathfrak E}
\def\mfK{\mathfrak K}  \def\mfN{\mathfrak N}  
\def\mfQ{\mathfrak Q}    

\def\hbA{{\bbA}}

\newcommand{\oti}{{\otimes}}
\newcommand{\LB}{\lambda_\rmB}
\newcommand{\Lio}{\mathrm{Lio}}
\newcommand{\CLE}{\mathrm{CLE}}
\def\div{\mathop{\mathrm{div}}\nolimits}
\def\tra{\mathsf{T}}

\newcommand{\FW}{\mathrm{fw}}
\newcommand{\BW}{\mathrm{bw}}
\def\bflog{\mathop{\mathbf{log}}}

\def\sfc{\mathsf c}

\def\PM{\mathscr{P}}
\def\Meas{\mathscr{M}}

\newcommand{\weaks}{\overset{*}{\rightharpoonup}}

\DeclareMathOperator{\Ker}{Ker}
\DeclareMathOperator{\Ran}{Ran}

\newcommand{\ip}[1]{\langle {#1}\rangle}
 
\newcommand{\INDIC}{1\!\!1}

\makeatletter
\def\rightharpoondownfill@{\arrowfill@\relbar\relbar\rightharpoondown}
\def\leftharpoonupfill@{\arrowfill@\leftharpoonup\relbar\relbar}

\newcommand{\xleftrightharpoons}[2][]{\mathrel{%
  \raise.22ex\hbox{$\ext@arrow 3095\leftharpoonupfill@{\phantom{#1}}{#2}$}%
  \setbox0=\hbox{$\ext@arrow 0359\rightharpoondownfill@{#1}{\phantom{#2}}$}%
  \kern-\wd0
  \lower.22ex\box0}}
\newcommand{\xrightleftharpoons}[2][]{\mathrel{%
  \raise.22ex\hbox{$\ext@arrow 3095\rightharpoonupfill@{\phantom{#1}}{#2}$}%
  \setbox0=\hbox{$\ext@arrow 0359\leftharpoondownfill@{#1}{\phantom{#2}}$}%
  \kern-\wd0 \lower.22ex\box0}}
\makeatother
 \def\REVER#1#2{\xleftrightharpoons[#2]{#1}}

\newcommand{\BBneu}[3]{\!\;\bbB_{#1}^{#2}\!(#3)\!\;}

\usepackage[colorlinks,linkcolor={black},%
            citecolor={black},urlcolor={black}]{hyperref}

\definecolor{jan}{rgb}{0.6,0.0,0.1}

\numberwithin{equation}{section}

\usepackage{tikz}
\usetikzlibrary{matrix}

\newcommand{\biggset}[2]{ \bigg\{\: #1 \; \bigg| \; #2 \:\bigg\} }

\begin{document}
\title{Modeling of chemical reaction systems with\\ 
detailed balance using gradient structures} 

\author{
Jan Maas\thanks{Institute of Science and Technology Austria (IST
  Austria), Am Campus 1, 3400 Klosterneuburg, Austria}\ \  
  and 
Alexander Mielke\thanks{WIAS Berlin, Mohrenstra\ss{}e 39, 10117
    Berlin and Institut f\"ur Mathematik, Humboldt-Universit\"at zu
    Berlin, Unter den Linden 6, 10119 Berlin, Germany} 
}

\date{6 April 2020 
} 
\maketitle
\unitlength1cm

\begin{abstract}
  We consider various modeling levels for spatially homogeneous
  chemical reaction systems, namely the chemical master equation, the
  chemical Langevin dynamics, and the reaction-rate
  equation. Throughout we restrict our study to the case where the
  microscopic system satisfies the detailed-balance condition. The
  latter allows us to enrich the systems with a gradient structure,
  i.e.\ the evolution is given by a gradient-flow equation. We present
  the arising links between the associated gradient structures that
  are driven by the relative entropy of the detailed-balance steady
  state.  The limit of large volumes is studied in the sense of
  evolutionary $\Gamma$-convergence of gradient flows.  Moreover, we
  use the gradient structures to derive hybrid models for coupling
  different modeling levels.
\end{abstract}

\section{Introduction}\label{s:Intro}
 
In this work we discuss different models for chemical reactions taking
place in a container of volume $V$.  Throughout we assume that the
spatial extent of the container and the position of the chemical
species are irrelevant, which means that we are looking at a
well-stirred system.  We assume that the system is composed of $I$
different species named $X_1$ to $X_I$, which may represent different
molecules, e.g., $X_1=\rmH_2$, $X_2=\rmO_2$, and $X_3=\rmH_2\rmO$.  We
assume that these $I$ species undergo $R$ different reactions of
mass-action type:
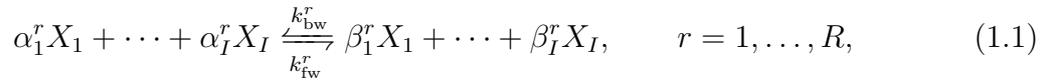
\begin{align}\label{eq:intro-reaction}
 	\alpha^r_1 X_1+\cdots +\alpha^r_I X_I 
		\REVER{k^r_\BW}{k^r_\FW} 
	\beta^r_1 X_1+\cdots +\beta^r_I X_I, 
		\qquad r = 1, \ldots, R,
\end{align}
where the vectors $\bfalpha^r, \bfbeta^r\in \N_0^I$ contain the
stoichiometric coefficients, and $k^r_\FW, k^r_\BW > 0$ are the
forward and backward reaction rates, see Section \ref{s:ReactKin}.
The reaction $2 \rmH_2 + \rmO_2 \REVER{}{} 2 \rmH_2\rmO$ would lead to
the vectors $\bfalpha=(2,1,0)$ and $\bfbeta=(0,0,2)$.

Denoting by $\bfc=(c_1,\ldots, c_I) \in \bfC:={[0,\infty[}^I$ the
vector of nonnegative densities, the simplest model is the macroscopic
\emph{\bfseries reaction-rate equation} (RRE), which is a system of
ODEs on the state space $\bfC$:
\begin{align}
  \label{eq:I.RRE} \tag{RRE}
  \dot \bfc = -\bfR(\bfc) \quad \text{with }\bfR(\bfc):=\sum_{r=1}^R 
  \big(k^r_\FW \bfc^{\bfalpha^r} -
  k^r_\BW\bfc^{\bfbeta^r}\big)\big( \bfalpha^r {-} \bfbeta^r\big).
\end{align}  
Here the monomials $\bfc^{\bfalpha^r}:=\Pi_{i=1}^I c_I^{\alpha^r_i}$
indicate that the probability for the right number of particles for
the $r$th reaction to meet is given by a simple product of the
corresponding densities, i.e., we assume that the positions of the
particles are independent.

A truly microscopic model can be obtained as a stochastic
process. Here we count the number of particles $N_i^V(t)$ for each
species $X_i$ and consider the random vector
$\bfN^V(t)=(N_1^V(t),\ldots,N_I^V(t)) \in \calN:= \N_0^I$.  A forward
or backward reaction of type $r$ is modeled as an instantaneous event
where the particle numbers jump from $\bfN^V(t) + \bfalpha^r$ to
$\bfN^V(t) + \bfbeta^r$ or vice versa.  The corresponding jump rates
in a volume of size $V > 0$ are given by
$k^r_\FW\BBneu V {\bfalpha^r} {\bfN^V(t)}$ and
$k^r_\BW \BBneu V {\bfbeta^r} {\bfN^V(t)}$ respectively; see
\eqref{eq:CME.3} for the definition of $\BBneu V {\bfalpha} {\bfn}$.

Here we study the vector of probabilities
\[
\bfu(t) \in
	\PM(\calN) := \bigset{ \bfv=(v_\bfn)_{\bfn \in \calN}}{v_\bfn \geq 0,\
  \sum\nolimits_{\bfn\in\calN} v_n=1}
\]
that describes the probability distribution of the random variable
$\bfN^V(t)$.  The time evolution of $\bfu(t)$ is given by the
\emph{\bfseries chemical master equation} (CME), i.e., the Kolmogorov forward
equation associated with the continuous time Markov chain above.  This
is a countable linear system of ODEs:
\begin{align}
  \label{eq:I.CME} \tag{CME}
    \dot \bfu(t) = \calB_V \bfu(t), \qquad \bfu(0) = \bfu_0,
\end{align}
where $\calB_V$ is an (unbounded) linear operator on $\ell^1(\calN)$,
see Section \ref{s:CME}, where also existence and uniqueness of
solutions is discussed.  We refer to \cite{MLSH11?HSDS} for a short
introduction to the CME and to \cite{Gill92RDCM} for a justification.

The basis of this work is the observation from \cite{Miel11GSRD} that
\eqref{eq:I.RRE} can be interpreted as a gradient flow if the reaction
system satisfies the detailed-balance condition, i.e., there exists a
positive equilibrium
$\bfc_*=(c_i^*)_{i=1,\ldots,I} \in {]0,\infty[}^I$ such that
\[
	\kappa_*^r 
		:= k^r_\FW \bfc_*^{\bfalpha^r}
		 = k^r_\BW \bfc_*^{\bfbeta^r}
	\ \text{ for } r = 1, \ldots, R. 
\]
Defining the Boltzmann entropy $E$ and the Onsager operator $\bbK$ via 
\begin{align*}
	&E(\bfc) 
	= \sum_{i=1}^I 
		\LB\big(\frac{c_i}{c_i^*} \big)c_i^* 
			\qquad
			\text{with }\LB(z)=z \log z - z +1 \text{ and}\\
&\bbK(\bfc)= \sum_{r=1}^R \kappa_*^r \,\Lambda\big( 
 \frac{\bfc^{\bfalpha^r}}{\bfc_*^{\bfalpha^r}}, 
   \frac{\bfc^{\bfbeta^r}}{\bfc_*^{\bfbeta^r}}\big) \,
 \big(\bfalpha^r{-}\bfbeta^r\big) \oti
 \big(\bfalpha^r{-}\bfbeta^r\big) \ \in \R^{I\ti I}_{\text{sym},\geq 0}
\end{align*}
where $\Lambda(a,b)=\int_0^1 a^s b^{1-s}\dd s$ is the logarithmic
mean, we see that \eqref{eq:I.RRE} is generated by the gradient system
$(\bfC,E,\bbK)$, namely
$\dot\bfc =-\bfR(\bfc)=- \bbK(\bfc)\rmD E(\bfc)$.  In Section
\ref{su:GenerGradStr} we also discuss further gradient structures, 
e.g.\ those used in \cite{MiPeRe14RGFL, MPPR17NETP, MieSte19?CGED}.

If \eqref{eq:I.RRE} satisfies the detailed-balance condition, then
\eqref{eq:I.CME} does so with an equilibrium distribution
$\bfw^V \in \PM(\calN)$ that is explicitly given as a product of
one-dimensional Poisson distributions with mean $c_i^* V$, namely
(cf.\ Theorem \ref{th:DB.CME}),
\[
	w_\bfn^V = 
		\prod_{i=1}^I \ee^{-c_i^*V} 
			\frac{(c_i^*V)^{n_i}}{n_i!}
	\	\text{ for all } 
		\bfn = (n_1,\ldots,n_I) \in \calN.
\]
Consequently, we are also able to interpret \eqref{eq:I.CME} as a
gradient flow induced by a gradient system
$(\PM(\calN), \calE_V,\calK_V)$, see \eqref{eq:CME.E.K}.  Here
$\calE_V(\bfu)$ is again the Boltzmann entropy with respect to
$\bfw^V$, but now divided by the volume $V$:
\begin{equation}
  \label{eq:I.calE.V}
  \calE_V(\bfu)=\frac1V \sum_{n\in \calN}
\LB\big(\frac{u_\bfn}{w^V_\bfn}\big)w^V_\bfn = \frac1V \sum_{n\in
  \calN} u_\bfn \log u_\bfn + \sum_{\bfn \in \calN} u_\bfn \frac{1}V
\log\frac1{w^V_\bfn}. 
\end{equation}

\subsection*{Large-volume approximations using
  gradient structures} 
\label{su:LargeVolAppr}
 
A major challenge in modeling chemical reactions is the question of
understanding the transition from small-volume effects to the
macroscopic behavior in large volumes.  The first breakthrough was
obtained in \cite{Kurt67COSA, Kurt69ETOS, Kurt70SODE, Kurt72RSDM} by
connecting the particle numbers $\bfN(t)\in \calN$ to the
concentrations $\bfc \in \bfC$ and showing that
\begin{equation}
  \label{eq:I.Kurtz}
  	\frac1V \bfN^V(0)\to \bfc_0 
		\text{ almost surely implies }
	\frac1V \bfN^V(t)\to \bfc(t) 
		\text{ almost surely for all } t>0,
\end{equation}
where $t \mapsto \bfc(t)$ is the solution of \eqref{eq:I.RRE} with
$\bfc(0)=\bfc_0$.  This result may be interpreted as a justification
for the RRE in terms of the Markovian model.  In \cite{MiRePe16GORR,
  MPPR17NETP} a dynamic large deviation principle is applied to
$\frac1V \bfN^V(\cdot)$, which leads to a rate functional that
generates a gradient structure $(\bfC,\calE,\Psi_{\cosh})$; see
Section \ref{su:GenerGradStr}. Recent large deviation results for chemical reaction networks can be found in \cite{AgDeEc18OGCR,AgDeEc18LDTM}.

In this paper we study the limit $V\to \infty$ for the gradient system
$(\PM(\calN), \calE_V,\calK_V)$, and hence for \eqref{eq:I.CME}, in
the sense of evolutionary $\Gamma$-convergence for gradient systems,
as introduced in \cite{SanSer04GCGF, Serf11GCGF} and further developed
in \cite{Miel16EGCG, DoFrMi19GSWE}.  For this purpose we use a
suitable embedding $\iota_V: \PM(\calN) \to \PM(\bfC)$ (Section
\ref{s:Liouville}) and obtain the coarse grained gradient system
$(\PM(\bfC),\bfE,\bfK)$ with
\[
\bfE(\varrho)=\int_\bfC E(\bfc)\varrho(\rmd \bfc) \quad \text{and} \quad 
\big(\bfK(\varrho)\xi\big)(\bfc) = - \div_{\!\bfc}\big(\varrho(\bfc)
\bbK(\bfc) \nabla_{\bfc} \xi(\bfc) \big).
\] 
In particular, the coarse grained gradient flow equation is the
\emph{\bfseries Liouville equation}
\begin{align}
 \label{eq:I.Lio}\tag{Lio}
   \dot \varrho(t,\bfc) = \div_\bfc\big(\varrho(t,\bfc) \bfR(\bfc)\big),\qquad 
\varrho_{t = 0} = \varrho_0,
\end{align}
associated with \eqref{eq:I.RRE}; here we used that
$\xi = \rmD_\varrho \bfE = E$ and $\bfR=-\bbK\rmD_\bfc E$.  Thus, in
this scaling a pure transport equation remains, while all diffusion
disappears, as can be seen in the factor $1/V$ before the middle sum
in \eqref{eq:I.calE.V}.  In particular, our result is consistent with
Kurtz' result \eqref{eq:I.Kurtz}: by assuming
$\varrho(0) = \delta_{\bfc_0} \in \PM(\bfC)$ we obtain
$\varrho(t) = \delta_{\bfc(t)}$.  While Kurtz works directly on the
Markovian random variables, we work at the level of their
distributions:
\begin{center}
\begin{tikzpicture}
  \matrix (m) [matrix of math nodes,row sep=3em, 
               column sep=4em,minimum width=2em]
  {
     \substack{\textstyle \bfu \in \PM(\calN)\\[1\jot]\textstyle
       \text{CME: \ } \dot \bfu = - \calK_V(\bfu) \rmD \calE_V(\bfu)\qquad} & 
  \substack{\textstyle  \varrho \in \PM(\bfC)\\[1\jot]\textstyle
    \text{\qquad Liouville: } \partial_t \varrho = - \bfK(\varrho)
    \rmD \bfE(\varrho)}   
\\
\substack{\textstyle \bfN \in \calN \\[1\jot]\textstyle
  \text{Markovian model \quad} }& \substack{\textstyle \bfc \in \bfC
  \\[1\jot]\textstyle 
   \text{\quad RRE: \ } \dot\bfc = - \bbK(\bfc) \rmD E(\bfc).
} \\};
  \path[-stealth]
    (m-2-1) edge node [left] {} (m-1-1)
    (m-1-1)        edge node [below] {\small{here}}
            node [above] {\small{
$\iota_V(\bfu^V) \to \varrho$}} (m-1-2)
    (m-2-1.east|-m-2-2) edge node [below] {\small{Kurtz}}
            node [above] {\small{$\frac1V \bfN^V \to \bfc$}} (m-2-2)
 	(m-2-2) edge node [right] {}    (m-1-2);
\end{tikzpicture}
\end{center}

Our convergence result for the gradient systems
$(\PM(\calN), \calE_V,\calK_V)$ to the limiting gradient system
$(\PM(\bfC),\bfE,\bfK)$ can be seen as a concrete example of the
\emph{EDP convergence} of gradient systems as discussed in
\cite{LMPR17MOGG,DoFrMi19GSWE}.  Another example treating the
convergence of ``Markovian discretizations'' towards a Fokker--Planck
equation is studied in \cite{DisLie15GSMC}; see also
\cite{FatSim16GFAH,EFLS16GFSM,Schl19MLBD} for applications to interacting
particle systems.

In addition to the extreme cases $V$ finite and $V \to \infty$ it is
also important to study the case of intermediate $V$, where
$\frac1V \bfN^V(t)$ already behaves continuously but still shows some
fluctuations of standard deviation $1/\sqrt V$, see
\cite{WinSch17HMCR} for a numerical approach to treat the hierarchy
via a suitable hybrid method.  In \cite{Kurt78SATD} it is shown that
the random vector $t\mapsto \bfX^V(t)\in \bfC$ obtained by solving the
stochastic differential equation
\begin{align}
&\nonumber \rmd \bfX^V(t) = -\bfR(\bfX^V(t))\dd t + 
 \frac1{\sqrt{V}} \Big(\Sigma^\FW(\bfX^V(t))\dd \mathfrak B^\FW(t) +
 \Sigma^\BW(\bfX^V(t))\dd \mathfrak B^\BW(t)\Big) \\
\label{eq:RRE.SDE}
&\text{with independent Brownian vectors }\mathfrak B^\FW(t),\:
\mathfrak B^\BW(t) \in \R^R, \text{ and } \\ 
 &\nonumber \Sigma^\FW(\bfX)=
 \Big( \big(\tfrac{\kappa^r\bfX^{\bfalpha^r}}{\bfc_*^{\bfalpha^r}} \big)^{1/2}
   (\bfalpha^r{-}\bfbeta^r) \Big)_r, \; \ \
 \Sigma^\BW(\bfX)=
 \Big(\big( \tfrac{\kappa^r\bfX^{\bfbeta^r}}{\bfc_*^{\bfbeta^r}}\big)^{1/2}
   (\bfbeta^r{-}\bfalpha^r) \Big)_r \in \R^{I\ti R},
\end{align}
(see \cite[Eqn.\,(1.7)]{Kurt78SATD}) yields an improved approximation
because $\frac1V \bfN^V(t)= \bfX^V(t)+ O\big(({\log V})/V\big)$, while
$\frac1V \bfN^V(t)= \bfc(t)+O\big(1/\sqrt V\big)$. This model is a
so-called diffusion approximation, which in the reaction context also
is termed `chemical Langevin dynamics'. In
\cite[Eqn.\,(23)]{Gill00CLE} and \cite[Eqn.\,(7)]{WinSch17HMCR} the
stochastic differential equation \eqref{eq:RRE.SDE} is called
\emph{\bfseries chemical Langevin equation} (CLE). 

The associated Kolmogorov forward equation takes the form
\begin{equation}
  \label{eq:CLE.FP}
  \dot \rho
  	 = 
   \frac1V \!\sum_{i,j=1}^I \!\partial_{ij}^2 \big( \rho\:
   \wh\bbK_{\CLE}(\bfc)_{ij} \big) + \div\big( \rho
  \bfR(\bfc)\big)	 
		\text{   with } 
	\wh\bbK_\CLE =
	 \frac12\big(\Sigma^\FW (\Sigma^\FW)^\tra 
	 		   {+} \Sigma^\BW (\Sigma^\BW)^\tra \big).
\end{equation}
Here the diffusion matrix $\wh\bbK_\CLE$ can be written in  the explicit form 
\begin{equation}
  \label{eq:whK.CLE}
  \wh\bbK_\CLE (\bfc)=\sum_{r=1}^R{\kappa^r}\,
\frac1{2} \Big( \frac{\bfc^{\bfalpha^r}}{\bfc_*^{\bfalpha^r}}  
 {+} \frac{\bfc^{\bfbeta^r}}{\bfc_*^{\bfbeta^r} }
  \Big)\, \big(\bfalpha^r{-}\bfbeta^r\big)\oti \big(\bfalpha^r{-}\bfbeta^r\big)
\end{equation}
that is different from $\bbK(\bfc)$, because in the
former the arithmetic mean while in the latter the logarithmic mean is taken.  

One drawback of the chemical Langevin equation \eqref{eq:CLE.FP} is
that it cannot be written as gradient flow of the relative entropy, as
the Einstein relation for the drift flux and the diffusion flux is not
satisfied.  Therefore we propose other approximations that stay inside
the theory of gradient flows and seem to work sufficiently well if the
concentrations are not too large or small. Our simplest approximation
is given by the gradient system $(\PM(\bfC), \wt \bfE_V, \bfK)$ with
\[
 \wt \bfE_V(\varrho)  = \int_\bfC \Big( \frac1V \rho\log \rho \,+\, \rho E \Big)
\dd \bfc, \ \text{ where } \varrho = \rho \dd \bfc,
\]
which leads to the linear \emph{\bfseries Fokker--Planck equation} 
\begin{equation}
  \label{eq:I.FPE}\tag{FPE}
  \dot \rho = \div \Big( \frac1V \bbK(\bfc) \nabla \rho + \rho \bfR
\Big). 
\end{equation}
In Section \ref{su:FokkerPlanck} we show that by systematically
deriving higher-order corrections to $\wt\bfE_V$ and  $\bfK$
we can recover the asymptotically correct diffusion matrix
$\wh\bbK_\text{CLE}$ while keeping the gradient structure,
but have to accept several additional terms, or switch over
to the notion of \emph{asymptotic gradient flow structures} in the
sense of \cite{BBRW17CDSE}.

\subsection*{Hybrid modeling using gradient structures}

A major advantage of the gradient flow description is that the
different structures can be combined to obtain hybrid models, in which
the set of chemical species is divided into subclasses which may be
treated differently depending on the desired or needed accuracy.  Our
approach is based on the idea of model reduction for gradient
structures.  The idea is to approximate a complicated gradient
structure $(\sfX,\sfE_X,\sfK_X)$ by a simpler one
$(\sfY,\sfE_Y,\sfK_Y)$ via an embedding mapping $\sfx = \Phi(\sfy)$.
Staying within the class of gradient systems has the advantage that
the most important features of the original system can be preserved.
In particular, decay of the driving functional along the approximate
flow holds automatically.  By contrast, such crucial features could
get lost in a direct approach based on the evolution equation itself.

In Section \ref{s:Hybrid} we shall deal with three examples for hybrid
models where it is essential to keep $V$ as a large but finite
parameter.  First, we shall consider a hybrid model in which an RRE is
coupled to a Fokker--Planck equation.  Here the set of species is
divided into two classes: $\bfC = \bfC_\rms\ti\bfC_\rmm$.  Some of
them will be described \emph{stochastically} (s), while others are
described \emph{macroscopically} (m).  This leads to a gradient flow
structure on the hybrid state space
$\sfY=\PM(\bfC_\rms)\ti \bfC_\rmm$.  The resulting gradient flow
equation turns out to be a mean-field equation, in which the density
of the component $\bfc_\rms$ satisfies a linear equation which is
nonlinearly coupled to an ODE for the component $\bfc_\rmm$.

We also study the coupling of an RRE for macroscopic variables to a
CME for $n$ microscopic variables.  This leads to a hybrid system on
$\PM(\N_0^n)\ti {\bfC_\rmm}$. Finally we analyze a mixed CME /
Fokker--Planck model with state space $\PM(\mfN)$, in which the
underlying space $\mfN:=\{0,1,\ldots,N{-}1\}\cup {[N/V,\infty[}$
contains a mixture of discrete and continuous components.

The present work concentrates solely on the analytical underpinnings
of hybrid modeling for CME; for numerical approaches to CME and to
spatio-temporal CME we refer to \cite{ACTDH05ASHS, HelLot07HMCM,
  MunKha07MTIF, High08MSCR, Engb09SASC, Jahn11RMCM, DolKho14SSTA,
  WinSch17HMCR}.  \bigskip

\noindent{\itshape\bfseries Notational conventions.}
Throughout the paper we will consistently use the following notation
to distinguish the different modeling levels.\medskip

\noindent
\textbf{Reaction-rate equation:} The RRE is denoted by $(\bfC,E,\bbK)$:\\ 
state and state space $\bfc \in \bfC:={[0,\infty[}^I$, 
steady state $\bfc^*=\bfc_*$,
dual variable $\bfzeta$ \\
energy functional $E(\bfc)$, 
Onsager operator $\bbK(\bfc)$\\ 
conserved quantities $\bbQ \bfc = \bfq$, 
stoichiometric subsets
$\bfI(\bfq)=\set{\bfc\in \bfC}{ \bbQ\bfc = \bfq}$.\medskip

\noindent
\textbf{Chemical master equation:} The CME is denoted by 
$(\PM(\calN),\calE_V,\calK_V)$:\\  
state and state space $\bfu=(u_\bfn)_{\bfn\in \calN} \in
 \PM(\calN)\subset$, 
steady state  $\bfw^V$,
dual variable $\bfmu$\\
energy functional $\calE_V(\bfu)$, Onsager operator $\calK_V(\bfu)$
\\
invariant subsets $\calI(\ol\bfn)=\set{\bfn\in \calN}{ \bbQ\bfn =
  \bbQ\ol\bfn}$.\medskip 

\noindent\textbf{Liouville equation:} The LE is denoted by 
$(\PM(\bfC),\bfE,\bfK)$:\\ 
state and state space $\varrho = \rho \dd \bfc \in \PM(\bfC)$, 
steady state $\delta_{\bfc_*}$,
dual variable $\xi$ \\
energy functional $\bfE(\varrho) = \int_\bfC E(\bfc) \dd
\varrho(\bfc) $,  Onsager 
operator $\bfK(\varrho)=-\div\!\big(
\varrho(\cdot) \bbK(\cdot)\nabla \Box\big)$\medskip

\noindent
\textbf{Fokker--Planck equation:} The FPE is denoted by
$(\PM(\bfC), \wt\bfE_V ,\bfK)$:\\
state and state space $\varrho \in \calP(\bfC):=\PM(\bfC)$, steady
state $W_V$, dual variable $\xi$\\
energy functional $\wt\bfE_V(\varrho) =
\int_\bfC \big(\frac1V \rho(\bfc)\log \rho(\bfc)
{+}\rho(\bfc) E(\bfc) \big)\dd \bfc $, Onsager operator
$\bfK$.\medskip

\noindent
\textbf{Hybrid systems} are denoted by ``mathfrak'' letters:\\
$(\PM(\bfC_\rms)\ti \bfC_\rmm, \mfE^\text{FP-RR}_V,
\mfK^\text{FP-RR}_V)$ for coupling FPE and RRE\\
$(\PM(\N_0^{J})\ti \bfC_\rmm, , \mfE^\text{CM-RR}_V,
\mfK^\text{CM-RR}_V)$ for coupling CME and RRE\\
 $(\PM(\mfN_{V,N}), \mfE_{V,N},\mfK_{V,N})$ for merging discrete and
 continuous modeling for one species.
\medskip

\noindent
The space of all signed Borel measures of bounded variation on $\bfC$
is denoted by $\Meas(\bfC)$.

\section{Reaction rate equations}
\label{s:ReactKin}

We denote by
$\bfc=(c_1,\ldots,c_I)\in \bfC:= \left[ 0, \infty\right[^I$ the
concentrations of $I$ different chemical species $X_1,\ldots, X_I$
reacting according to the mass action law, i.e., the reactions
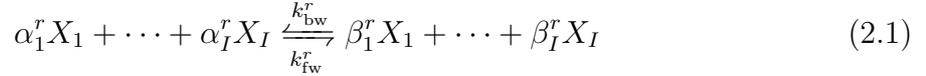
\begin{equation}
  \label{eq:ReactNetwork}
  \alpha^r_1 X_1+\cdots +\alpha^r_I X_I \REVER{k^r_\BW}{k^r_\FW} 
\beta^r_1 X_1+\cdots +\beta^r_I X_I 
\end{equation}
for $r=1,\ldots,R$, where $R$ is the number of possible reactions,
$\bfalpha^r,\bfbeta^r\in \N_0^I$ are the vectors of the stoichiometric
coefficients, and $k^r_\FW, k^r_\BW > 0$ are the forward and backward
reaction-rates. In general these rates may depend on $\bfc$, but for
simplicity we keep them as constants in this work.  A typical example
is the splitting of water into hydrogen and oxygen, namely
$2\rmH_2 + \rmO_2 \; \REVER{}{} \; 2 \rmH_2\rmO$.
  
The corresponding \emph{reaction-rate equations} (RRE) are given via
the ODE system
\begin{equation}
  \label{eq:ReactKin}
  \dot \bfc = -\bfR(\bfc) \quad \text{with }\bfR(\bfc):=\sum_{r=1}^R 
  \big(k^r_\FW \bfc^{\bfalpha^r} -
  k^r_\BW\bfc^{\bfbeta^r}\big)\,\big( \bfalpha^r {-} \bfbeta^r\big),
\end{equation}
where $\bfc^\bfalpha=c_1^{\alpha_1}\cdots c_I^{\alpha_I}$, see
\cite{FeiHor77CMSC,Grog83ABSC,ErdTot89MMCR}.

\subsection{Stoichiometry, conservation, and decomposition of the
  state space}
\label{su:StoichCons}

The stoichiometric subspace $\bbS \subset \R^I$ and its orthogonal complement
$\bbS^\perp$ are defined via
\begin{equation}\label{e:Stoch}
  \bbS:=\mathrm{span}\set{\bfalpha^r - \bfbeta^r}{r=1,\ldots,R}, \ \ \
  \bbS^\perp:=\set{\bfxi\in \R^I}{ \bfxi{\cdot}\bfmu=0\text{ for all
    }\bfmu\in \bbS}.
\end{equation}
For each $\bfxi\in \bbS^\perp$ the function $C_\bfxi(\bfc)=\bfxi\cdot
\bfc$ defines a first integral, which easily follows from $\bfxi \cdot
\bfR(\bfc)\equiv 0$.  These conservation laws often go under the name
\emph{conservation of atomic species}, see \cite{ErdTot89MMCR}.
Suppose now that $\bbS^\perp$ is a non-trivial subspace of $\R^I$. We
shall argue that the RRE induces a decomposition of the state space
$\bfC={[0,\infty[}^I$ into affine invariant subsets. (If $\bbS^\perp =
\{ {\bf 0}\}$, the only invariant set is $\bfC$ itself.)

Choosing a basis $\set{ \bfm_k \in \R^I}{ k=1, \ldots ,m_\bbW}$ of
$\bbS^\bot$ we define the matrix $\bbQ\in \R^{m_\bbW\ti I}$, which has
the rows $\bfm_k \in \R^I$. By construction we have
$\bbQ[\bbS]=\{\bf0\}$, and we conclude that the solutions $\bfc$ of
\eqref{eq:ReactKin} conserve $\bbQ\bfc$ as follows:
\begin{equation}
  \label{eq:RR.cons}
  \dot \bfc= - \bfR(\bfc)  \quad \Longrightarrow \quad \bbQ \bfc(t) =
\bbQ \bfc(0) \text{ for }t>0.
\end{equation}
By construction every affine conserved quantity is of the form
$\bfxi \cdot \bfc + q$ for some $\bfxi \in \bbS^\perp$ and $q \in \R$.
This allows us to decompose the full state space $\bfC={[0,\infty[}^I$
into the invariant, affine subsets $(\bfc_0{+}\bbS)\cap \bfC$ for
$\bfc_0 \in \bfC$. Using the notation
\[
\mfQ:= \set{\bbQ\bfc \in \R^{m_\bbW} }{ \bfc\in \bfC}
\]
we define, for all $\bfq\in \mfQ$,  the sets 
\begin{align}
\label{eq:I(bfq)}
\bfI(\bfq):= \set{\bfc \in \bfC}{ \bbQ\bfc=\bfq}. 
\end{align}
Then, $\bfq_1\neq \bfq_2$ implies $\bfI(\bfq_1)\cap \bfI(\bfq_2)=\emptyset$,
and we have $\bfC=\bigcup_{\bfq\in \mfQ} \bfI(\bfq)$.  Let us note that
this decomposition does not depend on the choice of the orthonormal
basis which determines the matrix $\bbQ$, although the set $\bfI(\bfq)$
does depend on $\bbQ$. Note also that we can always write $\bfI(\bfq) =
(\bfc{+}\bbS)\cap \bfC$ for some arbitrary $\bfc \in \bfC$ satisfying
$\bbQ \bfc = \bfq$.

\subsection{Detailed balance and the Wegscheider matrix}
\label{su:WegDBC}

We say that the above reaction system fulfills the \emph{condition of detailed balance} if there exists a positive equilibrium density
vector $\bfc_* \in \left]0,\infty\right[^I $ such that all reactions are simultaneously in equilibrium, i.e., 
\begin{equation}\label{RRE:DBC}
\kappa^r_*:=k^r_\FW\bfc_*^{\bfalpha^r} = k_\BW^r \bfc_*^{\bfbeta^r} 
  \ \text{ for } \  r=1,\dots,R.
\end{equation}
This condition implies that $\bfR(\bfc_*)=0$, but we emphasize that this
condition is stronger in general cases. 
The condition of detailed balance is also called the condition of
\emph{microscopic reversibility}, see \cite[p.\,45]{ErdTot89MMCR} or 
\cite{DegMaz84NET} for a general discussion of these concepts.

We are looking for a characterization of detailed balance.
Let $\bbW\in \Z^{R\times I} $ be the matrix which has the row vectors $ \bfgamma^r:=\bfalpha^r {-} {\bfbeta^r} \in \Z^I$, $r=1, \dots, R$. 
We call $\bbW$ the \emph{Wegscheider matrix} because of the pioneering
work in \cite{Wegs02SGBT}. We then have
\[
	\bbS = \Ran \bbW^\tra\quad
	 \text{and} \quad 
	\bbS^\bot = \Ker\bbW,
\]
which explains the abbreviation $m_\bbW:= \dim \bbS^\bot=\dim
\Ker\bbW$. Since $\bfc_*$ is strictly positive, we can take the
logarithm of the polynomial conditions \eqref{RRE:DBC} and find the
equivalent linear system
\begin{equation}\label{bb2}
  \bbW \,\bflog \bfc_* = \big(\log  (k_\BW^r/k^r_\FW)\big)_{r=1,\dots,R}, \quad 
\text{where } \bflog \bfc = \big( \log c_i\big)_{i=1,\ldots,I}. 
\end{equation} 
By Fredholm's alternative, \eqref{bb2} is solvable if and only if
\begin{equation}\label{bb3}
  \bfy\cdot\big(\log (k_\BW^r/k^r_\FW)\big)_{r=1,\dots,R} =0\quad\text{for all }\bfy\in 
  \Ker\bbW^\tra.
\end{equation}
These conditions on the reaction coefficients $k^r_\FW$ and $k^r_\BW$
are called \emph{Wegscheider conditions} (see, e.g., \cite{Wegs02SGBT,SchSch89GWCI,VlaRos09TBCR,GliMie13GSSC}). By
choosing a basis of $\Ker \bbW^\tra$ and exponentiation
they can be rewritten as polynomial conditions without referring to the
equilibrium state $\bfc_*$. 

Let $n_\bbW := \dim(\Ker \bbW^\tra) \in \N_0 $ denote the
number of Wegscheider conditions. Then the following assertions hold:
\begin{itemize}\itemsep-0.3em
\item[(i)] If the stoichiometric vectors $\bfalpha^r-\bfbeta^r$, $r=1,\dots,R$, are linearly independent, then $\Ker \bbW^\tra = \{ {\bf 0} \}$, hence there is no Wegscheider condition.
\item[(ii)] If  $\bfalpha^r-\bfbeta^r$, $r=1,\dots,R$, are
linearly dependent, then $\dim(\Ker \bbW^\tra)>0$ and
non-trivial Wegscheider conditions appear.
\end{itemize}
Since $\dim(\Ran \bbW) = \dim(\Ran \bbW^\tra) = \dim \bbS$ by standard
linear algebra, the number of Wegscheider conditions can be expressed
as
\begin{align*}
	 n_\bbW = R - \dim \bbS 
 		    = R - I + \dim(\Ker \bbW) 
	    	= R - I  + m_\bbW. 
\end{align*}
Hence, if the number $R$ of reactions is smaller than the number $I$
of species, the Wegscheider conditions can usually be satisfied
easily.

\begin{remark}[Wellposedness of RRE] \label{rem:RRE-well-posed} We
  conclude this subsection with a statement concerning the
  well-posedness of the RRE given as in Theorem \ref{thm:RRE.GradSys}
  below. For all $\bfc(0)\in \bfC={[0,\infty[}^I$ there exists a
  unique global solution $\bfc: {[0,\infty[} \to \bfC$. Local
  existence for solutions starting in the interior of $\bfC$ is
  trivial, as $\bfR$ is a polynomial vector field. Since the relative
  entropy $E$ is a coercive Liapunov functional, the solutions cannot
  blow up and stay inside a region $B_R(0) \cap \bfC$ for some $R>0$.

  Moreover, solutions cannot leave this region via the boundary $\pl
  \bfC$, since the vector field is either tangential to $\pl\bfC$ or
  points inwards. Indeed, if  $c_j(t_0)=0$ for
  some $j$, then 
\[
\dot c_j(t_0)=- R_j(\bfc(t_0)) = -\sum_{r=1}^R \kappa_*^r 
\big(\tfrac{\bfc^{\bfalpha^r}(t_0)}{\bfc_*^{\bfalpha^r}} -
\tfrac{\bfc^{\bfbeta^r}(t_0)}{\bfc^{\bfbeta^r}_*} \big) \big(\alpha^r_j {-}
\beta^r_j\big)\geq 0,
\]
because each term in the sum is nonpositive: If
$\alpha^r_j = \beta^r_j$ or $\min\{\alpha^r_j ,\beta^r_j \} >0$, then the
term is 0.  Thus, we are left with the cases $(\alpha^r_j,\beta^r_j)
\in \{ (n,0),(0,n)\}$ for some positive $n$. In the first case
$c_j(t_0)=0$ implies $\bfc^{\bfalpha^r}(t_0)=0$ and the result
follows, and the second case is similar.
\end{remark}

\subsection{The reaction-rate equations as a gradient system}
\label{su:RRE.GradSys}

We show that a RRE satisfying the detailed-balance condition can be
generated by a gradient system $(\bfC,E,\bbK)$. Here, the
state space $\bfC:= {[0,\infty[}^I$ contains all possible
concentration vectors $\bfc$. The driving functional is the relative
entropy $E:\bfC \to {[0,\infty[}$ and the Onsager matrix $\bbK$ is
 chosen suitably (recall that $\LB(z)=z \log z - z +1 \geq 0$):
\begin{equation}
  \label{eq:RRE.E.K}
  E(\bfc) := \sum_{i=1}^I \LB\big(\frac{c_i}{c^*_i}\big) c^*_i
 \  \text{  and } \ 
\bbK(\bfc)= \sum_{r=1}^R \kappa^r_* 
\Lambda\big(\frac{\bfc^{\bfalpha^r}}{\bfc_*^{\bfalpha^r}},
\frac{\bfc^{\bfbeta^r}}{\bfc^{\bfbeta^r}_*}\big)\, 
\big(\bfalpha^r {-} \bfbeta^r\big) \otimes 
   \big(\bfalpha^r {-} \bfbeta^r\big),  
\end{equation}
where the logarithmic-mean
function $\Lambda$ is given via
\begin{equation}
  \label{eq:Lambda}
\Lambda(a,b)= \int_0^1 a^s b^{1-s} \dd s = \frac{a - b}{\log a - \log b}. 
\end{equation}
The following result shows that a RRE \eqref{eq:ReactKin} satisfying
the detailed-balance condition \eqref{RRE:DBC} is indeed generated by
the gradient system $(\bfC,E,\bbK)$. This was first established in
\cite[Sect.\,VII]{Yong08ICPD} to derive entropy bounds for hyperbolic
conservation laws in reactive flows and was rederived in
\cite{Miel11GSRD} in the context of reaction diffusion systems
including electric charge-interactions. It is interesting to note that
for continuous time Markov chains (CTMC), which form a special
subclass of RRE with linear reactions, there are several distinct
gradient structures, see \cite[Prop.\,4.2]{Maas11GFEF} and
\cite[Thm.\,3.1]{Miel13GCRE} and Section \ref{su:CTMCgrad}.
However, in the case of nonlinear reactions according to the mass-action
law, only the gradient structure with the Boltzmann entropy remains. The
key fact is the logarithm identity $(\bfalpha - \bfbeta) \cdot \bflog \bfc =
\log(\bfc^{\bfalpha - \bfbeta})$.

\begin{theorem}[Gradient structure for RRE]\label{thm:RRE.GradSys}
  If the RRE \eqref{eq:ReactKin}
  satisfies the detailed-balance condition \eqref{RRE:DBC} for a
  positive steady state $\bfc_*=(c^*_i)_{i=1,\ldots,I}$, then it has the
  gradient structure $(\bfC,E,\bbK)$ defined in
  \eqref{eq:RRE.E.K},  namely $\dot\bfc= -
  \bfR(\bfc)=- \bbK(\bfc) \rmD E(\bfc)$. 
\end{theorem}
\begin{proof} Multiplying
$\rmD E(\bfc)= (\log(c_i/c^*_i))_{i=1,\ldots,I}$ 
by $\bfalpha^r{-}\bfbeta^r \in \R^I$ we obtain 
\begin{equation}\begin{aligned}  \label{eq:RRE4}
  (\log(c_i/c^*_i))_{i=1,\ldots,I} \cdot \big(\bfalpha^r {-}
\bfbeta^r\big)
&= \sum_{i=1}^I \Big( \alpha_i^r \log(c_i/c^*_i) - 
\beta_i^r \log(c_i/c^*_i)\Big) 
\\& =  
 \log\big(\tfrac{\bfc^{\bfalpha^r}}{\bfc_*^{\bfalpha^r}}\big) - 
\log \big( \tfrac{\bfc^{\bfbeta^r}}{\bfc^{\bfbeta^r}_*} \big) ,
\end{aligned}\end{equation}
which is the denominator of
$\Lambda\big(\tfrac{\bfc^{\bfalpha^r}}{\bfc_*^{\bfalpha^r}},
\tfrac{\bfc^{\bfbeta^r}}{\bfc^{\bfbeta^r}_*}\big)$. Hence, using
$\Lambda(a,b)(\log a {-}\log b) = a{-}b$ gives
\[
\bbK (\bfc)\rmD E(\bfc)=\sum_{r=1}^R \kappa_*^r 
\big(\tfrac{\bfc^{\bfalpha^r}}{\bfc_*^{\bfalpha^r}} -
\tfrac{\bfc^{\bfbeta^r}}{\bfc^{\bfbeta^r}_*} \big) \big(\bfalpha^r {-}
\bfbeta^r\big) \overset{\rmD\rmB}=\sum_{r=1}^R   \big(k^r_\FW
\bfc^{\bfalpha^r}   {-} k^r_\BW\bfc^{\bfbeta^r}\big)\, 
\big( \bfalpha^r {-} \bfbeta^r\big) =  \bfR(\bfc),
\]
where we used the detailed-balance condition \eqref{RRE:DBC} in $\overset{\rmD\rmB}=$. Thus, the assertion is established.
\end{proof}

Summarizing the above derivations, we have rewritten the RRE in
thermodynamic form
\begin{equation}
  \label{eq:Onsager}
\dot \bfc = - \bfR(\bfc)= - \bbK(\bfc) \bfmu \quad \text{with }\ 
\bfmu = \rmD E(\bfc),
\end{equation}
which is also called the Onsager principle
\cite{Onsa31RRIP,OnsMac53FIP}. The latter states that the rate (flux)
of a macroscopic variable is given as the product of a symmetric
positive definite matrix $\bbK$ and the thermodynamic driving force
$-\bfmu$, see e.g.\ \cite[Ch.\,X, \S\,4]{DegMaz84NET}.  The symmetry
$\bbK = \bbK^\top$ is related to microscopic reversibility,
i.e.,~detailed balance, see also \cite{MiPeRe14RGFL, MiRePe16GORR}.
Subsequently, we refer to $\bbK$ as the Onsager operator or matrix.

Here we clearly see the advantage of using the Onsager operator $\bbK$
to write the RRE as a gradient system,as opposed to working with the
Riemannian tensor: we do not have to take care of the fact that $\bbK$
is not invertible except if $\bbS=\R^I$.

\subsection{Continuous time Markov chains as a gradient system}
\label{su:CTMCgrad}

The forward equation for a reversible CTMC on
a discrete space $\{1,2,\ldots,I\}$ is a special case of the RRE
considered above. In this case all reactions are of the form
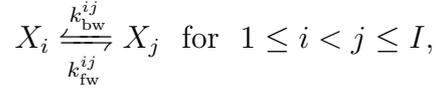
\begin{align*}
 X_i \REVER{k^{ij}_\BW}{k^{ij}_\FW} 
 X_j   \ \text{ for } \  1 \leq i < j \leq I,
\end{align*}
and the reaction rates $k^{ij}_\FW$ (resp. $k^{ij}_\BW$) are
interpreted as the transition rates from $i$ to $j$ (resp. from $j$ to
$i$).  The reaction-rate equation is the linear system of ODEs
\begin{align}\label{eq:Markov}
\dot\bfc = - \bfR(\bfc) = \calA \bfc \quad \text{with } \calA \bfc = -\sum_{i<j}
 \big(k^{ij}_\FW c_i - k^{ij}_\BW c_j\big)
 (\bfe_i{-}\bfe_j),
\end{align}
and the detailed-balance condition for the equilibrium state
$\bfc_*$ takes the form 
\begin{align}\label{eq:DB-Markov}
 \kappa_*^{ij} := c_i^* k^{ij}_\FW = c_j^* k^{ij}_\BW
  \ \text{ for } \  1 \leq i < j \leq I.
\end{align}
Using this condition, the RRE can be written
coordinate-wise as
\begin{align*}
 \frac{\dot c_i}{c_i^*} = \sum_{j < i}
k^{ji}_\BW \Big(\frac{c_j}{c_j^*} - \frac{c_i}{c_i^*}\Big)
 + \sum_{j > i}
k^{ij}_\FW \Big(\frac{c_j}{c_j^*} - \frac{c_i}{c_i^*}\Big),
\ \text{ or equiv., } \
 {\dot c_i} = \sum_{j \neq i}
\kappa_*^{ij} \Big(\frac{c_j}{c_j^*} - \frac{c_i}{c_i^*}\Big).
\end{align*}
Here we used the notational convention that $\kappa_*^{ij} :=
\kappa_*^{ji}$ for $j < i$.  The relative entropy $E$ is as above and
the Onsager matrix takes the form
\begin{align}\label{eq:KM}
\bbK_\rmM (\bfc)= \sum_{i<j}
\kappa_*^{ij} \Lambda\Big(\frac{c_i}{c_i^*} , \frac{c_j}{c_j^*}\Big)
(\bfe_i{-}\bfe_j) \otimes (\bfe_i{-}\bfe_j),
\end{align}
where $\bfe_i\in \R^I$ denotes the $i$-th unit vector.
We then have the gradient structure $(\bfC,E, \bbK_{\rm M})$, namely
\begin{align*}
\dot\bfc = \calA \bfc =- \bbK_{\rm M}(\bfc) \rmD E(\bfc).
\end{align*}

This gradient flow structure has been found in the independent works
\cite{Maas11GFEF} (which deals with Markov chains exclusively) and
\cite{Miel11GSRD} (in the setting of reaction-diffusion systems, in
which Markov chains are implicitly contained). The related work
\cite{CHLZ12FPEF} deals with discretizations of Fokker--Plank equations.

In fact, for the construction of gradient structures for Markov chains
$\dot \bfc = \calA \bfc$ we do not need the summation rule for
logarithms. Hence, following \cite{Maas11GFEF,Miel13GCRE} there are
more general gradient structures.  Choosing a strictly convex
function $\phi :[0,\infty[ \to \R$ that is smooth on $]0,\infty[$ we set 
\begin{equation}
  \label{eq:EK.general}
	E^\phi(\bfc) 
		:= \sum_{i=1}^n c_{i}^* 
			\,\phi\Big( \frac{c_i}{c_i^*}\Big), 
			\quad 
	\bbK_{\rm M}^\phi(\bfc) 
		= \sum_{j=2}^n \sum_{i=1}^{j-1} \kappa_*^{ij} 
			\:\Phi \Big( \frac{c_i}{c_i^*}\,,\, 
						 \frac{c_j}{c_j^*}\Big) \, 
			 (\bfe_i{-}\bfe_j) \otimes (\bfe_i{-}\bfe_j), 
\end{equation}
where $\Phi( a,b)=(a-b)/(\phi'(a){-}\phi'(b))$ for $0<a\neq b$ and
$\Phi(a,a)=1/\phi''(a)$. 
The gradient flow structure $(\bfC,E, \bbK_{\rm M})$ corresponds to the case where $\phi = \lambda_{\rm B}: z \mapsto z \log z - z +1$.

\begin{proposition}[Gradient structure for CTMC]\label{prop:CTMC.GradSys}
  If the CTMC \eqref{eq:Markov} satisfies the detailed-balance
  condition \eqref{eq:DB-Markov} for a positive steady state
  $\bfc_*=(c^*_i)_{i=1,\ldots,I}$, then it has the gradient structures
  $(\bfC,E^\phi, K_\rmM^\phi)$, namely $\dot\bfc= \calA \bfc =-
  K_\rmM^\phi(\bfc) \rmD E^\phi(\bfc)$.
\end{proposition}

\begin{remark}\label{rem:RRE-no-extension}
  The construction in Proposition \ref{prop:CTMC.GradSys} does not
  extend to general RRE. There one would need to replace
  the quantity
  $\Lambda\big(\tfrac{\bfc^{\bfalpha^r}}{\bfc_*^{\bfalpha^r}},
  \tfrac{\bfc^{\bfbeta^r}}{\bfc^{\bfbeta^r}_*}\big)$ in
  \eqref{eq:RRE.E.K} by
$\big(\tfrac{\bfc^{\bfalpha^r}}{\bfc_*^{\bfalpha^r}} -
\tfrac{\bfc^{\bfbeta^r}}{\bfc^{\bfbeta^r}_*}\big)/\big(
  (\bfalpha^r - \bfbeta^r) \cdot \bfphi'(\frac{\bfc}{\bfc_{*}})\big)$,
    but this quantity can be negative in general. As a consequence,
    the corresponding Onsager matrix would not be positive definite.
    This cannot happen for Markov chains (i.e., when $\bfalpha =
    \bfe_i$ and $\bfbeta = \bfe_j$), by virtue of the convexity of
    $\phi$.
\end{remark}

In the following we will mainly concentrate on the gradient structure
$(\bfC,E,\bbK_\rmM)$ with the logarithmic entropy, as it is the only
one that connects with the RRE.

\subsection{Generalized gradient structures} 
\label{su:GenerGradStr}

For Markov chains and RRE there are several families of \emph{generalized
gradient structures} $(\bfC,E,\Psi^*)$ where the quadratic 
function $\Psi^*(\bfc,\bfzeta) = \frac12\langle \bfzeta,
\bbK(\bfc)\bfzeta\rangle$ is replaced by a general \emph{dual
  dissipation potential}  $\Psi^*(\bfc,\,\cdot\,): \R^I \to \left[ 0, \infty
\right[$ that is continuous and
convex and satisfies $\Psi^*(\bfc,0)=0$.

In the case of RRE, the monomial terms
$\bfc^\bfalpha$ can only be generated by the logarithmic summation
rule $\sum_{i=1}^I\log(b_i) = \log\big( \Pi_{i=1}^I b_i\big)$. Hence,
we stick to the relative entropy $E$ defined in
\eqref{eq:RRE.E.K}, i.e.,\ $\phi(z) = \LB(z)$. However, we may replace the
linear Onsager principle 
$\dot \bfc = - \bbK(\bfc) \rmD E(\bfc)$ by the more general nonlinear form
$\dot\bfc= \pl_\bfzeta \Psi^*\big(\bfc, -\rmD E(\bfc)\big)$. 

To define $\Psi^*$ we choose an arbitrary family
of smooth dissipation functionals $\psi_r:\R \to {[0,\infty[}$, i.e.,\
$\psi_r(0)=\psi'_r(0)=0$ and $\psi''_r >0$ and define the
dissipation potential 
\begin{equation}
  \label{eq:GenGS1}
  \Psi^*(\bfc,\bfzeta) = \sum_{r=1}^R  L_r(\bfc)\, \psi_r\big(
 (\bfalpha^r {-} \bfbeta^r)\cdot \bfzeta\big) \  
   \text{ with } L_r (\bfc)= \kappa^r_*\,\frac{ 
   \tfrac{\bfc^{\bfbeta^r}}{\bfc^{\bfbeta^r}_*}
   -\frac{\bfc^{\bfalpha^r}}{\bfc_*^{\bfalpha^r}} }{\psi'_r\big( 
    \log \tfrac{\bfc^{\bfbeta^r}}{\bfc^{\bfbeta^r}_*} 
     - \log\frac{\bfc^{\bfalpha^r}}{\bfc_*^{\bfalpha^r}} \big)}.
\end{equation}
Using \eqref{eq:RRE4} we easily obtain
$-\bfR(\bfc) = \pl_\bfzeta \Psi^*(\bfc,- \rmD E(\bfc))$, i.e.,\
$\dot\bfc = -\bfR(\bfc)$ is generated by the generalized gradient
system $(\bfC,E,\Psi^*)$.

The case $\psi_r(\zeta)=\frac12\zeta^2$ leads to the quadratic
dissipation potential in \eqref{eq:RRE.E.K}, i.e.,\ the functions $L_r$
are given in terms of the logarithmic mean.  In \cite{AlGaHu02TDEM}
the choices $\psi_r(\pm\zeta) = \ee^\zeta-1-\zeta$ is used. 
Based on a derivation via the large deviation principle (see
\cite{MiPeRe14RGFL, MiRePe16GORR, MPPR17NETP}) a special role is
played by the choice of a ``cosh-type'' function $\psi_r$: 
\begin{equation}
  \label{eq:CoshGS}
  \psi_r(\zeta)=\sfC^*(\zeta):=4 \cosh\big(\frac12\zeta\big)- 4 \quad
  \text{giving} \quad 
L_r(\bfc) =   \kappa^r_* 
		\Big(\frac{\bfc^{\bfalpha^r}}{\bfc_*^{\bfalpha^r}}\,
		      \frac{\bfc^{\bfbeta^r}}{\bfc_*^{\bfbeta^r}}
		\Big)^{1/2}. 
\end{equation}
Here $\sfC^*$ is normalized such that $\sfC^*(\zeta)=\frac 12\zeta^2 +
O(\zeta^4)$. Hence, the dual dissipation potential takes the form 
\begin{equation}
  \label{eq:Psi*cosh}
	\Psi^*_{\cosh}(\bfc, \bfzeta) 
	:= \sum_{r=1}^R \kappa^r_* 
\,\Big(\frac{\bfc^{\bfalpha^r}}{\bfc_*^{\bfalpha^r}}\, 
   \frac{\bfc^{\bfbeta^r}}{\bfc_*^{\bfbeta^r}} \Big)^{1/2} \:
   \sfC^*\big( (\bfalpha^r{-} \bfbeta^r)\vdot \bfzeta \big). 
\end{equation}
It is shown in \cite[Prop.\,4.1]{MieSte19?CGED} that this generalized
gradient structure is distinguished as the only tilt-invariant
gradient structure for CTMCs.

\section{The chemical master equation}
\label{s:CME}

\subsection{Modeling discrete particle numbers via CME}
\label{su:ModCME}

The chemical master equation (CME) is a CTMC that is defined on the
set $\calN=\N_0^I$ where $\bfn=(n_1,\ldots,n_I) \in \calN$ is the
vector of particle numbers, see \cite{MLSH11?HSDS} for an
introduction. This means that $n_i\in \N_0$ denotes the 
number of particles of species $X_i$ in a sufficiently big volume,
whose size is denoted by $V>0$.
The modeling assumes that all particles move randomly in this big volume (well-stirred tank reactor) so that they can meet independently.
The dynamics is formulated in terms of the probabilities
\[
u_\bfn(t) = \text{probability that at time $t$ there are $n_i$
  particles of species $X_i$ for $i=1,\ldots,I$.}
\]

All the $R$ reaction pairs may happen independently of each other
according to the number of the available atoms needed for the
reactions and the reaction coefficients $k^r_\FW \geq 0$ and
$k^r_\BW\geq 0$, respectively. Moreover, the jump intensities
\begin{align*}
  k^r_\FW \BBneu V{\bfalpha^r}{\bfn} 
  \text{ from $\bfn + \bfalpha^r$ to $\bfn + \bfbeta^r$} 
  \quad \text{ and } \quad
  k^r_\BW \BBneu V{\bfbeta^r}{\bfn}
  \text{ from $\bfn + \bfbeta^r$ to $\bfn + \bfalpha^r$} 
\end{align*}
also depend on the volume $V$, as $n_i$ denotes the absolute particle
number, while for the reaction the densities $c_i=n_i/V$ matter.  The
specific form of $\BBneu V\bfalpha \bfn$ (cf.\ \cite{Kurt70SODE,
  MLSH11?HSDS}) reads
\begin{equation}
  \label{eq:CME.3}
  \BBneu V \bfalpha {\bfn} 
  	= \left\{ \ba{cl} \ds
  \frac{V (\bfn +\bfalpha)!}{V^{|\bfalpha|}\bfn!} 
  & \text{for } \bfn	 \in \calN, \\ 
0 & \text{for } \bfn \not\in \calN, 
	\ea\right.  \quad
  \text{where } \bfn! = \prod_{i=1}^I n_i !\:. 
\end{equation}
To avoid clumsy notation we defined $\BBneu V \bfalpha \bfn$ for all $\bfn \in \Z^I$, but $\BBneu V \bfalpha \bfn = 0$ if $\bfn\not\in \calN$. 
We also see that $\BBneu V \bfalpha \bfn\approx V \bfc^\bfalpha$ for
$\bfc = \frac{1}{V} \bfn$, where the factor $V$ indicates that the number of
reactions is proportional to the volume of the container, if the
densities are kept constant. 

The CME associated with the RRE \eqref{eq:ReactKin} is the
Kolmogorov forward equation for the probability distributions
$\bfu=(u_\bfn)_{\bfn \in \calN} \in \PM(\calN)$, namely 
\begin{equation}
  \label{eq:CME.1}
 \begin{aligned}
 \dot \bfu = 
  \sum_{r=1}^R \ol\calB_V^r\bfu \ \text{with }
\big(\ol\calB_V^r \bfu\big)_\bfn 
	& = k^r_\FW \big(
		\BBneu V{\bfalpha^r}{\bfn{-}\bfbeta^r} 
			  u_{\bfn{+}\bfalpha^r{-}\bfbeta^r} 
	  - \BBneu V{\bfalpha^r}{\bfn{-}\bfalpha^r} 
	  		  u_\bfn \big)   
\\ &\quad + 
		k^r_\BW \big( 
		\BBneu V{\bfbeta^r}{\bfn{-}\bfalpha^r}
			 u_{\bfn{-}\bfalpha^r{+}\bfbeta^r} 
	  - \BBneu V{\bfbeta^r}{\bfn{-}\bfbeta^r}  
	  		 u_\bfn\big).
\end{aligned}
\end{equation}
The $r$th forward reaction from $\bfn{+}\bfalpha^r$ to
$\bfn{+}\bfbeta^r$ can only happen 
(i.e.,\ $\BBneu V{\bfalpha^r}{\bfn}>0$) 
if $\bfn \geq \bm0$. Hence any
occurring $u_\bfm$ with $\bfm\not\in \calN$ is multiplied by intensity
$0$, so in \eqref{eq:CME.1} we may set $u_\bfm \equiv 0$ for all
$\bfm\not\in \calN$.  The operators $\ol\calB_V^r$ are the adjoints of
the Markov generators $\calQ^r_V$ given by
\begin{align}\label{eq:calL}
(\calQ_V^r \bfmu)_\bfn 
=  k^r_\FW \BBneu V{\bfalpha^r}{\bfn{-}\bfalpha^r} 
	(\mu_{\bfn-\bfalpha^r{+}\bfbeta^r} - \mu_\bfn)
 + k^r_\BW \BBneu V{\bfbeta^r}{\bfn{-}\bfbeta^r} 
 	(\mu_{\bfn+\bfalpha^r{-}\bfbeta^r} - \mu_\bfn)
\end{align}
for $\bfmu=(\mu_\bfn)_{\bfn \in \calN}$.

We emphasize that the RRE as well as the CME are uniquely specified if the reaction network \eqref{eq:ReactNetwork}, the reaction rates $k^r_\FW$ and $k^r_\BW$, and the volume $V>0$ are given. 
Hence, there are obviously close relations between both models, in particular for $V\gg 1$, see  \cite{Kurt70SODE, Gill92RDCM, ACTDH05ASHS, AnCrKu10PFSD, WinSch17HMCR}.

So far, we have not used the detailed-balance condition, i.e.,\ we can
even allow for $k^r_\BW = 0$ in the above considerations.  In all
cases, the Kolmogorov forward equation is an infinite-dimensional
linear ODE as in Section \ref{su:CTMCgrad}.  The following result
shows that the detailed-balance condition is inherited from the RRE to
the CME, and moreover a simple equilibrium $\bfw^V$ can be given
explicitly as a product distribution of individual Poisson
distributions, namely $m\mapsto \ee^{-V c_i^*} (V c_i^*)^m/m!$.  This
result can also be retrieved from \cite{AnCrKu10PFSD} by combining
Theorems 4.1 and 4.5 there, where it is shown that the weaker
``complex-balance condition'' is sufficient to guarantee that the
Poisson distribution $\bfw^V$ is an equilibrium for CME. 

For completeness we give a short and independent proof of the fundamental result that for RRE with detailed balance the associated CME satisfies detailed balance again.

\begin{theorem}[Detailed balance for CME]\label{th:DB.CME} 
  Let $\BBneu V \bfalpha \bfn$ be given in the form \eqref{eq:CME.3}.
  Assume that \eqref{eq:ReactKin} has the equilibrium $\bfc_* \in ]0, \infty[^I$
  satisfying the detailed-balance condition \eqref{RRE:DBC}.
Then the equilibrium $\bfw^V:=(w^V_\bfn)_{\bfn\in \calN} \in
\PM(\calN)$ given by
\[
w^V_\bfn =\frac1{Z^*_V} \: \frac{(V\bfc_*)^\bfn}{\bfn !} \quad
\text{with } Z^*_V:= \Pi_{i=1}^I \ee^{V c_i^*}
\]
satisfies the detailed-balance condition for the CME  
\eqref{eq:CME.1}, namely
\begin{align*}
\forall \, r=1,\ldots,R\ \forall \, \bfn\in
\calN: \ \  
 k^r_\FW \BBneu V {\bfalpha^r} {\bfn}
         w^V_{\bfn + \bfalpha^r} 
 = 
 k^r_\BW \BBneu V {\bfbeta^r} {\bfn} 
 		w^V_{\bfn + \bfbeta^r}
 =  \kappa^r_*V w^V_\bfn
 =: {\wh\nu}_V^{\bfn,r}  . 
\end{align*} 
\end{theorem}
\begin{proof}
For each reaction we obtain the relation 
\[
k_\FW^r \BBneu V {\bfalpha^r} {\bfn}
         w^V_{\bfn{+}\bfalpha^r}
	= k^r_\FW
		\frac{V(\bfn{+}\bfalpha^r)!}
			 {V^{|\bfalpha^r|}\bfn!}   
		\frac{(V\bfc_*)^{\bfn{+}\bfalpha^r}}
			 {Z^*_V (\bfn{+}\bfalpha^r)!} 
	= k^r_\FW V \bfc_*^{\bfalpha^r}
		\frac{(V\bfc_*)^\bfn}
			 {Z^*_V \bfn!} 	
	= V \kappa^r_* w_\bfn^V. 
\]
Analogously we obtain the same result for
$k_\BW \BBneu V {\bfbeta^r}{\bfn} w^V_{\bfn + \bfbeta^r}$, where the
detailed-balance condition \eqref{RRE:DBC} is used in the definition
of $\kappa^r_*$.
\end{proof}

Using the detailed-balance coefficients ${\wh\nu}_V^{\bfn,r}$ we
can rewrite the operator $\ol \calB^r_V$ from \eqref{eq:CME.1} in
a symmetrically balanced form as
\begin{equation}
  \label{eq:CME.5}
  \ol\calB^r_V\bfu = \sum_{n\in \calN}  {\wh\nu}^{\bfn,r}_V
  \Big(\frac{u_{\bfn+\bfalpha^r}} {w^V_{\bfn+\bfalpha^r}} -
  \frac{u_{\bfn+\bfbeta^r}} {w^V_{\bfn+\bfbeta^r}}  \Big)
  \big(\bfe^{(\bfn+\bfbeta^r)} {-} \bfe^{(\bfn+\bfalpha^r)} \big) ,
\end{equation}
where $\bfe^{(\bfm)}$ is the unit vector, i.e.,\ 
$\bfe^{(\bfm)}_\bfn=\delta_{\bfn{-}\bfm}$.

It is important to realize that in general the steady state for the
detailed-balance condition is highly non-unique, because of the
discrete versions 
\[
\calI(\ol\bfn):= \bigset{n\in \calN }{ \bbQ \bfn = \bbQ\ol\bfn} \
\subset \ \calN
\]
of the invariant stoichiometric subspaces
$\bfI(\bfq)=\set{\bfc\in \bfC}{ \bbQ\bfc =\bfq} \subset \bfC$.
Indeed, choosing $\ol n$ arbitrary and defining
$\ol\bfw=(\ol w_\bfn) \in \PM(\calN)$ via
$\ol w_\bfn = \frac1Z w_\bfn^V$ for $\bfn\in \calI(\ol\bfn)$ and
$\ol w_\bfn=0$ elsewhere, we obtain another equilibrium for the CME
\eqref{eq:CME.1}. Defining convex combination we obtain a rich family
of steady states.

The following counterexamples show that the above result, which is
central to our work, cannot be expected for systems not satisfying the
detailed-balance condition.

\begin{example}[Equation without detailed balance]
  \slshape For $a, b \in \N$ we consider the RRE 
\begin{equation}
    \label{eq:NoDBC}
      \dot c = 2a - 4b\, c + 2\,(1{-}c^2),
\end{equation}
which consists of two individual reaction pairs, namely
$ X\REVER{2a}{4b} \emptyset $ and $2X \REVER{1}{1} \emptyset$ with the
individual steady states $c_{(1)}=a/(2b)$ and $c_{(2)}=1$. The joint
steady state of \eqref{eq:NoDBC} is $c_*=(1{+}a{+}b^2)^{1/2}-b$, and
we have detailed balance if and only if $a=2b$.

Building the CME according to \eqref{eq:CME.1}
based on the two reaction pairs we obtain 
\[
  \dot u_n =2aVu_{n-1}-(2aV{+}4bn)u_n + 4b(n{+}1)u_{n+1}
  + V u_{n-2} -\big(V+\tfrac{n(n{-}1)_+}V \big) u_n +
    \tfrac{(n{+}2)(n{+}1)} V u_{n+2} .
\]
For the case $a=2$ and $b=1$, where the detailed-balance condition
holds with $c_*=1=c_{(1)}=c_{(2)}$, we obtain
\[
 \dot u_n =V u_{n-2} + 4Vu_{n-1}
   -\big(5V+4n+\tfrac{n(n{-}1)_+}V \big) u_n  + 4(n{+}1)u_{n+1}
    +\tfrac{(n{+}2)(n{+}1)} V u_{n+2} ,
\]
and it is easy to check that $\wt\bfw^V=(\ee^{-V} V^n/n!)_{n\in \N_0}$ is a steady
state.

However, for $a=7$ and $b=1$ the detailed-balance condition fails with
$c_{(1)}=7/2> c_*=2> c_{(2)}=1$. The CME reads
\[
  \dot u_n = V u_{n-2}+ 14Vu_{n-1} 
   -\big( 15V{+}4n+\tfrac{n(n{-}1)_+}V \big) u_n 
    + 4(n{+}1)u_{n+1}+ \tfrac{(n{+}2)(n{+}1)} V u_{n+2} .
\]
An explicit calculation shows that the Poisson distribution
$\wt \bfw^V$ based on $c_*=2$, i.e.,\ $\wt w^V_n=\ee^{-2V}
(2V)^n/n!$, is not a steady state.  Indeed, inserting $\wt \bfw^V$ into the
right-hand side of the last equation we find (for $n\geq 1$)
\[
{\dot u_n|_{\bfu = \wt\bfw^V} = \frac{\ee^{-2V}(2V)^{n-2}}{n!} \Big(
 {-}12 V^3 + 12nV^2 -3n(n{-}1) V\Big) \neq 0 \text{ for general } n
 \in \N.  }
\]
\end{example}

\begin{example}[Microscopic versus macroscopic detailed balance]
  \slshape
  \ \ We may also consider a RRE that looks macroscopically
  as being in detailed balance,
  but is generated by a microscopic model that is not in detailed
  balance.  The two reactions $\emptyset \overset{2}{\rightharpoonup}
  X$ and $2X \overset{1}{\rightharpoonup} \emptyset$ produce the RRE 
$\dot c = 2 (1 {-}c^2)$ that has the equilibrium
$c_*=1$. However the CME reads 
\[
  \dot u_n = 2V u_{n-1} - \big( 2V + \frac{n(n{-}1)}{ V}\big) u_n +
  \frac{(n{+}2)(n{+}1)}{V} u_{n+2} .
\] 
Again, the Poisson distribution $\wt\bfw^V$ with
$w^V_n=\ee^{-V}V^n/n!$ is \emph{not} the equilibrium:
\[
  \dot u_n|_{\bfu=\wt\bfw^V} = \frac{\ee^{-V}V^{n-1}}{n!}\,
  \Big( 2Vn - V^2 {-} n(n{-}1) \Big) \ \neq \ 0.
\]
Note that the reversible reaction pair
$2 X \REVER{1}{1} \emptyset$ yields the same RRE, and its
associated CME satisfies the detailed-balance condition.
\end{example}

\subsection{Existence and uniqueness of solutions of CME} 
\label{su:ExUniCME}

In this part we establish well-posedness for the CME. We do this by
combining classical results from the theory of Markov chains with
abstract semigroup theory.

For fixed $\bfn_0 \in \calN$ we construct a special Green's
function $p_t(\bfn_0,\cdot)$. General Markov chain theory
(e.g., \cite[Ch.\,2]{Ligg10CTMP}) implies that there exist a unique
minimal solution $[0,\infty[ \times \calN \ni (t, \bfn) \mapsto
p_t(\bfn_0, \bfn)  $ to the backward equation
\begin{align*}
  \dot p_t(\bfn_0, \bfn) 
	&  = \sum_{r = 1}^R\Big(
  k^r_\FW \BBneu V {\bfalpha^r} {\bfn_0{-}\bfalpha^r}
     		\big( p_t(\bfn_0 {-} \bfalpha^r {+} \bfbeta^r , \bfn) 
	 - p_t(\bfn_0, \bfn)\big)
	\\&\qquad\qquad+
   k^r_\BW \BBneu V {\bfbeta^r} {\bfn_0{-}\bfbeta^r}
     		\big( p_t(\bfn_0 {+} \bfalpha^r {-} \bfbeta^r, \bfn) 
		 - p_t(\bfn_0, \bfn) \big)
		\Big)
\end{align*}
associated with the CME with initial condition
$p_0(\bfn_0, \bfn) = \delta_{\bfn_0}(\bfn)$.  This minimal solution is
non-negative and satisfies $p_t(\bfn_0,\bfn)\geq 0$ and
$\sum_{\bfn \in \calN} p_t(\bfn_0, \bfn) \leq 1$, but for general CTMC
it can happen that the latter inequality is strict, which means that
the corresponding Markov chain explodes in finite time.  We will show
that explosion does not happen for CME with detailed balance.

For the functional analytic existence and uniqueness result we use the
sequence spaces
\[
  \ell^p(\calN):=\bigset{\bfu=(u_\bfn)_{\bfn\in \calN}}{
    \sum\nolimits_{\bfn\in \calN}|u_\bfn|^p <\infty}
\]
as well as the weighted spaces 
\[
\rmL^p(\calN,\bfw^V):=\bigset{\bfv=(v_\bfn)_{\bfn\in \calN}}{ \sum\nolimits_{\bfn\in
    \calN}\big|\frac{v_\bfn}{w^V_\bfn} \big|^p <\infty} 
\]
with the corresponding norms and the usual modification for $p=\infty$. 
Now, we consider the transition semigroup $(\calP_t)_{t \geq 0}$ defined by 
\[
(\calP_t \bfv)_\bfn := \sum_{\bfm \in \calN}
p_t(\bfn,\bfm)v_\bfm,\quad  \bfv=(v_\bfm) \in \ell^\infty(\calN),
\]
which we shall study by induction over the number $R$ of reactions
using the Trotter-Kato formula, where the detailed-balance condition
guarantees that each subsystem is a contraction semigroup on
$\rmL^2(\calN,\bfw^V)$.

\begin{theorem}\label{th:semigroup} Assume that the detailed
  balance condition \eqref{RRE:DBC} holds. Then, the semigroup
  $(\calP_t)_{t \geq 0}$ extends to a $\rmC_0$-semigroup of
  contractions on $\rmL^p(\calN, \bfw^V)$ for all $1 \leq p <
  \infty$. Moreover, the semigroup is selfadjoint on $\rmL^2(\calN,
  \bfw^V)$ and Markovian, i.e.,\ $\calP_t {\bf 1} = {\bf 1}$ for all $t
  \geq 0$.
\end{theorem}
A related existence result for the Markov semigroup of the CME
was established in \cite{GauYse14RASC}, which however does not apply
to the case of reversible RRE, because of the restrictions on the
growth of the transition rates. \medskip

\noindent\begin{proof} 
All of the above statements follow from the general theory of
continuous time Markov chains, except for the Markovianity.
To show the latter, we first consider the case of a single reaction,
thus $R = 1$. Each of the irreducible components of the state space
$\calN$ is then one-dimensional (see also \cite{MaaMie15?GSRC2}), and
the Markov chain is a birth-death chain on a countable (possibly
finite) set.

If there exist two components of $\bfalpha - \bfbeta$ with opposite
sign, then each of the irreducible components of the state space
$\calN$ is finite.  Therefore it is clear that the Markov chain does
not explode in finite time.  Suppose now that all components of
$\bfalpha - \bfbeta$ have equal sign, say $\alpha_i - \beta_i \geq 0$
for all $i =1, \ldots, I$, and at least one component is strictly
positive.  Then each of the infinite irreducible components of $\calN$
is of the form
\begin{align*}
\set{ \bfn^{(k)} := \bfn^{(0)} + k (\bfalpha - \bfbeta) }{k \in \N_0} 
\end{align*}
for some $\bfn^{(0)} \in \calN$, and the restricted Markov
process is a birth-death process with birth rate $\mathsf{b}_k$ and death rate
$\mathsf{d}_k$ given by
\begin{align*}
\mathsf{b}_k & := k_\BW \BBneu V {\bfbeta} {\bfn^{(k)}{-}\bfbeta} 
 \ \text{ from
  $\bfn^{(k)}$ to $\bfn^{(k+1)}$} \quad \text{and} \quad \\
\mathsf{d}_k &:= k_\FW \BBneu V {\bfalpha} {\bfn^{(k)}{-}\bfalpha} 
 \ \text{ from $\bfn^{(k)}$ to $\bfn^{(k-1)}$}.
\end{align*}
Reuter's criterion (\cite[Thm.\,11]{Reut57DMPA}) gives a
characterization of non-explosion for birth-death chains; it asserts
that the chain is non-explosive if and only if 
\begin{align*}
  \sum_{k \geq j \geq 0} \mathsf{r}_{j,k} =\infty, 
  \quad \text{where }\mathsf{r}_{j,k}:= 
  \frac{\mathsf{d}_k \cdot \ldots \cdot \mathsf{d}_{j+1}} 
       {\mathsf{b}_{k} \cdot \ldots \cdot \mathsf{b}_j} .
\end{align*}
In our setting we have
$\frac{\mathsf{d}_{k+1}}{\mathsf{b}_k} = (V \bfc^*)^{\bfbeta -
  \bfalpha} \frac{\bfn^{(k+1)}!}{\bfn^{(k)}!}$, so that
$\mathsf{r}_{0,k} = \frac{1}{\mathsf{b}_k}(V \bfc^*)^{k(\bfbeta -
  \bfalpha )} \frac{\bfn^{(k)}!}{\bfn^{(0)}!}$, and therefore
\begin{align*}
  \sum_{k \geq j \geq 0} \mathsf{r}_{j,k} \geq  \sum_{k \geq 0}
  \mathsf{r}_{0,k} 
  \geq \frac{V^{|\bfbeta|-1}}{k_\BW \bfn^{(0)}!}
  \sum_{k \geq 0} \frac{(\bfn^{(k)} {-}\bfbeta)!}
    {(V \bfc^*)^{k(\bfalpha - \bfbeta)}}.
\end{align*}
Since the summands tend to $\infty$ as $k \to \infty$, we infer that
the latter sum is infinite; hence the Markov chain is non-explosive,
or equivalently $\calP_t {\bf 1} = {\bf 1}$ (see
\cite[Thm.\,2.33]{Ligg10CTMP}).

\medskip

The case of multiple reactions follows by induction on the number of
reactions $R$. Indeed, for $\calR \subseteq \{1, \ldots, R \}$, let
$(\calP^{\calR}_t)_{t \geq 0}$ denote the semigroup corresponding to
the reactions $r \in \calR$. Then the Trotter product formula for
contraction semigroups on $\rmL^2(\calN,\bfw^V)$ (see
e.g., \cite{Davi80OPS}) asserts that
\begin{align*}
 \calP^{\{1, \ldots, R+1\}}_t
  = \lim_{n \to \infty} \Big( \calP^{\{1, \ldots, R\}}_{t/n}
   \calP^{\{R+1\}}_{t/n} \Big)^n
\end{align*}
strongly in $\rmL^2(\calN,\bfw^V)$.  Note that we can apply this formula,
since the detailed-balance conditions hold for all reactions
simultaneously, hence all of the semigroups are contractive on the
same space $\rmL^2(\calN,\bfw^V)$. We also observe that the class of
finitely supported functions is a core for each of the generators.
The Markovianity of $\calP^{\{1, \ldots, R+1\}}$ thus follows from the
Trotter formula and the Markovianity of $\calP^{\{1, \ldots, R\}}$ and
$\calP^{\{R+1\}}$.
\end{proof}

\begin{remark}\label{rem:CME-existence}
The mere existence of a probability distribution satisfying the
detailed-balance equations is not sufficient to guarantee
non-explosion of a continuous time Markov chain. It might happen that
the chain jumps infinitely often in a finite time interval,  see
\cite[Sec.\,3.5]{Norr97MC} for an example. The previous result shows
that this phenomenon does not occur in CME satisfying the
detailed-balance condition.
\end{remark}
 
It remains to transfer the results from $\rmL^1(\calN,\bfw^V)$ to
$\ell^1(\calN)$.  Denoting by $\calQ$ the generator of the
$C_0$-semigroup $(\calP_t)_{t\geq 0}$ on $\rmL^1(\calN,\bfw^V)$, we
define the operator
$\calB : \text{Dom}(\calB) \subseteq \ell^1(\calN) \to \ell^1(\calN)$
by
\begin{align*}
  \calB \bfu = \bfw^V \calQ(\bfu / \bfw^V), \quad \text{Dom}(\calB) =
  \set{ \bfu \in \ell^1(\calN) }{\bfu / \bfw^V \in \text{Dom}(\calQ)}.
\end{align*}
This definition of $\calB$ is consistent with the explicit formula for
$\calB$ given above. Since $\calQ$ generates a $C_0$-semigroup of
contractions on $\rmL^1(\calN,\bfw^V)$, it follows that $\calB$ generates
a $C_0$-semigroup $(\calP_t)_{t \geq 0}$ of contractions on
$\ell^1(\calN)$. Furthermore, since $\calP_t$ preserves positivity and
$\calP_t {\bf 1} = {\bf 1}$, it follows that $\PM(\calN)$ is
invariant under the semigroup generated by $\calB$.

As an immediate consequence we obtain global well-posedness for the
CME in $\PM(\calN)$.

\begin{theorem}[Global well-posedness of the CME]\label{thm:CME-well-posed}
  Let the detailed-balance condition \eqref{RRE:DBC}
  hold. Then, for all $\bfu_0 \in \PM(\calN)$ there exists a
  unique mild solution $\bfu : [0,\infty) \to \PM(\calN)$ to the CME
  \eqref{eq:CME.1} satisfying $\bfu(0) = \bfu_0$.
\end{theorem}

\subsection{Gradient structures for CME}
\label{su:grad-CME}

Since the CME is the forward equation associated with a reversible
CTMC, we can formulate it as a gradient flow in view of Proposition
\ref{prop:CTMC.GradSys}. Indeed, for a strictly convex function
$\phi :[0,\infty[ \to \R$ that is smooth on $]0,\infty[$, let us write
\begin{align}
  \label{eq:CME.E.K.gen}
& \calE_V^\phi(\bfu):=  \sum_{\bfn \in \calN} w^V_\bfn \,\phi\big(
\frac{u_\bfn}{w^V_\bfn}\big), \\
&\nonumber
\calK_V^\phi(\bfu):= \sum_{\bfn \in \calN} \sum_{r=1}^R 
  \wh\nu_V^{\bfn,r} \,\Phi\Big(
  \frac{u_{\bfn+\bfalpha^r}}{w^V_{\bfn+\bfalpha^r}}\, , \,
  \frac{u_{\bfn+\bfbeta^r}}{w^V_{\bfn+\bfbeta^r}}  \Big) 
  \, (\bfe^{(\bfn+\bfalpha^r)} {-} \bfe^{(\bfn + \bfbeta^r)} )
\oti (\bfe^{(\bfn+\bfalpha^r)} {-} \bfe^{(\bfn + \bfbeta^r)} ),
\end{align}
where $\Phi$ is defined after \eqref{eq:EK.general},
$\wh\nu^{\bfn,r}_V$ is given in Theorem \ref{th:DB.CME}, and
$\bfe^{(\bfm)}$ denotes the $\bfm$-th unit vector in $\ell^1(\calN)$.  The
following result is then a special case of Proposition
\ref{prop:CTMC.GradSys}.

\begin{proposition}[Quadratic gradient structures for
  CME]\label{prop:CME-GradFlow} 
  If the RRE \eqref{eq:ReactKin} satisfies the detailed-balance
  condition \eqref{RRE:DBC} for a positive steady state
  $\bfc_*=(c^*_i)_{i=1,\ldots,I}$, then the associated CME has the
  gradient structure $(\PM(\calN), \calE_V^\phi, \calK_V^\phi)$ 
  defined in \eqref{eq:CME.E.K.gen}, namely
\begin{align*}
\dot\bfu =
  \calB_V\bfu=- \calK_V^\phi(\bfu) \rmD \calE_V^\phi(\bfu).
\end{align*}
\end{proposition}

In the following we will mainly be concerned with the case that
$\calE_V^\phi$ is the logarithmic entropy, where $\phi$ is the
Boltzmann function $ \LB(z) = z\log z -z +1$. In that case we obtain
\begin{align}
\label{eq:CME.E.K}
&\calE_V(\bfu) := \frac1V \sum_{\bfn\in \calN} w_\bfn^V \LB\big(
                 \frac{u_\bfn} {w_\bfn^V} \big) = 
\frac1V  \sum_{\bfn\in \calN}\big( u_\bfn \log u_\bfn - u_\bfn\log
w^V_\bfn\big)  ,\\   \nonumber 
&\calK_V(\bfu):= V  \sum_{r=1}^R  \sum_{\bfn \in \calN} 
 \wh\nu_V^{\bfn,r} \, \Lambda\Big( \frac{u_{\bfn+\bfalpha^r}}{w^V_{\bfn+\bfalpha^r}} ,
 \frac{u_{\bfn+\bfbeta^r}}{w^V_{\bfn+\bfbeta^r}} \Big) 
  \, (\bfe^{(\bfn+\bfalpha^r)} {-} \bfe^{(\bfn + \bfbeta^r)} )
\oti (\bfe^{(\bfn+\bfalpha^r)} {-} \bfe^{(\bfn + \bfbeta^r)} ),
\end{align}
where the logarithmic mean $\Lambda(a,b)$ is defined in
\eqref{eq:Lambda}. The above definitions do not only restrict to
the entropy function $\phi=\LB$, but also introduce a normalization
with respect to the volume $V$. Hence, $\calE_V$ can be seen as an
entropy per unit volume. The corresponding scaling of $\calK_V$ was
chosen such that the evolution equation
$\dot \bfu = - \calK_V(\bfu)\rmD \calE_V(\bfu)$ is the same as
$\dot\bfu =- \calK_V^\phi(\bfu) \rmD \calE_V^\phi(\bfu) $.

For later purposes we also provide the cosh-type gradient structure
for CME, whose relevance and usefulness is discussed in
\cite{MiPeRe14RGFL, MPPR17NETP, FreLie19?EDTS, MieSte19?CGED}.
Recall the definition of $\sfC^*$ in \eqref{eq:CoshGS} and note
the special scaling via the volume $V$ in \eqref{eq:Psi.V.cosh} below,
which is needed because $\bfPsi^*_{\cosh,V}(\bfu,\cdot)$ is not
scaling invariant.

\begin{proposition}[cosh-type gradient structure for CME] 
 \label{prop:CME-coshGS} 
  If the RRE \eqref{eq:ReactKin} satisfies the detailed-balance
  condition \eqref{RRE:DBC} for a positive steady state
  $\bfc_*=(c^*_i)_{i=1,\ldots,I}$, then the associated CME has the
  gradient structure $(\PM(\calN), \calE_V, \bfPsi^*_{\cosh,V})$ with
  $\calE_V$ from \eqref{eq:CME.E.K} and 
  \begin{equation}
    \label{eq:Psi.V.cosh}
   \bfPsi^*_{\cosh,V}(\bfu,\bfmu):= \frac1V \sum_{r=1}^R \sum_{n\in \calN} 
 \wh\nu_V^{\bfn,r}  \; \Big( \frac{u_{\bfn+\bfalpha^r}} {w^V_{\bfn+\bfalpha^r}} \,
 \frac{u_{\bfn+\bfbeta^r}}{w^V_{\bfn+\bfbeta^r}} \Big)^{1/2}\: 
  \sfC^*\Big( V (\mu_{\bfn + \bfbeta^r}{-} \mu_{\bfn+\bfalpha^r})\Big) .
  \end{equation}
\end{proposition}
\begin{proof} The desired formula
  $\sum_{r=1}^R \ol\calB^r_V \bfu = \rmD_\mu
  \Psi_{\cosh,V}^*(\bfu,-\rmD\calE_V(\bfu)\big)$ follows easily by
  recalling $\ol\calB_V^r$ from \eqref{eq:CME.5} and by using
  $\sqrt{ab}\:(\sfC^*)'\big(\log a - \log b\big)=a{-}b$ and
  $\rmD \calE_V(\bfu)=\frac1V \big( \log (u_\bfn/
  w^V_\bfn)\big)_{\bfn\in \calN}$.
\end{proof}

\section{Liouville and Fokker--Planck equations}
\label{s:Liouville}

For general evolutionary equations one can define a measure-valued
flow in the phase space that is given by transporting the measures
according to the semiflow of the original equation. The evolution
equation describing this measure-valued flow is the Liouville
equation.  For our RRE  $\dot \bfc = -\bfR(\bfc)$ in
$\bfC:={[0,\infty[}^I$ we assume that we have a global semiflow 
$\bfc(t)=\Phi_t(\bfc(0))$ and consider probability measures $\varrho(t,\cdot)
\in \PM(\bfC)$ that are obtained by transporting $\varrho_0$ with
$\Phi_t$, namely
\[
\varrho(t,\cdot)=\Phi_t^\#\varrho_0, \quad \text{i.e.,\ } \forall\,\psi\in 
\rmC_b(\bfC): \int_\bfC \psi(\bfc) \varrho(t,\rmd\bfc) = \int_\bfC
\psi(\Phi_t(\bfc)) \varrho_0(\rmd\bfc).
\]
In particular, if $\varrho_0 = \sum_{k=1}^m a_k \delta_{\bfc^k_0}$,
then
$\varrho(t,\cdot)= \sum_{k=1}^m a_k \delta_{\Phi_t(\bfc^k_0)}(\cdot)$.

It is now easy to see that $t\mapsto \varrho_t \in \PM(\bfC)$ satisfies
the  Liouville equation 
\begin{align}\label{eq:Liouville}
  \partial_t\varrho(t,\bfc) = \div\big(\varrho(t,\bfc) \bfR(\bfc)\big),
\end{align}
in the sense of distributions. We will regard \eqref{eq:Liouville} as
an evolution equation in the space $\PM(\bfC)$.  We will not always
notationally distinguish between an absolutely continuous probability
measure and its density, but if we want to distinguish them
we will write $\varrho(\rmd\bfc) = \rho(\bfc) \dd \bfc$ with
$\rho \in \rmL^1(\bfC)$.

The goal of this section is to give a rigorous connection between the
CME for $V\to \infty$ and the Liouville equation in terms of the
associated gradient structures.  

\subsection{The Liouville equation as a gradient system}
\label{su:Lio-GradSyst}

We show that the gradient structure $\dot \bfc = -\bfR(\bfc)=-
\bbK(\bfc) \rmD E(\bfc)$ for the RRE, which was discussed in Section
\ref{su:RRE.GradSys}, induces a natural gradient structure for
the Liouville equation.
Consider the ``Otto-Wasserstein-type'' Onsager operator
$\bfK(\varrho)$ that acts on functions $\xi : \bfC \to \R$ via
\begin{align*}
	\bfK(\varrho)\xi = -\div\big( \varrho \,  \bbK \, \nabla \xi\big),
\end{align*}
where $\div$ and $\nabla$ are taken with respect to $\bfc \in \R^I$.
We also consider the affine potential energy functional
$\bfE : \PM(\bfC) \to [0,+\infty]$ defined by
\begin{align}
  \label{eq:bfE.def}
 \bfE(\varrho) = \int_\bfC E(\bfc) \dd \varrho(\bfc).
\end{align}
In the next result we identify the formal gradient structure for the
Liouville equation.

\begin{proposition}[Gradient structure for the Liouville equation]
\label{prop:Liouville-GradFlow}
  If the RRE \eqref{eq:ReactKin} satisfies the detailed-balance
  condition \eqref{RRE:DBC} for a positive steady state
  $\bfc_*=(c^*_i)_{i=1,\ldots,I}$, then the associated Liouville
  equation has the gradient structure $(\PM(\bfC) , \bfE, \bfK)$,
  namely
\begin{align}\label{eq:KE-gradflow}
\dot\varrho = - \bfK(\varrho) \rmD \bfE(\varrho)= \div\big( \varrho
\bbK \nabla E\big) = \div(\varrho \bfR) .
\end{align}
\end{proposition}
\begin{proof}
  Let $\varrho \in \PM(\bfC)$ and let $\sigma \in \Meas(\bfC)$ be a
  signed measure of finite total variation such that
$\sigma(\bfC) = 0$ and $\varrho + h \sigma \in \PM(\bfC)$ for
  $|h|$ sufficiently small. Then we have
\begin{align*}
 \frac{\bfE(\varrho + h \sigma) - \bfE(\varrho)}{h} 
      = \int_\bfC E(\bfc) \dd \sigma(\bfc),
\end{align*}
hence $\rmD \bfE(\varrho) = E$ for all $\varrho$. Therefore,
$
- \bfK(\varrho) \rmD \bfE(\varrho) = \div\big( \varrho \, \bbK \, \nabla E\big)
= \div( \varrho \bfR).  $
The gradient flow equation $\dot\varrho = - \bfK(\varrho) \rmD \bfE(\varrho)$
is thus given by the Liouville equation \eqref{eq:Liouville}.
\end{proof}

\subsection{Passing to the limit from CME to Liouville}
\label{su:CME-Liouv}

In this section we shall demonstrate that the gradient flow structure
for the CME converges in a suitable sense to the gradient structure
for the Liouville equation if $V \to \infty$.

More precisely, we will show that after a suitable $V$-dependent
embedding of $\PM(\calN)$ into $\PM(\bfC)$ the proper scalings of the
functionals $\calE_V$ and
$\Psi_V^*:(\bfu,\bfmu)\mapsto \frac12\bfmu \cdot \calK_V(\bfu)\bfmu$
converge in the sense of $\Gamma$-convergence to the corresponding
structures for the Liouville equation given by the gradient system
$(\PM(\bfC),\bfE,\bfK)$, see Section \ref{su:GamEntro} to
\ref{su:GamDiss}.  Following the approach in \cite{SanSer04GCGF,
  Serf11GCGF, Miel16EGCG}, and in particular \cite{LMPR17MOGG}, we are
then able to establish the convergence for $V\to \infty$ of solutions
$\bfu^V:{[0,\infty[}\to \PM(\calN) $ of the CME
$\dot \bfu^V = - \calK_V(\bfu^V) \rmD \calE_V(\bfu^V)$ to the solution
$\varrho:{[0,\infty[}\to \PM(\bfC)$ of the Liouville equation
$\dot \varrho = - \bfK(\varrho) \rmD \bfE(\varrho)$, thereby
recovering Kurtz' result \eqref{eq:I.Kurtz}, see Section
\ref{su:GamCME2Lio}.

The main tool for proving this evolutionary $\Gamma$-convergence for
gradient systems is the so-called \emph{energy--dissipation principle}, 
 cf.\ \cite[Sec.\,3.3]{Miel16EGCG}, which
states that $\bfu^V$ solves the CME if and only if for all $T>0$ the following
energy-dissipation estimate holds: 
\begin{equation}
  \label{eq:EDP-CME}
  \calE_V(\bfu^V(T)) + \int_0^T \!\! \Big(\Psi_V(\bfu^V,\dot \bfu^V) +
  \Psi_V^*\big(\bfu^V, {-}\rmD\calE_V(\bfu^V)\big)\Big) \dd t
  \leq \calE_V(\bfu^V(0)),
\end{equation}
where we use the quadratic dissipation potential $\Psi_V$ and its
Legendre dual $\Psi^*_V$ defined via $\Psi_V^* (\bfu,\bfmu) := \frac12\langle
\bfmu,\calK_V(\bfu)\bfmu\rangle$ with $\calK_V$ from
\eqref{eq:CME.E.K}, namely
\begin{subequations}
  \label{eq:Psi*Vu-xi}
 \begin{align}
\label{eq:Psi*Vu-xi.a}
 \Psi_V^*(\bfu,\bfmu)
&= \frac V2 \sum_{\bfn\in \calN}\sum_{r=1}^R \wh\nu^{\bfn,r}_V 
      \Lambda\big( \frac{u_{\bfn+\bfalpha^r}}{w^V_{\bfn+\bfalpha^r}} , 
               \frac{u_{\bfn+\bfbeta^r}}{w^V_{\bfn+\bfbeta^r}} \big)  
\big( \mu_{\bfn+\bfalpha^r} {-} \mu_{\bfn+\bfbeta^r}\big)^2
\\ 
\label{eq:Psi*Vu-xi.b}
& =
\frac V2 \sum_{\bfn\in \calN}\sum_{r=1}^R \Lambda 
\big( k^r_\FW \BBneu V{\bfalpha^r}{\bfn} u_{\bfn+\bfalpha^r} ,
      k^r_\BW \BBneu V{\bfbeta^r} {\bfn} u_{\bfn+\bfbeta^r}\big)  
\big( \mu_{\bfn+\bfalpha^r} {-} \mu_{\bfn+\bfbeta^r}\big)^2,
\end{align}
\end{subequations}
where the second form uses Theorem \ref{th:DB.CME} and is especially
useful to perform the limit $V \to \infty$, see the proof of Proposition
\ref{pr:Psi*Limsup}.

We refer to \cite{DanSav14LNGF,Miel16EGCG} for this equivalence and
general methods for proving such results on evolutionary
$\Gamma$-convergence. In \cite{DisLie15GSMC} a similar approach was
used to establish the convergence of CTMC to a Fokker--Planck
equation. However, there the convergence of a parabolic equation is
established, where upper and lower bounds of the density can be
used. Here, the importance is that our limit measures $\varrho(t)$ may
not have densities; indeed, because we want to recover the Kurtz
result \eqref{eq:I.Kurtz} we are interested in the ``deterministic
case'' $\varrho(t)=\delta_{\bfc(t)}$. So our analysis has to be more
careful in dealing with general limit measures.  For this, we use the
dualization approach introduced in \cite{LMPR17MOGG} where
$t\mapsto \Psi_V(\bfu^V,\dot\bfu^V)$ is estimated from below by
$\langle \dot\bfu^V, \bfmu^V \rangle- \Psi^*_V(\bfu^V,\bfmu^V)$ for
suitably chosen recovery functions $t\mapsto \bfmu^V(t)$.

In order to compare probability measures on different spaces $\calN$
and $\bfC$, we consider a suitable embedding $\iota_V : \PM(\calN)
\to \PM(\bfC)$. Here $\iota_V(\bfu)$ is simply obtained by assigning
the mass of $\bfu$ at $\bfn \in \calN$ uniformly to the cube
\begin{align*}
  A_\bfn^V := \left[\tfrac{n_1}{V}, \tfrac{n_1 + 1}{V} \right[ \times \cdots
  \times \left[\tfrac{n_I}{V}, \tfrac{n_I + 1}{V}\right[ \subseteq \bfC.
\end{align*}
More explicitly, $\iota_V(\bfu)$ is given by
\begin{align}
  \label{eq:iota-V}
\iota_V : \PM(\calN)\to \PM(\bfC) ; \quad \bfu \mapsto
\iota_V(\bfu)=\varrho =\rho\dd\bfc \text{ with }
\rho(\bfc) := V^I \sum_{\bfn \in \calN} u_\bfn \INDIC_{A_\bfn^V}(\bfc),
\end{align}
where $\INDIC_A$ denotes the indicator function with $\INDIC_A(b)=1$
for $b\in A$ and $0$ otherwise. The corresponding dual operation
acting on functions $\xi \in \rmC_\rmb(\bfC)$ is given by
\begin{equation}
  \label{eq:EDP-iota*}
\iota_V^* : \rmC_\rmb(\bfC) \to \ell^\infty(\calN); \quad 
	(\iota_V^* \xi)(\bfn) = V^I
\int_{\bfc\in A^V_\bfn} \xi(\bfc)\dd \bfc.\medskip
\end{equation}

The final convergence result will be formulated in Theorem
\ref{th:CME2Lio}, which will be a direct consequence of the following
three estimates
\\[0.2em]
\begin{tabular}{@{}l@{\quad}c@{\;}c@{\;}r@{\:}c@{\:}l}
Section \ref{su:GamEntro}&  $\iota_V(\bfu^V) \weaks \varrho $ & $
\Rightarrow $&$ \bfE(\varrho) $&$\leq$&$ \liminf\limits_{V\to \infty} \calE_V(\bfu^V);$\\
Section \ref{su:GamDual}&  $\iota_V(\bfu^V) \weaks \varrho $ & 
$\Rightarrow$ &   $\Psi^*_\Lio (\varrho,\rmD\bfE(\varrho))$&$ \leq
$&$\liminf\limits_{V\to\infty} \Psi_V^*(\bfu^V,\rmD \calE_V(\bfu^V))$; \\
Section \ref{su:GamDiss}& $\iota_V(\bfu^V)\weaks \varrho,\ \xi\in\rmC^1_\rmc(\bfC)$ &  $\Rightarrow $ & 
$\Psi^*_\Lio (\varrho,\xi) $&$\geq $&$\limsup\limits_{V\to \infty}
\Psi^*_V(\bfu^V,\iota_V^*\xi)$; \\
\end{tabular}\\
where the dual dissipation potential $\Psi^*_\Lio $ is defined
via 
\[
\Psi^*_\Lio (\varrho, \xi) = \frac12\int_\bfC \nabla \xi(\bfc) \cdot
\bbK(\bfc) \nabla\xi(\bfc) \dd\varrho(\bfc). 
\]
We will see in Section \ref{su:GamCME2Lio} that the limsup estimate
for the dual potential $\Psi^*_V$ in Section \ref{su:GamDiss} provides
a weak form of a liminf estimate for the primal potential $\Psi_V$.
   
A fundamental fact of the chosen gradient structures of the underlying
Markov processes is that all the three terms in the energy-dissipation
principle define convex functionals, which is of considerable help
in proving the desired liminf estimates. 
Note that the convergence $\iota_V(\bfu^V) \weaks \varrho $ is rather weak. However, we can use that the coefficients of the transition rates defining the CME are quite regular, so that the other parts in the integral converge in a much better sense. Moreover, the functionals $\varrho\mapsto \bfE(\varrho)$ and $\varrho \mapsto \Psi^*_\Lio (\varrho,\rmD\bfE(\varrho))$ are in fact linear in $\varrho$.

\subsection{\texorpdfstring{$\Gamma$}{Gamma}-limit of the relative entropies}
\label{su:GamEntro}

We also
define $\bfX_V:=\iota_V(\PM(\calN)) \subset \PM(\bfC)$
and $W_V = \iota_V(\bfw^V) \in \bfX_V$ and consider the functionals
\[ \wh \calE_V : \PM(\bfC) \to [0,\infty], \qquad
\wh \calE_V(\varrho)=\left\{ \ba{cl}  \wt\calE_V(\varrho)& \text{ if }\varrho
\in \bfX_V,\\ \infty&\text{ otherwise},   \ea \right. 
\]
where $\wt \calE_V : \PM(\bfC) \to [0,\infty]$ is defined via 
\[
\wt\calE_V(\varrho) =  \frac{1}{V} \text{Ent}(\varrho| W_V\rmd\bfc ) =
\left\{ \ba{cl}
 \frac{1}{V} \int_\bfC \LB(\rho/W_V) W_V \dd \bfc& \text{for
 }\varrho=\rho\dd c, \\ \infty &\text{otherwise} . \ea\right. 
\]
These definitions are chosen such that $\calE_V(\bfu)=
\wt\calE_V(\iota_V(\bfu))= \wh\calE_V(\iota_V(\bfu))$ for all $\bfu
\in \PM(\calN)$.

Finally we define a natural inverse of $\iota_V $, namely 
\begin{align}
\label{eq:kappa-V}
\varkappa_V:\ \PM(\bfC)\to \PM(\calN);\  \varrho \mapsto \Big( \varrho\big(
A_\bfn^V \big) \Big)_{\bfn  \in \calN}, 
\end{align}
such that $P_V := \iota_V \circ \varkappa_V$ is a projection from
$\PM(\bfC)$ onto $\bfX_V \subset \PM(\bfC)$.

To understand the limit of $\wh\calE_V$ for $V\to \infty$ we will use
the representation
\begin{align}
  \wt\calE_V(\rho\rmd\bfc)  \nonumber
  &=\frac1V\int_\bfC \LB(\rho/W_V) W_V \dd \bfc 
	= \int_\bfC \Big(\frac1V \rho \log\rho 
	   + \rho E_V(\bfc) \Big)  \dd \bfc  
  \\
  \label{eq:def.E.V}
  \text{with }\ \
  &E_V(\bfc)= \frac1V \log \Big( \frac1{W_V(\bfc)}\Big)
    = - I\,\frac{\log V}{V} - \frac1V \log \bfw^V_\bfn \ \ 
    \text{ for } \bfc\in A^V_\bfn. 
\end{align}
In Lemma \ref{lem:E-compare} below we will show that $E_V$
converges pointwise to $E$ as defined in \eqref{eq:RRE.E.K}. To
quantify the latter convergence, we use the classical lower and upper
bounds of \cite{Nanj59NSF} for Stirling's formula:
\begin{equation}
  \label{eq:Stirling}
  \begin{aligned}  \forall \, n\in \N_0 : \quad 
    &n! = \sqrt{2\pi k_n} \Big(\frac{n}\ee\Big)^n \quad
    \text{with }k_0=\frac{1}{2\pi}\\ 
    &\text{and }
    k_n=n+\frac16 +\frac{\gamma_n}{124/5+72 n} \text{ with } 
    \gamma_n\in [0.9,1]  \text{  for } n\geq 1.
\end{aligned}
\end{equation}
Using this estimate and recalling $E$ from \eqref{eq:RRE.E.K} we
obtain the following estimate.

\begin{lemma}[Pointwise bound for $E_V$]
\label{lem:E-compare}
For all $\bfc^*>0$ there exist $K_*>0$ and $V_*>0$ such that for all
$V\geq V_*$ the following bounds hold:
\begin{align}
\label{eq:E-compare}
| E_V(\bfc) -  E(\bfc) |
\leq \frac{K_*}{V} \big(\log V + E(\bfc)\big) \text{ \ for all }\bfc \in
\bfC.  
\end{align}
\end{lemma}
\begin{proof} 
We decompose the error via
\begin{equation}
    \label{eq:ErrorDecomp}
    E_V(\bfc)-E(\bfc)=\big(E_V(\bfc){-} E(\frac1V \bfn)\big) +
    \big(E(\frac1V \bfn){-} E(\bfc) \big)
\end{equation}
with $\bfn$ defined by $\bfc \in A_\bfn^V$.  For the second term we
use the convexity of $\lambda_{\rmB}$ and the estimate
$\log z \leq 1+\LB(z)$. Hence, we have
\begin{equation*}
\begin{aligned}
c_i^* \Big[ \LB\big(\frac{c_i}{c_i^*}\big) - \LB\big(\frac{n_i}{V
c_i^*}\big) \Big] &\leq \big(c_i {-} \frac{n_i}{V}\big)
\log\big(\frac{c_i}{c_i^*}\big) \leq
\frac1V\Big(1+\LB\big(\frac{c_i}{c^*_i} \big) \Big) \\
& \leq \frac{\max\{I,1/c^*_i\}} V \Big(\frac1I
+c^*_i\LB\big(\frac{c_i}{c^*_i} \big) \Big).
\end{aligned}
\end{equation*}
Summing this inequality over $i=1, \ldots,I$ we obtain the upper bound 
\begin{equation}
\label{eq:E-1} E(\bfc) - E\big(\tfrac1V {\bfn}\big)  \leq 
\frac{K_1}{V} \big(1 {+} E(\bfc)\big) \text{ \ with }K_1=\max\{I,
1/c^*_1, ... ,1/c^*_I\}. 
\end{equation}
For the opposite direction we use (a) that $\LB$ decreases on
$[0,1]$ and the convexity of $\LB$ which implies
(b) $\LB(z_1)-\LB(z_2)\leq \LB(0)-\LB(z_2{-}z_1)$ for
$0 \leq z_1 \leq z_2$. This yields
\begin{equation}
\begin{aligned}
 \label{eq:E-2}
 E\big(\tfrac1V {\bfn}\big) - E(\bfc) & = \sum_{i=1}^I c_i^*
 \Big[\lambda_{\rmB}\big(\frac{n_i}{Vc_i^*}\big) -
 \lambda_{\rmB}\big(\frac{c_i}{c_i^*}\big) \Big]
 \;\overset{\text{(b)}}\leq \; \sum_{i=1}^I
 c_i^* \Big[\lambda_{\rmB}(0) - \lambda_{\rmB}\big(\frac{c_i}{c_i^*}
 {-} \frac{n_i}{Vc_i^*}\big) \Big]
 \\
 & \overset{\text{(a)}}\leq \, \sum_{i=1}^I c_i^*
 \big[\lambda_{\rmB}(0) - \lambda_{\rmB} \big(\tfrac{1}{Vc_i^*}\big)
 \big] = \frac1V \sum_{i=1}^I\! \big( 1 {+} \log (V c_i^*)\big)
 \, \overset{\text{(c)}}\leq \, 2I\, \frac{\log V}{V},
\end{aligned}
\end{equation}
if $V\geq V_1^*:= \max\bigset{\max\{1/c^*_i, \ee c^*_i
  \}}{i=1,...,I}$, where $Vc^*_i\geq 1$ and $V\geq \ee c^*_i$ are
needed in (a) and (c), respectively. Together with \eqref{eq:E-1} this
controls the second error term in \eqref{eq:ErrorDecomp}, viz.\ 
\begin{equation}
  \label{eq:ErrorTerm1}
  \big|E(\bfc)- E\big( \tfrac1V \bfn) \big| \leq \frac{K_2}V \big(
  \log V + E(\bfc)\big) \ \text{ for } V\geq V^*_2 = \max\{ \ee,
  V^*_1\},
\end{equation}
where $K_2=\max\{2I, K_1 \}$.\medskip 
 
For controlling the first error term in \eqref{eq:ErrorDecomp}
we use  \eqref{eq:def.E.V} and obtain the identity 
\begin{align}
\label{eq:E-3}
 E_V(\bfc) - E\big(\tfrac1{V}{ \bfn }\big)
  =  - I \frac{\log V}{V} + \frac{1}{2V} \sum_{i=1}^I \log (2\pi k_{n_i})
  \quad\text{for all } \bfc\in A^V_\bfn,
\end{align}
with $k_n$ from \eqref{eq:Stirling}.  Because of
$2\pi k_n\geq 1$ we obtain, for all $V\geq 1$, the lower bound 
\[
  E_V(\bfc) - E\big(\tfrac1{V}{ \bfn }\big) \geq - I \frac{\log V}{V}
  \geq - \frac IV \big( \log V +E(\bfc)\big).
\]

For the upper bound we use $2\pi k_0=1$ and
$2\pi k_n \leq  8 n $ for $n \geq 1$.  Hence for $n_i \geq 1$
we obtain, using again the estimate $\log z \leq 1+\LB(z)$,
\[
  \log(2\pi k_{n_i}) \leq \log(8n_i) \leq \log(8c^*_i V) +
  \log\big(\frac{c_i}{c^*_i}\big) \leq \log V + \log(8 \ee c^*_i) +
  \frac{1}{c^*_i} c^*_i\LB\big( \frac{c_i}{c^*_i}\big).
\]
Summation over $i=1,\ldots,I$ yields, for all $\bfc\in A^V_\bfn$
and $V \geq V_3^*:= 8\ee \max\{c^*_1,\ldots,c^*_I\}$,
the upper bound  
\begin{align*}
 E_V(\bfc) - E\big(\frac{ \bfn }{V}\big)
  \leq  \frac{K_3}{V}  E(\bfc) \text{ \ with } 
K_3 = \max\Big\{ \frac1{2c^*_1},\ldots,\frac1{2c^*_I}  \Big\}.
\end{align*}

Together with the lower estimate we control the first error term in
\eqref{eq:ErrorDecomp} via 
\begin{align*}
	\big| E_V(\bfc) -  E(\tfrac1V \bfn)| \leq \frac{K_4}V \big(
  \log V + E(\bfc)\big) \ \text{for } V \geq V^*_4=\max\{1, V^*_3\},
\end{align*}
where $K_4=\max\{I,K_3\}$.

Adding the estimates for first and the second error term
\eqref{eq:ErrorDecomp} we obtain the desired estimate
\eqref{eq:E-compare} with the choices $K_*=K_2+ K_4$ and
$V_*=\max\{V_2^*,V_4^* \} $.
\end{proof}

For consistency of notation we remark that $\wt\calE_V(\varrho)$ can be
rewritten as
\begin{align*}
  \wt\calE_V(\rho \rmd\bfc) =\int_\bfC \big(\frac1V \log \rho(\bfc) + E_V(\bfc)
  \big) \rho(\bfc) \dd \bfc\;,
\end{align*}
provided that this integral exists. The limit functional $\bfE$ is
given by
\begin{equation}
  \label{eq:bfE-def}
  \bfE:\PM(\bfC)\to [0,\infty]; \quad \varrho \mapsto \int_\bfC E(\bfc)
  \dd \varrho(\bfc), 
\end{equation}
where we use that $E$ is a continuous and non-negative function, so
that $\bfE$ can be defined everywhere but attains the value $+\infty$
if $\varrho$ does not decay suitably at infinity. We will use the
following semi-continuity result. 

\begin{lemma}[Lower semi-continuity of $\bfE$]\label{lem:bfE} 
For sequences $(\varrho_k)_k \subset \PM(\bfC)$ with $\varrho_k \weaks
\varrho_\infty$, we have 
$\bfE(\varrho_\infty) \leq \liminf_{k\to \infty}\bfE(\varrho_k) $.  
\end{lemma} 
\begin{proof} For cut-off functions $\chi\in \rmC_\rmc(\bfC)$
  with $\chi(\bfc) \in [0,1]$ we have $\bfE(\chi \varrho_k)\to
  \bfE(\chi\varrho_\infty)$ by weak* convergence and continuity of
  $E$. Using $\chi\leq 1$ yields $ \bfE(\chi\varrho_\infty) \leq
  \liminf_{k\to \infty} \bfE(1 \varrho_k)$. Choosing a non-decreasing
  sequence $\chi_n$ with $\chi_n(\bfc)\to 1$
  for all $\bfc\in \bfC$ we have $\bfE(\chi_n \varrho_\infty) \to \bfE(1
  \varrho_\infty)$ by Beppo Levi's monotone convergence, and the
  assertion follows.
\end{proof} 

The following result gives the $\Gamma$-convergence of $\calE_V$ to
$\bfE$ with respect to the sequential weak* convergence as well as the
equi-coercivity.

\begin{theorem}[$\Gamma$-convergence of $\calE_V$ to $\bfE$]\label{th:EV-Gamma}
Let $\wh\calE_V$ and $\bfE$ be defined on $\PM(\bfC)$ as above. Then we
have the following properties:

(a) \emph{Compactness / equi-coercivity:}
\begin{equation}
  \label{eq:EV-equicoerc}
  \exists\, V_*,C, c>0\ \forall\, V\geq V_*\ \forall \, \varrho \in
  \PM(\bfC): \quad \wh\calE_V(\varrho) \geq -C + c \bfE(\varrho).
\end{equation}

(b) \emph{Weak* liminf estimate:}
\begin{equation}
  \label{eq:EV-liminf}
  \varrho_V \weaks \varrho \text{ in } \PM(\bfC) \quad
  \Longrightarrow \quad \liminf_{V\to \infty} \wh\calE_V(\varrho_V) \geq
  \bfE(\varrho). 
\end{equation}

(c) \emph{Limsup estimate / recovery sequence:} 
\begin{equation}
  \label{eq:EV-recovery}
 \forall\, \wh \varrho \in \PM(\bfC) \ \exists\, (\wh\varrho_V)_{V\geq 1}: 
  \quad \wh\calE_V(\wh\varrho_V) \to \bfE(\wh\varrho) \text{ and } \wh\varrho_V
  \weaks \wh\varrho,  
\end{equation}
where we may take $\wh\varrho_V= P_V \wh\varrho = \iota_V\big(
\varkappa_V(\wh\varrho)\big)$.
\end{theorem}
\begin{proof}
Obviously it is sufficient to show the lower bound (a) and the liminf
estimate (b) for the smaller functional $\wt\calE_V$, and for $\varrho = \rho \dd \bfc$ with $\rho \in \rmL^1(\bfC)$
(resp. $\varrho_V = \rho_V \dd \bfc$ with $\rho_V \in \rmL^1(\bfC)$).
We use the elementary convexity estimate 
\[
\forall\, r\geq 0,\ a,w>0: \quad w\LB(r/w) = r\log (r/w) - r
+ w \geq r \log(a/w)  -a + w. 
\]
We choose $r(\bfc)=\varrho(\bfc)$, $w(\bfc)=W_V(\bfc)$, and
$a(\bfc) = \ee^{-|\bfc|_1} = \Pi_{i=1}^I \ee^{-c_i} >0$. Note that
$a\in \rmL^\infty(\bfC)\cap \PM(\bfC)$ and $W_V/a$ is bounded from
above, for any fixed $V$. Hence,
$\bfc \mapsto \log(a(\bfc) / W_V(\bfc))= -|\bfc|_1 + V E_V(\bfc) $ is
bounded from below, and we can integrate the above estimate to obtain
the lower bound
\[
  \wt\calE_V(\varrho) \geq \frac1V\int_\bfC \log\big(
  a(\bfc)/W_V(\bfc)\big) \dd \varrho(\bfc) = \int_\bfC \Big(E_V(\bfc)
  - \frac{|\bfc|_1} V \Big) \dd \varrho(\bfc).
\]
Since there exists a constant $K_1 > 0$ such that $|\bfc|_1\leq
K_1\big(1{+}E(\bfc)\big)$ and since $E_V$ satisfies the lower bound in
\eqref{eq:E-compare}, we obtain the lower bound 
\begin{align*}
 \wt\calE_V(\varrho)
  & \geq \int_\bfC  E(\bfc) \dd\varrho(\bfc)
	 - \frac{K_*{+}K_1}V
	    \int_\bfC \big(\log V   {+}E(\bfc)\big)\dd \varrho(\bfc)
  \\
  &
  = \bfE(\varrho) - \frac{K_*{+}K_1}V\big( \log V + \bfE(\varrho)\big).
\end{align*} 
This immediately implies \eqref{eq:EV-equicoerc} in part (a) with
$V_* / \log V_* =2(K_*{+}K_1)$.  Moreover, if
$\varrho_V \weaks \varrho$ then we have the lower bound
$\wh\calE_V(\varrho_V) \geq \bfE(\varrho_V) - \frac{K_*{+}K_1} V \big(
\log V + \bfE(\varrho_V)\big)$ and the liminf estimate
\eqref{eq:EV-liminf} follows from Lemma \ref{lem:bfE}.

To show part (c) we use the indicated recovery sequence and the upper
bounds for $E_V$ from \eqref{eq:E-compare}. For a given $\wh \varrho\in
\PM(\bfC)$ we define $\wh\varrho_V=\iota_V(\varkappa_V(\wh\varrho))$. For an
arbitrary continuous and bounded test function $\psi$ we define the
piecewise constant approximation $\psi_V$ via averaging over
$A_\bfn^V$. We obtain
\[
\int_\bfC \psi (\bfc)\dd \wh\varrho_V (\bfc) = \int_\bfC \psi_V (\bfc)\dd
\wh\varrho_V (\bfc)  = \int_\bfC \psi_V (\bfc)\dd
\wh\varrho (\bfc) \to \int_\bfC \psi(\bfc) \dd \wh\varrho(\bfc),
\] 
where the convergence follows via Lebesgue's dominated
convergence from the pointwise convergence $\psi_V \to \psi$ and the
uniform boundedness of $\psi_V$. Thus, we conclude 
$\wh\varrho_V\weaks \wh\varrho$.

To show convergence of $\wh\calE_V(\wh\varrho_V)$ it suffices to prove the
upper bound $\limsup_{V\to \infty} \wh\calE_V(\wh\varrho_V) ${\linebreak[3]}$\leq
\bfE(\wh\varrho)$. For this we use the bound $\wh\rho_V(\bfc)\leq
V^I=1/\text{vol}(A^V_\bfn) $ and the fact that $\wh\rho_V$ and $E_V$
are constant on the same cubes to obtain
\[
\wh\calE_V(\wh\varrho_V) 
	= \int_\bfC \big( \frac{\log \wh\rho_V(\bfc)}V 
			+ E_V(\bfc)\big) \dd \wh\rho_V(\bfc) 
	\leq   \frac{I\log V}V + \int_\bfC E_V(\bfc)
	 \dd \wh\varrho(\bfc),
\]
where now only the measure $\wh\varrho$ is left.  The first term tends to
0 for $V\to \infty$, and the second can be estimated from above using
the upper estimate in \eqref{eq:E-compare}, which yields 
\begin{align*}
\ts
\wh\calE_V(\wh\varrho_V) 
	& \leq  \frac{I\log V}V 
			+  \int_\bfC 
					\! E(\bfc) + \frac{K_*}{V} \big(\log V {+} E(\bfc)\big)
				\dd \wh\varrho(\bfc)
		=   \big(1 + \frac{K_*}{V}\big) \bfE(\wh\varrho)
			+ \frac{I{+} K_*}{V}\log V .
\end{align*}
This implies the desired upper bound for $V\to \infty$, and the proof
is complete.
\end{proof}

\subsection{A liminf estimate for the dual dissipation functional}
\label{su:GamDual}

Here we provide the liminf estimate for the dual dissipation potential
$\Psi_V^*(\bfu^V,\rmD \calE_V(\bfu^V))$ based on the lower bound 
\begin{equation}
\label{eq:Lio.slope}
	\Psi^*_\Lio (\varrho,\rmD\bfE(\varrho)) 
		= \frac12 \int_\bfC 
						\nabla E(\bfc)\cdot \bbK(\bfc) \nabla E(\bfc) 
					\dd \varrho(\bfc).
\end{equation}
We observe that the latter term is linear in $\varrho$ while the former
term is convex in $\bfu^V$. Indeed, introducing the convex function
$G(a,b)=(a{-}b)(\log a - \log b)$ for $a,b>0$ and noting the
relation $\Lambda(a,b)(\log a - \log b)^2 = G(a,b)$ we have 
\begin{equation}
\label{eq:Psi*V-all}
 \Psi_V^* \big( \bfu^V,\rmD \calE_V(\bfu^V) \big)
  = \frac{1}{ 2  V} 
   \sum_{r=1}^R 
     \sum_{n\in \calN} 
 	\wh\nu^{\bfn,r}_V 
         G \Big(
	  \frac{u_{\bfn+\bfalpha^r}^V}{w^V_{\bfn+\bfalpha^r}}, 
	   \frac{u_{\bfn+\bfbeta^r}^V}{w^V_{\bfn+\bfbeta^r}} 
	 \Big).
\end{equation}
To establish the linear lower bound we use the elementary, affine
lower bound
\begin{equation}
  \label{eq:G-estimate}
  \forall\, a,b >0, \omega \in \R: \quad
  G(a,b) \geq g(\omega) \, a + g(-\omega)\,  b,
  \quad \text{where } g(\omega) := 1-\ee^{-\omega} + \omega. 
\end{equation}
This estimate follows easily by convexity,
$G(a,b)\geq G(\ee^\omega,1)+\rmD G(\ee^\omega,1)\cdot
(a{-}\ee^\omega,b{-}1)$, and 1-homogeneity giving
$G(\ee^\omega,1)=\rmD G(\ee^\omega,1)\cdot (\ee^\omega,1)$. Note that
equality holds in \eqref{eq:G-estimate} if $\omega = \log(a/b)$.
Moreover, we have $ g(\omega) + g(-\omega)
= 2-\ee^\omega-\ee^{-\omega} \leq 0$, so a
careful choice of $\omega$ depending on $\bfn$ will be necessary to
obtain a good lower bound with a positive leading term.
 
\begin{proposition}
 \label{pr:Psi*V} 
We have the liminf estimate 
\[
\iota_V(\bfu^V) \weaks \varrho \text{ in } \PM(\bfC)
\quad \Longrightarrow \quad 
  \Psi^*_\Lio (\varrho,\rmD\bfE(\varrho)) \leq
  \liminf_{V\to \infty} \Psi^*_V(\bfu^V,\rmD\calE_V(\bfu^V)).
\] 
\end{proposition}
\begin{proof} The special forms of $\bbK(c)$, $E(\bfc)$, and
$\Psi^*_\Lio $ in \eqref{eq:Lio.slope} give the formula
\begin{equation}
\label{eq:PsiLio-G}
 \Psi^*_\Lio (\varrho,\rmD\bfE(\varrho))
  =  \frac12 \int_\bfC
       \sum_{r=1}^R\,\kappa_*^r \,G\Big(
        \frac{\bfc^{\bfalpha^r}}{\bfc_*^{\bfalpha^r}} \, , \, 
       \frac{\bfc^{\bfbeta^r}}{\bfc_*^{\bfbeta^r}}  \Big) \dd \varrho(\bfc)  .  
\end{equation}
Since $\Psi_V^*$ and $\Psi^*_\Lio $ are defined as sums
over $r=1,\ldots,R$ of nonnegative terms, it suffices to show the
result for each $r$ separately, where we suppress
the index $r$. 

Inserting \eqref{eq:G-estimate} into \eqref{eq:Psi*V-all} yields, with
$\omega_\bfn \in \R$ to be fixed afterwards,
\begin{align*}
\Psi_V^*(\bfu^V,\rmD \calE_V(\bfu^V)) 
	& \geq  
	\frac1{2V} \sum_{\bfn\in \calN } \wh\nu^{\bfn}_V
	 \,\Big(g(\omega_\bfn) \frac{u_{\bfn+\bfalpha}^V}{w^V_{\bfn+\bfalpha}}
		  +g(-\omega_\bfn) \frac{u_{\bfn+\bfbeta}^V}{w^V_{\bfn+\bfbeta}} 
	   \Big)
 \\& 
 	= \frac{\kappa_*}{2} \sum_{\bfn\in \calN } 
	 \Big(g(\omega_\bfn)
	 		\hbA^\bfalpha_V(\bfn) u_{\bfn+\bfalpha}^V 
		+ g(-\omega_\bfn) 
			\hbA^\bfbeta_V(\bfn) u_{\bfn+\bfbeta}^V \Big)
\\ &\qquad			 \ \text{ with } 
 	\hbA^\bfdelta_V(\bfn) 
	:= \frac{w_{\bfn}^V}{w_{\bfn+\bfdelta}^V}
	 = \frac{(\bfn {+} \bfdelta)!}{ (\bfc_* V)^\bfdelta \bfn!},
\end{align*}
where we used the detailed-balance conditions from Theorem
\ref{th:DB.CME} for the last identity. Rearranging the sum and
recalling that $\hbA^\delta(\bfn)=0$ for $\bfn \not\in \calN$ we find 
\[
\Psi_V^*(\bfu^V,\rmD \calE_V(\bfu^V)) 
	\geq  \frac{\kappa_*}{2}\sum_{\bfn\in \calN} h^V_\bfn u_\bfn^V
		 \quad \text{with } 
		h^V_\bfn:=
			g(\omega_{\bfn-\bfalpha})
				\hbA^\bfalpha_V(\bfn-\bfalpha)  +
		    g(-\omega_{\bfn-\bfbeta})
		    	\hbA^\bfbeta_V(\bfn-\bfbeta).
\] 

We now choose
$\omega_\bfn=\log\big( \hbA^\bfalpha_V(\bfn) / \hbA^\bfbeta_V(\bfn)
\big)$ for $\bfn\in \calN$ and $\omega_\bfn=0$ otherwise and find, for
all $\bfn$ with $\bfn\geq \bfalpha$ or $\bfn\geq \bfbeta$, the
relation
\begin{align*}
  h^V_\bfn 
  &= G\big(\hbA^\bfalpha_V(\bfn {-} \bfalpha),
    \hbA^\bfbeta_V(\bfn {-} \bfbeta) \big) 
    + f^V_\bfn     \quad \text{with }
  \\
  f^V_\bfn
  &:=   \hbA^\bfalpha_V(\bfn {-} \bfalpha) 
            - \hbA^\bfalpha_V(\bfn {-} \bfbeta) 
            + \hbA^\bfbeta_V(\bfn {-} \bfbeta) 
            - \hbA^\bfbeta_V(\bfn {-} \bfalpha) 
  \\& \quad +  \hbA^\bfalpha_V(\bfn {-} \bfalpha) 
  \log\bigg(
  \frac{\hbA^\bfbeta_V(\bfn {-} \bfbeta)}
  {\hbA^\bfbeta_V(\bfn {-} \bfalpha)}
  \bigg)
  +  \hbA^\bfbeta_V(\bfn {-} \bfbeta) 
  \log\bigg(
  \frac{\hbA^\bfalpha_V(\bfn {-} \bfalpha)}
  {\hbA^\bfalpha_V(\bfn {-} \bfbeta)}
  \bigg)
\end{align*}

The idea is now that as $\frac1V \bfn \to \bfc > \bm0$ we have the
convergences
\[
  \hbA^\bfdelta_V(\bfn {-} \bfalpha) \to  \bfc^\bfdelta/\bfc_*^\bfdelta
  \quad \text{and} \quad 
\hbA^\bfdelta_V(\bfn {-} \bfbeta)  \to \bfc^\bfdelta/\bfc_*^\bfdelta,
\]
which yields $f^V_\bfn \to 0$ and
$h^V_\bfn \to G(\bfc^\bfalpha / \bfc_*^\bfalpha,
\bfc^\bfbeta/\bfc_*^\bfbeta)$ as desired.  To be more precise we
define, for all $\eps \in {]0,1[}$, the functions
\[
G_\eps(a,b)= -\eps + \min\{ (1{-}\eps)G(a,b), 1/\eps\},
\]
which converge monotonely to $G(a,b)$ for $\eps \searrow 0$. A
lengthy calculation using the explicit structure of
$\hbA^\bfdelta(\bfn)$ shows that for all $\eps>0$ there exists
$V_\eps \gg 1$ such that $h^V_\bfn\geq G_\eps
(\hbA^\bfalpha_V(\bfn),\hbA^\bfbeta_V(\bfn))$ for all $V \geq V_\eps$ and
all $\bfn$. 
Even more, if we define the functions $H_V:\bfC\to \R; \ \bfc \mapsto
\sum_{\bfn\in \calN} h_\bfn^V \,\INDIC_{A^V_\bfn}(\bfc)$, then, for all
$\eps>0$ there exists $\wt V_\eps\gg1$ such that 
\[
\forall \,V\geq \wt V_\eps \ \forall \, \bfc\in \bfC:\ 
H_V(\bfc)\geq \mathfrak H_\eps(\bfc):= 
 G_\eps\Big( \frac{\bfc^\bfalpha}{\bfc_*^\bfalpha}, 
\frac{\bfc^\bfbeta}{\bfc_*^\bfbeta} \Big). 
\] 
Hence, using the definition of $\iota_V$ we find the lower bound 
\[
\Psi^*_V(\bfu^V,\rmD\calE_V(\bfu^V)) 
	\geq \frac{\kappa_*}{2}\sum_\bfn h^V_\bfn u_\bfn^V
	   = \frac{\kappa_*}{2}\int_\bfC H_V(\bfc) \dd \iota_V(\bfu^V)(\bfc) 
	\geq \frac{\kappa_*}{2}\int_\bfC \mathfrak H_\eps(\bfc)  \dd \iota_V(\bfu^V)(\bfc) .
\]
Since $\mathfrak H_\eps$ is lower semi-continuous and bounded, this
implies the liminf estimate 
\[ 
	\iota_V(\bfu^V) 
		\weaks \varrho \quad \Longrightarrow \quad 
	\liminf_{V\to \infty} 
		 \Psi^*_V\big(\bfu^V,\rmD\calE_V(\bfu^V)\big) 
	\geq \frac{\kappa_*}{2} \int_\bfC \mathfrak H_\eps(\bfc) \dd \varrho(\bfc).
\]
Because $\eps >0$ was arbitrary we can use the monotone convergence
$\mathfrak H_\eps (\bfc) \nearrow G\big(
\tfrac{\bfc^\bfalpha}{\bfc_*^\bfalpha},  
\tfrac{\bfc^\bfbeta}{\bfc_*^\bfbeta} \big)$ to conclude the desired
result for each of the $R$ reactions 
\begin{align*}
\liminf_{V\to \infty} \Psi^{r,*}_V(\bfu^V,\rmD\calE_V(\bfu^V)) 
& \geq \frac{\kappa_*}{2}  \int_\bfC  G\big(
\tfrac{\bfc^{\bfalpha^r}}{\bfc_*^{\bfalpha^r}},  
\tfrac{\bfc^{\bfbeta^r}}{\bfc_*^{\bfbeta^r}} \big) \dd\varrho(\bfc)
\\& = \frac{\kappa_*}{2}  \int_\bfC \Lambda\big(
\tfrac{\bfc^{\bfalpha^r}}{\bfc_*^{\bfalpha^r}},  
\tfrac{\bfc^{\bfbeta^r}}{\bfc_*^{\bfbeta^r}} \big) \big( \nabla E(\bfc){\cdot}
(\bfalpha^r-\bfbeta^r)\big) ^2 \dd \varrho(\bfc). 
\end{align*}
Summation over $r=1,\ldots,R$ yields the full result for $\Psi^*_V$.
\end{proof}

\subsection{A liminf estimate for the dissipation functional}
\label{su:GamDiss}

In the evolutionary $\Gamma$-convergence method of
\cite{SanSer04GCGF, Serf11GCGF, Miel16EGCG} it is standard to provide
a liminf estimate for the \emph{primal  dissipation potential}
$\Psi_V$ which in our case is defined  via the Legendre transform 
\[
\Psi_V(\bfu,\bfv) 
	= \sup \Bigset{  {\ts\sum\limits_{\bfn\in \calN}} u_\bfn \xi_\bfn 
		- \Psi^*_V(\bfu,\bfxi)}{ \bfxi=(\xi_\bfn)_{\bfn\in \calN} }. 
\]
However, as our theory relies on the dualization $\Psi_V(\bfu,\bfv)
\geq \sum_{\bfn\in \calN} u_\bfn \xi_\bfn - \Psi_V^*(\bfu,\bfxi)$ it
will be sufficient to have the following limsup estimate for
$\Psi_V^*$, which crucially relies on the concavity of the map
$(a,b)\mapsto \Lambda(a,b)$.

\begin{proposition}\label{pr:Psi*Limsup}
Consider any pair $(\varrho,\xi)\in \PM(\bfC) \ti \rmC^1_\rmc(\bfC)$
and set 
$\bfxi^V=\iota_V^* \xi: \calN \to \R$ with $ \iota_V^* $ 
defined in \eqref{eq:EDP-iota*}.
Then, for every family $(\bfu^V)_{V>1}$ we have the limsup estimate 
\begin{equation}
  \label{eq:Psi*V.conv}
  \iota_V(\bfu^V)\weaks \varrho \quad \Longrightarrow \quad 
\limsup_{V\to \infty} \Psi_V^*(\bfu^V,\bfxi^V) \leq 
\Psi_\Lio ^*(\varrho,\xi)=\frac12\int_\bfC \nabla\xi\cdot
\bbK \nabla\xi \dd\varrho(\bfc).
\end{equation}
\end{proposition}
\begin{proof} 
As in the proof of Proposition \ref{pr:Psi*V} we can exploit that $\Psi_V^*$ is a sum of non-negative terms over  $r=1,\ldots, R$. 
Hence, it is sufficient to show the desired limsup  estimate for each reaction individually. 
For notational simplicity we drop the reaction index $r$. 

Defining $\varrho^V=\rho^V\rmd\bfc =\iota_V(\bfu^V)$, relation
\eqref{eq:Psi*Vu-xi.b} leads us to the integral representation
\[
\Psi_V^*(\bfu^V,\bfxi^V)
		= \frac{\kappa_*}2\int_{\bfc\in \bfC} 
			\Lambda\big( \,\rho^{V,a}(\bfc) \,,
						  \rho^{V,b}(\bfc)
			\big) 
			M_V^\xi(\bfc) \dd \bfc, 
\]
where 
\begin{align*}
	\rho^{V,a}(\bfc) = a_V(\bfc) \rho^V(\bfc{+}\tfrac1V\bfalpha), \quad
	\rho^{V,b}(\bfc)= b_V(\bfc) \rho^V(\bfc{+}\tfrac1V\bfbeta),
\end{align*}
and the functions $a_V$, $b_V$, and $M_V^\xi$ are given 
\begin{align*}
	&a_V(\bfc) = \frac{\BBneu V \bfalpha {\bfn} }{V\bfc_*^\bfalpha} ,\quad 
	b_V(\bfc)  = \frac{\BBneu V \bfbeta {\bfn}}{V\bfc_*^\bfbeta},
\quad 
  M_V^\xi(\bfc)
  =  V^2   \big(\bfxi^V_{\bfn + \bfalpha}
      {-}\bfxi^V_{\bfn + \bfbeta}  \big)^2 
  \quad \text{for }  \bfc \in A^V_\bfn.
\end{align*}

Using $\xi \in \rmC_\rmc^1(\bfC)$ there exists $R>0$ such that
$\mathop{\mathrm{sppt}}M_V^\xi \subset \bfC_R:=B_R(0)\cap \bfC$, and
we have uniform convergence
\[
\| a_V{-} a_\infty\|_{\rmL^\infty(\bfC_R)} + 
\| b_V{-} b_\infty\|_{\rmL^\infty(\bfC_R)} + 
\| M_V^\xi - (\nabla \xi \cdot \bfgamma)^2\|_{\rmL^\infty(\bfC_R)} \to
0 \ \text{  as }V\to \infty,
\] 
where $a_\infty(\bfc)=\bfc^\bfalpha/\bfc_*^\bfalpha$,
$b_\infty(\bfc)=\bfc^\bfbeta/\bfc_*^\bfbeta$, and $\bfgamma= \bfalpha{-}\bfbeta$. 
Using $\Lambda(r,t)\leq \frac12(r{+}t)$, the uniform boundedness of $a_V$ and $b_V$ on $\bfC_R$, and that $\varrho_V$ is a probability measure, we see that in the limsup of $\Psi_V^*(\bfu^V,\bfxi^V)$ we can replace
$M_V^\xi$ by $\big((\bfalpha{-}\bfbeta)\cdot\nabla \xi\big)^2$ without changing the limsup in the left-hand side of \eqref{eq:Psi*V.conv}.

Next we consider the functionals $F : \Meas(\bfC_R) \times \Meas(\bfC_R) \to [0, + \infty]$ given by
\begin{align*}
 F(\varrho_1,\varrho_2) 
 = \!\int_{\bfC_R} \!\!
     f(\rho_1(\bfc), \rho_2(\bfc)) 
   \big(\bfgamma\cdot\nabla \xi(\bfc)\big)^2 \dd \bfc 
\ \text{ with } 
 f(r,t) = 
	 \left\{\begin{array}{cl}
		r{+}t{-}\Lambda(r,t),
		 & \text{for }r,t \geq 0,\\
		+ \infty,
		 & \text{else}.\end{array} \right.		 
\end{align*}
Note that $f(r,t) \geq \frac12(r{+}t)$.  Moreover, $f$ is convex and
positively homogeneous of degree $1$.  Thus, $F$ is weak* lower
semi-continuous on $\Meas(\bfC_R)\ti\Meas(\bfC_R)$, cf.\
\cite[Thm.\,6.57]{FonLeo07MMCV}. Now using the convergences
\[
	\varrho^{V,a}
	\weaks a_\infty\varrho\big|_{\bfC_R}  
	\quad \text{ and }  \quad
	\varrho^{V,b}
	\weaks b_\infty\varrho\big|_{\bfC_R} 
	\quad \text{as } V\to \infty, 
\]
we obtain the liminf estimate \
$ \liminf_{V\to \infty} F(\varrho^{V,a},\varrho^{V,b}) \geq F(a_\infty
\varrho,b_\infty\varrho)$.

Thus, in the view of the identity
\begin{align*}
	&\int_{\bfC_R}
		    \Lambda(\rho^{V,a},\rho^{V,b})
			(\bfgamma\cdot\nabla \xi(\bfc))^2
	 \dd \bfc
	= \int_{\bfC_R} 
			(\rho^{V,a}{+} \rho^{V,b}) 
			\big(\bfgamma\cdot\nabla \xi(\bfc)\big)^2 
	  \dd \bfc
	- F(\varrho^{V,a},\varrho^{V,b}),
\end{align*}
and observing that the first term on the right-hand side is weak$^*$
continuous, the limsup for $V\to \infty$ gives
\begin{align*}
\limsup_{V\to \infty} \Psi_V^*(\bfu^V,\bfxi^V) 
	& = \frac{\kappa_*} 2 \limsup_{V \to \infty}  
		\int_{\bfC_R} \Lambda(\rho^{V,a},\rho^{V,b})
			(\bfgamma\cdot\nabla \xi(\bfc))^2
		\dd \bfc 
 \\ & \leq 
 	   \frac{\kappa_*} 2 
	   	\int_{\bfC_R}
			 (a_\infty{+}b_\infty)
			\big(\bfgamma\cdot\nabla\xi)^2 
		\dd \varrho 
		-
		\frac{\kappa_*} 2 
			F(a_\infty\varrho,b_\infty\varrho) 
\\ & 	 = \frac{\kappa_*}2 \!\int_{\bfC_R}\!\!\!\!
\Lambda(a_\infty,b_\infty) \big(\bfgamma {\cdot} \nabla \xi(\bfc)\big)^2
\dd \varrho(\bfc) =\Psi_\Lio ^*(\varrho,\xi). 
\end{align*}  
This is the desired result for one reaction, and the full result
follows by summation over $r=1,\ldots,R$ and the definition of $\bbK$, namely 
$\nabla \xi\cdot \bbK\nabla \xi= \sum_{r=1}^R
\kappa_*^r \Lambda\big( \frac{\bfc^{\bfalpha^r}}{\bfc_*^{\bfalpha^r}},   
\frac{\bfc^{\bfbeta^r}}{\bfc_*^{\bfbeta^r}} \big)(\bfgamma^r{\cdot}
\nabla \xi)^2$.    
\end{proof}

\subsection{Convergence of solutions}
\label{su:GamCME2Lio}

Here we provide the general convergence result as $V\to \infty$ for
the appropriately embedded solutions
$\bfu^V:{[0,\infty[} \to \PM(\calN)$ of the CME to the solutions
$\varrho:{[0,\infty[} \to \PM(\bfC)$ of the Liouville equation, which
is a simple transport along the solutions of the RRE
$\dot\bfc = -\bfR(\bfc) = -\bbK(\bfc) \rmD E(\bfc)$.  Our approach
follows the strategy of evolutionary $\Gamma$-convergence as initiated
in \cite{SanSer04GCGF,Serf11GCGF} with the new idea of dualization as
introduced in \cite{LMPR17MOGG}.

\begin{theorem}[Evolutionary $\Gamma$-convergence of CME to Liouville]
\label{th:CME2Lio} For all $V>1$ consider a solution $\bfu^V:{[0,\infty[}
\to \PM(\calN)$ of the CME \eqref{eq:CME.1}. 
Assume that the initial conditions are well-prepared in the sense that 
\[
	\iota_V(\bfu^V(0)) \weaks \varrho^0 \text{ in }\PM(\bfC)
			 \quad \text{and} \quad 
	\calE_V(\bfu^V(0)) \to \bfE(\varrho^0).
\]
Then, for all $t > 0$, we have the convergence
\[
	\iota_V(\bfu^V(t)) \weaks \varrho(t) \text{ in }\PM(\bfC) 
		\quad \text{and} \quad 
	\calE_V(\bfu^V(t)) \to \bfE(\varrho(t)),
\]
where $\varrho:{[0,\infty[} \to \PM(\bfC)$ is the unique solution of
the Liouville equation \eqref{eq:KE-gradflow} starting at 
$\varrho(0)=\varrho^0$, i.e.,\ for all $\varphi\in
\rmC^1_\rmc([0,T]\ti \bfC)$ with $\varphi(T,\cdot) = 0$ we have 
\begin{equation}
  \label{eq:Lio.weak}
  \int_\bfC \varphi(0,\bfc) \varrho^0(\rmd\bfc) + \int_0^T\int_\bfC
  \Big(\pl_t \varphi(t,\bfc) - \nabla\varphi(t,\bfc){\cdot} \bbK(\bfc)
  \nabla E(\bfc) \Big) \varrho(t,\rmd\bfc) \dd t =0. 
\end{equation}
Moreover, for all $r,s\in [0,T]$ with $r<s$ we have the energy identity
\begin{equation}
  \label{eq:Lio.EnBal}
  \bfE(\varrho(s)) 
  	+ 2 \int_r^s \!
		 \Psi_\Lio^*\big(\varrho(t), - \rmD\bfE(\varrho(t))\big)
		 \dd t 
	= \bfE(\varrho(r)).
\end{equation}
\end{theorem}

For the proof we use the energy-dissipation principle for $V \geq 1$ and pass to the limit in each of the terms. 
If $\bfu^V$ is a solution of the CME, then for all $T>0$ we have 
\begin{equation}
  \label{eq:EDP-V}
  \calE_V(\bfu^V(T)) 
  + \int_0^T \Psi_V(\bfu^V,\dot \bfu^V) 
  + \Psi_V^*\big(\bfu^V,{-}\rmD \calE_V(\bfu^V)\big) \dd t 
  = \calE_V(\bfu^V(0)). 
\end{equation}
Following the ideas in \cite{DisLie15GSMC} for the passage from a
Markov chain to the Fokker--Planck equation or the general methods in
evolutionary $\Gamma$-convergence, we want to pass to the limit in
each of the four terms.  As a general fact, it will be sufficient to
obtain liminf estimates on the left-hand side, since by a chain-rule
argument an estimate with ``$\leq$'' instead of equality can be turned
back into an equality.  Moreover, by the assumptions of the theorem we
see that the right-hand side converges to the desired limit.

However, it is rather delicate to pass to the limit in the integral
$\int_0^T \Psi_V(\bfu^V,\dot\bfu^V) \dd t$, because the potential
$\Psi_V$ is only implicitly defined and we expect the limit to be
given in terms of the Benamou-Brenier formula for the Wasserstein
distance induced by the metric on $(\bfC,\bbK)$.  A major difficulty
is even to obtain a suitable equi-continuity for the solutions $\bfu^V$
to be able to extract a subsequence converging at all times. In
particular, it is unclear how to pass to the limit in
$\iota_V(\dot \bfu^V(t))$ by a direct argument.

Hence, following \cite{LMPR17MOGG}, we estimate the primal dissipation
potential $\Psi_V$ from below using the definition in terms of the
Legendre transform of $\Psi_V^*$.  Using additionally an integration
by parts we have 
\begin{align*}
&\int_0^T \Psi_V(\bfu^V,\dot \bfu^V)\dd t \geq \mathfrak J_V 
(\bfu^V,\bfeta) \quad \text{ for all } \bfeta \in
\rmC^1([0,T];\ell^\infty(\calN))  \text{ with}
\\ 
&\mathfrak J_V(\bfu,\bfeta):=  \langle \bfu(T), \bfeta(T)\rangle -
\langle \bfu(0), 
\bfeta(0)\rangle - \int_0^T  \langle \bfu(t),\dot\bfeta(t) \rangle
+ \Psi_V^*(\bfu(t),\bfeta(t)) \dd t,
\end{align*} where $\langle\bfu,\bfeta\rangle := \sum_{n\in \calN}
u_\bfn \eta_\bfn$.  With this argument we can replace the
energy-dissipation principle \eqref{eq:EDP-V} by the estimate
\begin{equation}
  \label{eq:EDP-V2} \calE_V(\bfu^V(T)) + \mathfrak
J_V(\bfu^V, \bfeta ) + \int_0^T \Psi_V^*\big(\bfu^V,
{-}\rmD\calE_V(\bfu^V)\big) \dd t \leq \calE_V(\bfu^V(0)),
\end{equation} which holds for all differentiable $\bfeta$. In this
equation we are then able to pass to the limit $V\to \infty$, when choosing
$\bfeta= \bfeta^V = \iota^*_V(\xi)$ for a smooth function $\xi$. 

At the end we are then able to calculate the supremum over all $\xi$
by using the especially simple quadratic structure in $\xi$, which
mirrors the fact that the Liouville equation is a simple transport
equation. \bigskip

\noindent
\begin{proof}[Proof of Theorem \ref{th:CME2Lio}]
\vspace{-1em}

\paragraph*{Step 1: Embedding and uniform a priori bounds.} 
We now consider the family $\bfu^V:[0,T] \to \PM(\calN)$ and embed it
into $\PM(\bfC)$ via $\iota_V$ from \eqref{eq:iota-V}. As in
\cite{DisLie15GSMC} we show an equi-continuity in a 1-Wasserstein
distance, but introduce an additional weight accounting for our
unbounded domain $\bfC$.  We define the maximal order $p$ of all
reactions via
\[
p:= \max\set{|\bfalpha^r|_1,|\bfbeta^r|_1}{ r=1,\ldots, R}.
\]
For $\mu \in \Meas(\bfC)$ and for $\varrho_0, \varrho_1 \in \PM(\bfC)$  we set  
\[
  \| \mu \|_{1\rmW}
  := \sup\biggset{ \int_\bfC f(\bfc) \dd \mu(\bfc) }{f \in \bbF}
\quad  \text{and} \quad 
	d_{1\rmW}(\varrho_0,\varrho_1) 
		=  \| \varrho_0 - \varrho_1 \|_{1\rmW},
\]
where $\bbF:= \set{f\in \rmC^1(\bfC)}{ \sup_\bfC
  (1{+}|\bfc|^p) |\nabla f(\bfc)| \leq 1}$.  

Using the definition of the Markov generators $\calQ^r_V$ in terms of
the coefficients $\BBneu V{\bfdelta^r}\bfn$, see \eqref{eq:calL}, it is easy
to derive the uniform estimate $\|\iota_V(\dot\bfu^V(t))\|_{1\rmW}\leq
C_{1\rmW}$ independently of the initial conditions and $V\geq 1$ (one simply
needs $\sum u^V_\bfn\equiv 1$). Hence, we obtain the uniform Lipschitz
bound 
\[
d_{1\rmW}\big(\iota_V(\bfu^V(t)),\iota_V(\bfu^V(s))\big)\leq C_{1\rmW}|t{-}s|\quad
\text{ for all }s,t\in [0,T] \text{ and all }V\geq 1. 
\]
Moreover, as
$\calE_V(\bfu^V(t)) \leq \calE_V(\bfu^V(0)) \leq
\bfE(\varrho^0)+o(1)_{V\to \infty}$ by well-preparedness, the
equi-coercivity of $\calE_V$ established in \eqref{eq:EV-equicoerc}
yields the uniform bound
\begin{equation}
  \label{eq:EDP-V3}
  \exists\, V_* \geq 1, C_\rmB < \infty \ \forall\, t>0,\ V\geq V_*:  \int_\bfC (1{+}|\bfc|)
  \iota_V(\bfu^V(t)) \dd \bfc \leq C_\rmB. 
\end{equation}

\paragraph*{Step 2: Extraction of a subsequence.}
The subset of $\PM(\bfC)$ defined by the boundedness of the above
first moment is a compact subset of the metric space
$(\PM(\bfC),d_{1\rmW})$. Indeed, using Prokhorov's theorem one
finds that this set is weak$^*$ sequentially compact.  Since
$d_{1\rmW}$ is dominated by the bounded Lipschitz metric (which
metrizes weak$^*$ convergence), the compactness of
$(\PM(\bfC),d_{1\rmW})$ follows. 

Hence, we can apply the abstract Arzel\`a-Ascoli theorem in
$(\PM(\bfC),d_{1\rmW})$ to extract a subsequence $V_k\to \infty$ and a
limit function $\varrho:[0,T]\to \PM(\bfC)$ such that
\begin{subequations}\label{eq:EDP.conv}
\begin{align}
\label{eq:EDP.conv.a}
&  \forall\, t\in [0,T]:
&& \iota_V(\bfu^V(t)) \weaks \varrho(t) \text{ in } \PM(\bfC),\\
\label{eq:EDP.conv.b}
&  \forall\, s,t\in [0,T]: 
&& d_{1\rmW}(\varrho(t),\varrho(s)) \leq C_{1\rmW}|t{-}s|,\\
\label{eq:EDP.conv.c}
& \forall\, t\in[0,T]: 
&& \bfE(\varrho(t)) \leq \bfE(\varrho^0),\\
\label{eq:EDP.conv.d}  
& 
&& \text{the mapping } t \mapsto \varrho(t) \text{ is weak* continuous}. 
\end{align}
\end{subequations}
At first, in place of \eqref{eq:EDP.conv.a} one obtains
$d_{1\rmW}\big(\iota_V(\bfu^V(t)),\varrho(t)\big)\to 0$.  To derive
\eqref{eq:EDP.conv.a}, we use the bound \eqref{eq:EDP-V3} together
with the fact that any bounded continuous function can be uniformly
approximated on compact sets by (multiples of) functions in $\bbF$.
Similarly, \eqref{eq:EDP.conv.d} follows from
\eqref{eq:EDP.conv.b}. In particular, combining
\eqref{eq:EDP.conv.d} and the assumption
$\iota_V(\bfu^V(0)) \weaks \varrho^0$ we conclude
$\varrho(0)=\varrho^0$.  Finally, \eqref{eq:EDP.conv.c} follows via
\eqref{eq:EDP.conv.a} from Theorem \ref{th:EV-Gamma}:
\[
\bfE(\varrho(t)) \leq \liminf_{V\to \infty}
\calE_V\big(\bfu^V(t)\big) \leq \liminf_{V\to \infty} 
\calE_V\big(\bfu^V(0)\big) = \bfE(\varrho^0).
\]

\paragraph*{Step 3: Limit passage in \eqref{eq:EDP-V2}.} Combining
\eqref{eq:EDP.conv.a} for $t=T$ and Theorem \ref{th:EV-Gamma} (cf.\
\eqref{eq:EV-liminf}), the first term satisfies the liminf estimate 
$\liminf_{V\to \infty} \calE_V(\bfu^V(T)) \geq \bfE(\varrho(T))$. 
For the last term we use the assumption $ \calE_V(\bfu^V(0)) \to
\bfE(\varrho^0)=\bfE(\varrho(0))$.  

For the third term we employ Proposition \ref{pr:Psi*V} for each $t\in
[0,T]$ based on \eqref{eq:EDP.conv.a}. Using Fatou's lemma we conclude
the liminf estimate 
\begin{align*}
\liminf_{V\to \infty} \int_0^T \Psi_V^*\big(\bfu^V(t),
   {-}\rmD\calE_V(\bfu^V(t))\big) \dd t 
   &\geq \int_0^T \liminf_{V\to
   \infty}\Psi_V^*\big(\bfu^V(t),{-}\rmD\calE_V(\bfu^V(t))\big) \dd t 
\\ &\geq \int_0^T \Psi^*_\Lio \big(\varrho(t),{-}\rmD\bfE(\varrho(t))\big)\dd t.
\end{align*}

Thus, it remains to pass to the limit in $\mathfrak
J_V(\bfu^V,\eta)$. For this we choose an arbitrary
$\xi\in \rmC^1_\rmc([0,T]\ti \bfC)$ and define
$\bfxi^V(t)=\iota_V^*(\xi(t))$, cf.\ \eqref{eq:EDP-iota*}. With this
choice we can apply Proposition \ref{pr:Psi*Limsup} for all
$t\in [0,T]$ based on \eqref{eq:EDP.conv.a}. Now, Fatou's lemma yields
\begin{align*}
	& \liminf_{V\to \infty} \mathfrak J_V(\bfu^V ,\bfxi^V) 
		\geq 
			\mathfrak J_\Lio  (\varrho,\xi)
	\quad \text{where} \\ 
	& \mathfrak J_\Lio(\varrho,\xi) 
		:=  \int_{\bfC} \xi(T,\bfc) \varrho(T,\rmd \bfc)
		   - \int_\bfC \xi(0,\bfc)\varrho(0,\rmd \bfc)\\
&\phantom{\mathfrak J_\Lio  (\varrho,\xi):= }  - \int_0^T\int_\bfC
 \Big( \pl_t\xi(t,\bfc)+ \frac12\nabla \xi(t,\bfc)\cdot \bbK(\bfc)\nabla
\xi(t,\bfc) \Big) \varrho(t,\rmd \bfc) \dd t.
\end{align*}  
In summary, we conclude that the limit function $\varrho:[0,T]\to
\PM(\bfC)$ satisfies 
\begin{equation}
  \label{eq:EDP-Lio.xi}
  \bfE(\varrho(T)) + \mathfrak J_\Lio  (\varrho,\xi) + \int_0^T
 \!\! \Psi^*_\Lio 
 \big(\varrho(t),{-}\rmD\bfE(\varrho(t))\big) \dd t \leq \bfE(\varrho(0))
\end{equation}
for all $\xi\in \rmC^1_\rmc([0,T]\ti \bfC)$.

\paragraph*{Step 4: Energy balance.} 
By inserting $\xi\equiv 0$ in \eqref{eq:EDP-Lio.xi} we obtain the
upper bound
\[
\calD(\varrho;0,T):= \int_0^T\!\!\int_\bfC \nabla E(\bfc){\cdot}
\bbK(\bfc) \nabla E (\bfc) \varrho(t,\rmd\bfc) \dd t \leq 2\big(
\bfE(\varrho(0)) - \bfE(\varrho(T))\big) .
\]
We want to show energy balance, i.e.,\ equality when the factor $2$ is
omitted.  For this purpose, we
observe that the measures $\varrho(t,\cdot)\in \PM(\bfC)$ decay at
infinity such that \eqref{eq:EDP.conv.c} holds. Hence, we may also use
$\xi(t,\bfc)=\lambda E(\bfc)$ as testfunctions in
\eqref{eq:EDP-Lio.xi}. Writing shortly $e(t):=\bfE(\varrho(t))$ we
find $\mathfrak J_\Lio  (\varrho,\lambda E )= \lambda \big(e(T) - e(0)\big) -
\frac{\lambda^2}2 \calD(\varrho;0,T)$ and obtain 
\[
  -\lambda\big(e(0)- e(T)\big) - \frac{\lambda^2}2 \calD(\varrho;0,T)
  = \mathfrak J_\Lio (\varrho,\lambda E ) \leq e(0)- e(T)
  -\frac12\calD(\varrho;0,T) \quad \text{for all } \lambda \in \R.
\]
Maximizing with respect to $\lambda$ leads to $(e(0){-} e(T))^2/\calD
\leq 2(e(0){-} e(T)) - \calD$ which implies $e(0){-} e(T)=\calD$, or
more explicitly $ \calD(\varrho;0,T)= \bfE(\varrho(0)) -
\bfE(\varrho(T))$, which is 
the desired energy balance \eqref{eq:Lio.EnBal} for $r=0$ and $s=T$. 

Moreover, we can repeat the calculation on $[0,s]$ with $0 < s < T$
instead of $[0,T]$. The full result \eqref{eq:Lio.EnBal}
follows by subtracting the identity on $[0,r]$ from that on $[0,s]$.

\paragraph*{Step 5: Weak form of gradient flow equation.} With Step 4
we rewrite \eqref{eq:EDP-Lio.xi} as
\[
\mathfrak J_\Lio(\varrho,\xi) \leq
 \bfE(\varrho(0))- \bfE(\varrho(T)) -
 \frac12\calD(\varrho;0,T)=\frac12\calD(\varrho;0,T), 
\]
and know that the left-hand side is maximized by
$\xi:(t,\bfc)\mapsto - E(\bfc)$. Inserting the test functions
$\xi(t,\bfc) = \delta \varphi(t,\bfc) - E(\bfc)$  with small
$\delta>0$ and $\varphi \in \rmC^1_\rmc([0,T]\ti \bfC)$ with 
$\varphi(T,\cdot)=0$ we arrive, after some cancellations and after
dividing by $\delta >0$, at
\[
-\int_\bfC \varphi(0,\bfc) \varrho(0,\dd \bfc)-
\int_0^T \!\int_\bfC \Big(\pl_t\varphi -\nabla\varphi \cdot\bbK
\nabla \big( E{-}\tfrac\delta2\varphi\big)
\Big)\varrho(t,\rmd\bfc) \dd t \leq 0.
\]
Taking the limit $\delta \searrow 0$ and replacing $\varphi$ by
$-\varphi$, we obtain the desired result \eqref{eq:Lio.weak}. 

With this, Theorem \ref{th:CME2Lio} is established. 
\end{proof}

\section{Approximation via Fokker--Planck equations}
 \label{su:FokkerPlanck}

 In the above section we have seen that the Liouville equation is the
 proper limit of the CME for $V \to \infty$. However, for finite but
 large $V$ it can still be advantageous to replace the discrete CME by
 a continuous PDE with $V$ as a large parameter. In this range the
 stochastic modeling is done by the so-called \emph{Langevin
   dynamics}, see \cite{Kurt78SATD, Gill00CLE, WinSch17HMCR}, which is
 based on a stochastic perturbation of the reaction-rate equation
 \eqref{eq:I.RRE}, see \eqref{eq:RRE.SDE}.  At the level of
 probability distributions the corresponding model is the associated
 Fokker--Planck equation (FPE).  We will discuss two different gradient
 flow approximations: in the first we simply add a suitable ``entropic
 term'' to the driving functional, but keep the dissipation fixed
 (cf.\ Section \ref{su:SimpleFP}), while in the second we expand
 $\calE_V$ and $\calK_V$ such that all terms of order $1/V$ are
 correct (cf.\ Section \ref{su:ImprovedFP}).

\subsection{Improved approximation of the relative entropy}
\label{su:ImprovedEntropy}

We interpret the sum in the definition of $\calE_V$ as a Riemann sum
and replace it by a corresponding integral. The main point of the
improvement is that we keep the entropy term
$\frac1V \sum u_\bfn \log u_\bfn$ in the definition of
$\calE_V(\bfu)$, which is in contrast to the limit $\bfE$ obtained in
Theorem \ref{th:EV-Gamma}.  Working with absolutely continuous
probability measures $\varrho(\rmd \bfc) = \rho(\bfc) \dd \bfc$ with
$\rho \in \rmL^1(\bfC)$, we can define the $V$-dependent entropy by
\begin{equation}
  \label{eq:4.EV.neu}
  \bfE_V(\varrho) = \frac1V \int_\bfC \rho(\bfc)
\log\Big(\frac{\rho(\bfc)}{W_V(\bfc)}\Big) \dd \bfc,
\end{equation}
where the equilibrium density $W_V\in \rmL^1(\bfC)$ has to be chosen
suitably. A first simple approximation is
$\wt W_V(\bfc)=\frac1{\wt Z_V} \,\ee^{-V E(\bfc)} $ with
$E(\bfc)=\sum_{i=1}^I c^*_i\LB(c_i/c^*_i) $ as above and
$\wt Z_V=\int_\bfC \ee^{-V E(\bfc)} \dd \bfc$. However, a better and
more refined $W_V$ is obtained using the next order of expansion in
Stirling's formula \eqref{eq:Stirling} as well. For this we use the
approximation $k_n \approx n+1/6$, i.e.,\
$\log(n!) = n\log n - n +\frac12 \log\big(2\pi(n{+}\frac16)\big) +
O(1/n^2)$ for $n\to \infty$.  Hence, taking the limits
$V, \:|\bfn| \to \infty$ such that $\frac{\bfn}{V} \to \bfc$, we
obtain
\begin{align*}
  - \frac{1}{V} \log w_\bfn^V
  \approx E(\bfc) + \frac1V G_V(\bfc)
\end{align*}
with the $V$-dependent correction $G_V(\bfc) := \frac12\sum_{i=1}^I
\log \big( 2\pi(V c_i + \tfrac16) \big)$ for $E$.

We now take a probability measure
$\varrho = \rho \rmd\bfc \in \PM(\bfC)$ and a discrete approximation
$\bfu \approx \varkappa_V (\varrho) \in \PM(\calN)$, where
$\varkappa_V : \PM(\bfC) \to \PM(\calN)$ is the natural projection
defined in \eqref{eq:kappa-V}. Then the Riemann-sum approximation
results in
\begin{align*}
\calE_V(\bfu) 
  & =  \frac{1}{V}\sum_{\bfn \in \calN} u_\bfn \log u_\bfn 
  	 - \frac{1}{V}\sum_{\bfn \in \calN} u_\bfn \log w_\bfn^V
\\& \approx \frac{1}{V} 
 	\int_{\bfC} \rho(\bfc) \log \rho(\bfc) \dd \bfc
	 - I \frac{\log V}{V}
	 + \int_{\bfC}\Big(  E(\bfc) + \frac1V G_V(\bfc)\Big) \rho(\bfc) \dd \bfc
\\& = \frac1V \int_\bfC \rho(\bfc)
\log\Big(\frac{\rho(\bfc)}{\widehat W(V,\bfc, \bfc^*)}\Big) \dd \bfc,
\\
\text{where }\ &
\widehat W(V,\bfc, \bfc_*)
   = \prod_{i=1}^I \widehat\sfW(V,c_i,c_i^*) \  \text{ with }
  \widehat\sfW(V,c,c^*)
   = \frac{V \ee^{-V c^*\LB(c/c^*)}}{
   			\sqrt{2\pi(V c + 1/6) }}.
\end{align*}
The probability density $W_V$ is then defined by normalizing
$\widehat W(V, \cdot, \bfc_*)$.  We thus set
$\sfZ(V, c^*) := \int_0^\infty \widehat\sfW(V, c, c^*) \dd c$ and
\begin{equation}
  \label{eq:sfW}
W_V (\bfc) := \prod_{i=1}^I \sfW(V,c_i,c_i^*) \quad \text{ with }
\sfW(V,c,c^*) :=
\frac{\widehat\sfW(V,c,c^*)}{\sfZ(V, c^*)}.
\end{equation}
This yields the expansion 
\[
-\frac1V \log W_V(\bfc) 
= E(\bfc) + \frac1V E^V_1(\bfc)
\quad \text{where} \quad E^V_1(\bfc) =  
\wh{\mathsf z}(V,\bfc_*) + 
\frac12\sum_{i=1}^I \log\big(V c_i{+}\tfrac16 \big)
\]
with $\wh{\mathsf z}(V,\bfc_*) 
	  := \sum_{i=1}^I
		\log\big(\sqrt{2\pi}\,\sfZ(V,c_i^*)/V\big)$.  
In summary, for $\bfE_V$ defined via \eqref{eq:4.EV.neu} and
\eqref{eq:sfW}  we have
\begin{equation}
  \label{eq:whE.V}
  	\bfE_V(\varrho) = \bfE(\varrho)+ \frac1V \int_\bfC \big(\rho \log
  \rho +E^V_1 \rho \big) \dd \bfc , 
\end{equation}
and
$\rmD \bfE_V(\varrho)(\bfc)=\frac1V \log \rho(\bfc) -\frac1V \log
W_V(\bfc)$.

\subsection{Simple Fokker--Planck approximation}
\label{su:SimpleFP}

Here we keep the $V$-independent Onsager operator
$\bfK(\varrho):\xi\mapsto -\div\big(\varrho \bbK\nabla \xi\big)$ of
the Liouville equation and obtain the $V$-dependent continuous
gradient system $(\PM(\bfC), \bfE_V,\bfK)$.  The associated
gradient-flow equation
$\dot\varrho = - \bfK(\varrho)\rmD \bfE_V(\varrho)$ is the FPE
\begin{equation}
  \label{eq:FP-V}
 \dot\rho
    = \div\Big( \frac1V
\bbK(\bfc) \nabla \rho + \rho \bfR(\bfc)+
\rho \bfA_V(\bfc)\Big), 
\end{equation}
where we used $\bbK(\bfc) \rmD E(\bfc)=\bfR(\bfc)$ and set
$\bfA_V(\bfc):= \frac12\bbK(\bfc) \big(
\frac1{Vc_i{+}1/6}\big)_{i=1,\ldots,I}$. 

We expect that this FPE is a good approximation
to the CME for all sufficiently large $V$. 
In particular, \eqref{eq:FP-V} has the steady state $\rho= W_V$,
which is close to the discrete steady state
$\bfw^V\in \PM(\calN)$ using the embedding as above. In contrast, the only steady states of the Liouville equation \eqref{eq:Liouville} are concentrated on the equilibria of $\dot \bfc=-\bfR(\bfc)$. 
Of course, the FPE still respects the invariant
sets $\bfI(\bfq)$, because the mobility $\bbK$ of the Onsager operator $\bfK$ is the same as for the Liouville equation. 
In particular, $\rho = W_V$ is the unique equilibrium density if and only if $\bbK$ has full rank, i.e.,\ $\bfI(\bfq)=\bfC$ for all $\bfq\in \mfQ$.

The simpler choice $\wt W_V(\bfc)=\frac1{\wt Z(V)} \ee^{-V E(\bfc)} $ for the equilibrium yields the relative entropy
\[
\wt\bfE_V (\varrho)
	= \frac1V 
		\int_\bfC 
			\rho(\bfc) 
			\log\Big(\frac{\rho(\bfc)}{\wt W_V(\bfc)}\Big)
		 \dd \bfc
 =
	  \int_\bfC \!\! \Big(\frac1V \log \rho(\bfc)  
  {+} E(\bfc) \Big) \rho(\bfc) \dd \bfc + \frac{\log(\wt Z(V))}{V}. 
\]
The flow equation $\dot\varrho = - \bfK(\varrho)\rmD \wt\bfE_V(\varrho)$ induced by the gradient system $(\PM(\bfC),\wt\bfE_V,\bfK)$ is the simplified FPE
\begin{equation}
  \label{eq:FPsimple}
  \dot \rho = \div\Big( \frac1V
\bbK(\bfc) \nabla \rho + \rho \bfR(\bfc) \Big), 
\end{equation}
which is the same as \eqref{eq:FP-V} but with $\bfA_V\equiv 0$. The
simplified equation will be used below as well, since $\wt W_V$ has a
simpler explicit form.

We believe that this approximation is suitable for many purposes.
However, it does not produce the correct diffusion as derived in
\cite[Eqn.\,(1.7)]{Kurt78SATD}. This diffusion correction is used to
improve the RRE $\dot\bfc = - \bfR(\bfc)$ by replacing it by a
stochastic differential equation called the chemical Langevin
equations (CLE) in \cite{Gill00CLE,WinSch17HMCR}, see \eqref{eq:RRE.SDE}.
The associated Fokker--Planck equation  takes the
form
\begin{equation}
  \label{eq:FP.CLE}
  \dot \rho = 
  \frac1V \sum_{i,j=1}^I \partial_{ij}^2
  \big( \rho \wh\bbK_{\CLE}(\bfc)_{ij} \big)
  + \div\big( \rho
  \bfR(\bfc)\big) , 
\end{equation}
where $\wh\bbK_\CLE(\bfc) \in \R^{I\ti I}$ is given in
\eqref{eq:whK.CLE} and differs from $\bbK$ as the logarithmic mean
$\Lambda(a,b)$ between $a=\bfc^{\bfalpha^r}\!/\bfc^{\bfalpha^r}_*$ and
$b=\bfc^{\bfbeta^r}\!/\bfc^{\bfbeta^r}_*$ is replaced by the
arithmetic mean $\frac12(a{+}b)$.  Obviously, \eqref{eq:FP.CLE}
does not have a gradient structure with respect to $\wh\bbK_\CLE$,
because there is no function $\bfc \mapsto \wh E(\bfc)$ such that
$\bfR(\bfc)=\wh\bbK_\CLE(\bfc)\nabla \wh E(\bfc)$.

\subsection{Fokker--Planck equation with higher-order terms}
\label{su:ImprovedFP}

To derive a proper expansion for the term of order $1/V$ in the
evolution equation, we work with the $V$-dependent entropy $\bfE_V$
defined in Section \ref{su:ImprovedEntropy}.  Up to an irrelevant
$V$-dependent constant, this functional approximates
$\calE_V$ from \eqref{eq:CME.E.K} up to order $1/V^2$.

Similarly, we need to derive a suitable expansion for the dissipation
potential, which can be done for each reaction independently.  The
discrete dual dissipation potential is given by
\eqref{eq:Psi*Vu-xi.b}, namely
\[
\Psi_V^*(\bfu,\bfxi)= \frac V2 \sum_{\bfn \in \calN} \Lambda\big( 
	k_\FW \BBneu V\bfalpha{\bfn} u_{\bfn+\bfalpha}, 
	k_\BW \BBneu V \bfbeta{\bfn} u_{\bfn+\bfbeta} \big) \big( \mu_{\bfn{+}\bfalpha}{-}
\mu_{\bfn+\bfbeta}\big)^2.
\]
For a smooth function $\xi : \bfC \to \R$ we use the second-order
accurate midpoint approximation
$\bfmu = \wh\bfmu^\xi_V: \bfn\mapsto
\xi\big(\frac1V(\bfn{+}\bfdelta)\big)$ with
$\bfdelta = \frac12(1,\ldots,1)$ to obtain the expansion
\[
	V\big( \mu_{\bfn+\bfalpha} - \mu_{\bfn+\bfbeta}\big) 
		= \nabla \xi(\sfc_\bfn^V)\cdot\big(\bfalpha{-}\bfbeta)
			 + O(1/V^2)_{V \to \infty} 
			 \ \text{ with }
	\sfc^V_\bfn:=\frac1{V} \Big(\bfn{+}\frac{\bfalpha{+}\bfbeta}{2}{+}\bfdelta\Big),
\]
where we used symmetric difference quotients to obtain second order
accuracy.  Moreover, for a smooth and sufficiently fast decaying
$\varrho=\rho\rmd\bfc \in \PM(\bfC)$ we define the associated discrete $\bfu\in
\PM(\calN)$ via $\bfu = \varkappa_V(\varrho) = \iota^*_V\varrho$, which yields $V^Iu_\bfn=\rho
\big(\frac1V(\bfn{+}\bfdelta)\big) + O(1/V^2)$,
\begin{align*}
  &V^I u_{\bfn+\bfalpha} = \rho(\sfc^V_\bfn) + \frac1{2V} \nabla
  \rho(\sfc^V_\bfn) \cdot( \bfalpha{-}\bfbeta)+ O(1/V^2),
\end{align*}
and similarly for $V^I u_{\bfn+\bfbeta}$. 
Hence, for the arguments of $\Lambda$ we find the expansion
\begin{align*}
	& \frac{1}{V} \BBneu V\bfalpha {\bfn} 
		  V^{I} u_{\bfn+\bfalpha} 
	= (\sfc^V_\bfn)^\bfalpha \rho(\sfc^V_\bfn) 
		+ \frac1V  F^V_\bfn + O(1/V^2) \\ 
	&\text{ with } 
	 F^V_\bfn 
	= - (\sfc_\bfn^V)^\bfalpha \rho(\sfc_\bfn^V) \sum_{i=1}^I 
\frac{\alpha_i\beta_i}{2(\sfc_\bfn^V)_i} 
+ \frac12 (\sfc_\bfn^V)^\bfalpha \nabla\rho(\sfc_\bfn^V)\cdot 
    (\bfalpha{-}\bfbeta). 
\end{align*}
For all smooth functions 
$f, g : \bfC \to \R$ with compact support in int$(\bfC)$, the
trapezoidal rule for Riemann integrals gives 
\[
	\sum_{\bfn\in \calN}
		\Big(f(\sfc^V_\bfn)
			+ \frac1V g(\sfc^V_\bfn) \Big) \frac1{V^I} 
	= \int_\bfC \Big(f(\bfc) + \frac1V g(\bfc) \Big) \dd\bfc 
			+ O(1/V^2).
\]
Hence, for smooth $\rho$ and $\xi$ we find the expansion 
\begin{align*}
  & \Psi_V^{*}\big(\varkappa_V (\varrho),\wh\bfmu^\xi_V\big)
  = \Phi^*_V(\varrho,\xi) + O(1/V^2) 
     \text{ for }V\to \infty \text{ with}\\
  & \Phi_V^{*}(\varrho,\xi) 
  	= \frac12 \int_\bfC \!\Big( \Lambda (k_\FW \bfc^\bfalpha,
  k_\BW\bfc^\bfbeta) \rho(\bfc) + \frac1V \Upsilon
  \big(\bfc,\rho(\bfc),\nabla\rho(\bfc)\big)\Big)
  \big(\nabla\xi(\bfc)\cdot(\bfalpha{-}\bfbeta)\big)^2 \dd\bfc,
\end{align*}
where the correction term $\Upsilon$ takes the explicit form 
\begin{align*}
& \Upsilon(\bfc,\rho,\bfp) 
	= \Upsilon_0(\bfc) \rho + \Upsilon_1 (\bfc)\:\bfp
     \vdot(\bfalpha {-}\bfbeta) \quad \text{with}\\
& \Upsilon_0(\bfc) 
	=  - \frac12 \Lambda(k_\FW\bfc^\bfalpha, 
					  k_\BW \bfc^\bfbeta) \,
	\bfalpha{\vdot} \check\bfC \bfbeta 
	\ \text{ with }
	\check\bfC = \mathrm{diag}(c_i^{-1})_{i=1,\ldots,I}, \\ 
& \Upsilon_1(\bfc) 
	= \Lambda(k_\FW\bfc^\bfalpha, k_\BW \bfc^\bfbeta) \,
	\frac{ k_\FW\bfc^\bfalpha {+} k_\BW \bfc^\bfbeta 
	 {-} 2\Lambda(k_\FW\bfc^\bfalpha, k_\BW \bfc^\bfbeta) }
	{2(k_\FW\bfc^\bfalpha {-} k_\BW \bfc^\bfbeta) }. 
\end{align*}
Here we used the relation $\pl_a \Lambda(a,b)= \frac{\Lambda(a,b)}{a}\: 
\frac{a- \Lambda(a,b)}{a-b}$, giving
\[
a\pl_a \Lambda(a,b) {+} b\pl_b \Lambda(a,b)=\Lambda(a,b) \ \text{ and } \ 
a\pl_a \Lambda(a,b) {-} b\pl_b \Lambda(a,b)=\Lambda(a,b)
\frac{a {+} b{-} 2\Lambda(a,b)}{a-b} .
\]

Now we are in the position to calculate the first-order correction to
the Liouville equation from the approximate entropy $\bfE_V$ (cf.\
\eqref{eq:whE.V}) and the dual dissipation potential $\Psi^{*}_V$, namely
$\dot \varrho =\rmD_\xi \Phi^{*}_V\big(\varrho , {-}\rmD\bfE_V(\varrho) \big)$, which yields 
\begin{align*}
\dot \rho &= 
\div\left\{ \Big(\rho\,\wh a(\bfc) + \frac1V\big[ \rho\, 
 \wh b_0(\bfc) +  \wh b_1(\bfc) \nabla \rho \cdot (\bfalpha{-} 
 \bfbeta) \big] + O(1/V^2)\Big)
   \big(\bfalpha{-}\bfbeta \big) \right\} ,
\end{align*}
where the coefficients are given by 
\begin{align*}
  &\wh a(\bfc)= \Lambda(k_\FW\bfc^\bfalpha, k_\BW
  \bfc^\bfbeta)  \,(\bfalpha{-}\bfbeta)\vdot\nabla
  E(\bfc)\ = \ k_\FW\bfc^\bfalpha {-}  k_\BW \bfc^\bfbeta, 
\\ 
& \wh b_0(\bfc)=\Lambda(k_\FW\bfc^\bfalpha, k_\BW \bfc^\bfbeta)\,
    (\bfalpha{-}\bfbeta)\vdot \check \bfC \bfdelta 
   - \frac12 (k_\FW\bfc^\bfalpha {-} 
   k_\BW \bfc^\bfbeta) \,\bfalpha \vdot \check \bfC \bfbeta , \\
& \wh b_1(\bfc)= \Lambda(k_\FW\bfc^\bfalpha, k_\BW \bfc^\bfbeta) +
\Upsilon_1(\bfc) (\bfalpha{-}\bfbeta){\vdot}\nabla E(\bfc)
 \ = \ \frac12(k_\FW\bfc^\bfalpha{+}k_\BW\bfc^\bfbeta) .
\end{align*}
It is interesting to see the cancellation in the term $\wh b_1$, where
$\Upsilon_1$ did not have a sign, but after multiplication with
$(\bfalpha{-}\bfbeta){\cdot}\nabla E(\bfc)$ it becomes positive and
increases the logarithmic mean
$ \Lambda(k_\FW\bfc^\bfalpha, k_\BW \bfc^\bfbeta)$ to the arithmetic
mean $\frac12(k_\FW\bfc^\bfalpha{+}k_\BW\bfc^\bfbeta)$.  Moreover, the
coefficient $b_0$ consists of two terms, the first of which
corresponds (up to order $1/V^2$) to the correction $\bfA_V$ in
\eqref{eq:FP-V} arising from the improvement of $\bfE_V$, while the
second term arises from improving the dissipation potential
$\Phi^*_V$, namely via $\Upsilon_0$.

Putting these derivations together, summing over $r=1,\ldots,R$ different
reactions, and dropping all terms of order $1/V^2$, we find the
following approximative Fokker--Planck equation:
\begin{align}  
\label{eq:FokPla-3}
  \dot \rho(t,\bfc) &= \div_{\!\bfc}\Big( 
\frac1V    \wh\bbK_\CLE(\bfc) \nabla\rho(t,\bfc) 
+ \rho(t,\bfc) \bfR(\bfc)+  \frac1V \rho(t,\bfc)\bfB(\bfc) \Big)
\end{align}
where $\bfR(\bfc)=\bbK(\bfc)\rmD E(\bfc)$, \
$\bfB(\bfc)=\sum_{r=1}^R \wh b_0^r(\bfc) (\bfalpha^r{-}\bfbeta^r)$,
and $\wh\bbK_\CLE$ is given in \eqref{eq:whK.CLE}.

The big disadvantage of equation \eqref{eq:FokPla-3} is that it is
generally no longer a gradient system.  However, it may be considered
as an equation with an asymptotic gradient flow structure in the sense
of \cite{BBRW17CDSE}.  To find the simplest true gradient system that
is compatible with the Fokker--Planck equation \eqref{eq:FokPla-3}, we
have to find a true dual dissipation potential $\wh\Phi_V^*$ that is
non-negative and coincides with $\Phi^*_V$ from above to lowest order.
To keep the notation light, we again explain the construction for the
case of one reaction only and set
$\Lambda_0(\bfc)=\Lambda(k_\FW\bfc^\bfalpha, k_\BW \bfc^\bfbeta)$. Our
simplest choice is
\begin{align*}
\wh\Phi_V^*(\rho,\xi)=\int_\bfC \Big(
 & \Lambda_0(\bfc) \rho(\bfc)+ \frac1V \Upsilon_0(\bfc) \rho(\bfc) 
  + \frac1V \Upsilon_1(\bfc)\nabla\rho(\bfc){\cdot}(\bfalpha{-}\bfbeta)\\
 &+ \frac{\Upsilon_2(\bfc)}{V^2} \,\rho(\bfc)
  + \frac{\Upsilon_3(\bfc)}{V^2} \frac{\big(\nabla\rho(\bfc){\cdot} 
  (\bfalpha{-}\bfbeta)\big)^2}{ \rho(\bfc)}  
\Big)
  \big(\nabla\xi(\bfc){\vdot}(\bfalpha{-}\bfbeta)\big)^2 \dd\bfc,
\end{align*}
where the higher-order corrections $\Upsilon_2(\bfc) $ and
$\Upsilon_3(\bfc) $ need to be chosen such that
$\wh\Phi_V^*(\rho,\xi)$ is still coercive. Choosing $\theta_1,\theta_2
\in {]0,1[}$ with $\theta_1<\theta_2$, we may require
\[
  \Lambda_0(\bfc) +  \frac{\Upsilon_0(\bfc)}{V} 
  +\frac{\Upsilon_2(\bfc)}{V^2} \geq \theta_2 \Lambda_0(\bfc) \quad
  \text{and}\quad 
  4\theta_1 \Lambda_0(\bfc) \Upsilon_3(\bfc) \geq \Upsilon_1(\bfc)^2 
\]
for all $V>1$, so that
$\wh\Phi_V^*(\rho,\xi) \geq (\theta_2{-}\theta_1)\int_\bfC
\Lambda_0(\bfc)\rho(\bfc) \big(\nabla\xi(\bfc){\vdot}
(\bfalpha{-}\bfbeta)\big)^2 \dd\bfc$.  The bounds for
$\Upsilon_2(\bfc) $ and $\Upsilon_3(\bfc) $ hold for the following
choices (or any bigger ones)
\[
  \Upsilon_2(\bfc) = \frac{\Lambda_0(\bfc)}{16(1{-}\theta_2)}
  \big(\bfalpha{\vdot}\check\bfC \bfbeta\big)^2 \quad\text{and}\quad
  \Upsilon_3(\bfc) = \frac{1}{4\theta_1\Lambda_0(\bfc)}
  \,\Upsilon_1(\bfc)^2. 
\]

Of course, we fix the energy functional to be the improved
entropy functional $\bfE_V$ from \eqref{eq:whE.V}, and the gradient
system $(\PM(\bfC),\bfE_V,\wh\Phi^*_V)$ has the associated
gradient-flow equation
$\dot \varrho = \rmD_\xi \wh\Phi^*_V(\varrho,{-}\rmD
\bfE_V(\varrho))$. With   $\rmD \bfE_V(\varrho)= \frac1V (1 {+} \log
\rho) + E + \frac1V E^V_1$  we find
\begin{equation}
  \label{eq:FokPla-5}
  \begin{aligned}
\dot \rho =\div&\left(\Big[ \wh a_0^V(\bfc)\rho +
\frac{\wh a_1^V(\bfc)}V \nabla_\bfgamma\rho + 
\frac{\wh a_2^V(\bfc)}{V^2} \frac{(\nabla_\bfgamma \rho)^2}\rho 
+ \frac{\wh a_3^V(\bfc)}{V^3}  
\frac{(\nabla_\bfgamma \rho)^3}{\rho^2}  
\Big]\bfgamma  \right) \\
\text{with }\ &\wh a_0^V= \Lambda_\Upsilon^V \:
 \big(\nabla_\bfgamma E + \frac1V \nabla_\bfgamma E^1_V\big),
\quad
\wh a_1^V= \Lambda_\Upsilon^V +  \Upsilon_1\:
 \big(\nabla_\bfgamma E + \frac1V \nabla_\bfgamma E^1_V\big),
\\
&\wh a_2^V=\Upsilon_1{+}\Upsilon_3 \:
\big(\nabla_\bfgamma E  +\frac1V \nabla_\bfgamma E^1_V\big),
 \quad \text{ and }\  \wh a_3^V(\bfc)= \Upsilon_3,
 \\
 \text{where }
 &\Lambda_\Upsilon^V(\bfc)=\Lambda_0(\bfc)
  +\frac{\Upsilon_0(\bfc)}V + \frac{\Upsilon_2(\bfc)}{V^2},
  \quad \bfgamma=\bfalpha{-}\bfbeta, \quad \text{and }
  \nabla_\bfgamma f= \nabla f\vdot \bfgamma.
\end{aligned}
\end{equation}
Because $\nabla_\gamma E^V_1$ is of order $1/V$, we see that this
equation involves terms up to order $1/V^4$, namely through $\wh
a{}^V_0$ and through $\wh a{}^V_2/V^2$.

Clearly, our gradient-flow equation \eqref{eq:FokPla-5} is much more
complicated than those generated by the asymptotic gradient-flow
structures in the sense of \cite{BBRW17CDSE}, where higher order terms
are simply dropped.

There is also the question of well-posedness for equation
\eqref{eq:FokPla-5}. To have parabolicity of the leading terms we need
that the mapping $p \mapsto \frac1V\wh a^V_1 p + \frac1{V^2}\wh a^V_2
p^2 + \frac1{V^3}\wh a^V_3 p^3$ is monotone, which amounts to asking
that $\wh a^V_1 + 2 \wh a^V_2 q + 3\wh a^V_3 q^2 \geq 0$ for all $q
\in \R$. This can be always be achieved by making $\Upsilon_2$
very big while keeping $\Upsilon_3$ constant, since $\Upsilon_2$ only
enters once via $\wh a{}^V_1$.

\subsection{Comparison of models}
\label{su:Comp.CME.Lio.FP}

To appreciate the positive and negative aspects of the different
approximations of the CME, we treat the simplest example, namely the
linear RRE on $\bfC={[0,\infty[}$:
\begin{equation}
  \label{eq:linRRE}
  \dot c = 1 -c \qquad \text{corresponding to the reaction pair }
  X\REVER{1}{1} \emptyset . 
\end{equation}
Obviously, we have the explicit solution $c(t)=1+(c(0){-}1\big)
\ee^{-t}$. 

The associated CME for $\bfu=\PM(\N_0)$ is given by 
\begin{equation}
  \label{eq:SimpleCME}
  \dot u_n = V u_{n-1} -\big(V{+}n\big) u_n +(n{+}1) u_{n+1} \quad
  \text{for } n \in \N_0,
\end{equation}
where $u_{-1}=0$. Using the linearity in \eqref{eq:linRRE}, which leads
to the linearity in $n$ of the coefficients in \eqref{eq:SimpleCME},
we obtain explicit closed form relations of the evolution of the
rescaled expectation $\wh e(t):= \frac1V\sum_{n \in \N_0} nu_n(t)$ and
variance $\wh v(t):= \frac1{V^2} \sum_{n\in \N_0} n^2 u_n - \wh
e(t)^2$, namely
\begin{equation}
  \label{eq:MomSimpleCME}
  \dot{\wh e}(t)= 1 - \wh e(t) \quad \text{ and } \quad 
\dot {\wh v}(t) = -2 \wh v(t) + \frac{1 + \wh e(t)}V.
\end{equation}
Moreover, it can be easily checked that for any solution $t \mapsto
c(t)$ of the RRE \eqref{eq:linRRE} the following formula provides the
explicit solution of the CME \eqref{eq:SimpleCME}:
\begin{equation}
  \label{eq:Sol.CME}
   u_n(t) = \frac{\ee^{-c(t)V}}{n!} \,\big(c(t)V\big)^n 
       \quad \text{for } n \in \N_0. 
\end{equation}
Note that this is expression is compatible with the ODEs \eqref{eq:MomSimpleCME} for
the moments, since for these Poisson distributions we have $\wh
e(t)=c(t)$ and $\wh v(t)=c(t)/V$. 

The Liouville equation and the simple Fokker--Planck equation read
\[
\text{(Lio)}\quad \dot\varrho =\pl_c\big( (c{-}1)\varrho\big) \quad
\text{\ \ and \ \ (FP)} \quad
\dot \rho = 
\pl_c \big(   \Lambda(1,c) \,\frac{ \pl_c\rho}V +
 (c{-}1)\rho \big).
\]
The Fokker--Planck equation for the chemical Langevin equation (cf.\
\eqref{eq:FP.CLE}) takes the form
\[
\text{(FP$_\CLE$)} \quad \dot \rho
 = \pl_c^2 \Big( \frac{1{+}c}{2V} \rho \Big)\,
		+ \pl_c \big( (c{-}1) \rho \big).
\]
To compare the solutions of (FP) and (FP$_\CLE$) with the true
solutions of the CME \eqref{eq:SimpleCME}, we assume that the
solutions can be approximated by Gau{\ss}ians.  In general, for
multidimensional Fokker--Planck equations of the form
$\dot \rho = \frac1V \sum_{ij} \partial_{ij}^2\big( \rho
\bbM_{ij}(\bfc) \big) + \div\big( \rho \bfR_V\big)$ the ansatz
$\rho(t,\cdot) \sim \rmN(\bfa(t), \frac1V \bbA(t))$ with
$\bfa(t)\in \R^d$ and $\bbA(t) \in \R^{d\ti d}_{\mathrm{spd}}$ leads
to the necessary conditions
\[
\dot\bfa(t) = - \bfR(\bfa(t)) \quad \text{ and } \quad 
\dot \bbA(t)= - \rmD \bfR(\bfa(t)) \bbA(t) - \bbA(t) \rmD \bfR(\bfa(t))^\tra +
2 \bbM(\bfa(t)),
\] 
see \cite{SanStu17GASN} for rigorous results of this type. 
Applying these formulas to (FP$_\CLE$) we obtain 
\begin{align}\label{eq:mean-variance:CLE}
	\dot a = 1 - a 
	\quad \text{and} \quad
	\dot A = -2 A + 1 + a,
\end{align}
hence the ODEs for $a$ and $A/V$ coincide with those for ${\wh e}$ and
${\wh v}$ in \eqref{eq:MomSimpleCME}.

A similar argument indicates that solutions to (FP) are well
approximated by Gau{\ss}ians with mean $a_V$ and variance $A_V$
satisfying
\begin{align}\label{eq:mean-variance}
	\dot a_V = 1 - a_V + \frac1V \partial_2 \Lambda(1, a_V) 
		\quad \text{and} \quad
	\dot A_V = -2 A_V + 2\Lambda (1, a_V).	
\end{align}

On the one hand, this clearly indicates that (FP$_\CLE$)
provides a better approximation to the CME for $t\in [0,T]$.  By
formally passing to the limit $V \to \infty$ in
\eqref{eq:mean-variance}, we see that the ODE for $a_V$ is
asymptotically correct.  This is not the case for the ODE for $A_V$,
since the arithmetic mean in \eqref{eq:mean-variance:CLE} is replaced
by the logarithmic mean in \eqref{eq:mean-variance}.  However, the
error of $\Lambda(1,c)$ compared to $\frac12(1{+}c)$ is less than
10\,\% for $c\in [1/3,3]$ and it converges to $0$ for $c\to 1$, i.e.,\
in the limit $t\to \infty$.  Equations \eqref{eq:mean-variance:CLE}
are consistent with Kurtz' central limit theorem, which asserts that
the normalized process $\frac1V \bfN^V(t)$ has fluctuations around
$\bfc(t)$ of order $1/\sqrt{V}$, and the rescaled process
$\sqrt{V}\big( \frac1V \bfN^V(t) -\bfc(t)\big)$ converges to a
Gau{\ss}ian process $t\mapsto \bfV(t)$ with covariance matrix $\bbA$
satisfying
$\dot \bbA(t) = -\rmD\bfR(\bfc(t))\bbA(t) -
\bbA(t)\rmD\bfR(\bfc(t))^\tra + 2 \wh\bbK_\CLE(\bfc(t))$, see, e.g.,\
\cite[Eqn.\,(1.9)]{Kurt78SATD}.
 
On the other hand, (FP) makes a better prediction for the equilibrium
distribution that is attained for $t\to \infty$. For (FP$_\CLE$) we
have the unique steady state 
\[
\rho^\text{eq,CLE}_V(c)= \frac1{Z^\CLE_V} \:\ee^{-V \wt E(c)} \quad
\text{with } \wt E(c)= \int_1^c \!\tfrac{2b{-}2{+} 1/V}{b+1} \dd b
=2c{-}2{-}(4{-}\tfrac1V) \log\tfrac{1{+}c}2.
\]
Thus, $\wt E$ grows only like $c$, such that $\rho^\text{eq,CLE}_V$
decays exponentially only. In contrast, the equilibrium
$\rho^\text{eq,FP}_V= Z_V^{-1} \ee^{-V E(c)}$ of  (FP) produces the
correct super-exponential 
decay of the stationary Poisson distribution equation for the CME
\eqref{eq:SimpleCME}.

\subsection{Approximation via cosh-type gradient structure}
\label{su:cosh.Approx}

The derivation of a gradient structure \eqref{eq:KE-gradflow} for the
Liouville equation \eqref{eq:Liouville} can be repeated very similarly
by starting from the cosh-type gradient structure introduced in
\cite{MiPeRe14RGFL}, see Proposition \ref{prop:CME-coshGS}. We do not
give the details here but provide the result only.

Starting from the cosh-type dual dissipation potential 
$\bfPsi^*_{\cosh,V}$ defined in \eqref{eq:Psi.V.cosh} instead of the
quadratic dual potential $\Psi^*_V$ defined in \eqref{eq:Psi*Vu-xi} 
we obtain the counterparts to
Propositions \ref{pr:Psi*V} and \ref{pr:Psi*Limsup} but now with 
\[
\bfPsi^*_{\cosh,\Lio }(\varrho,\xi):= \int_{\bfc\in \bfC} \sum_{r=1}^R
\kappa^r \Big(\frac{\bfc^{\bfalpha^r}}{\bfc_*^{\bfalpha^r}} \,
\frac{\bfc^{\bfbeta^r}}{\bfc_*^{\bfbeta^r}} \Big)^{1/2} 
\sfC^*\big( ( \bfbeta^r{-}\bfalpha^r) \vdot \nabla_\bfc
\xi(\bfc)\big) \rmd\varrho(\bfc) .
\]

Without any need to justify the approximation procedure in the sense of Section
\ref{su:CME-Liouv} we easily obtain the following result.

\begin{proposition}[cosh-type gradient structure for the Liouville
  equation] The Liouville equation \eqref{eq:Liouville} has the
  gradient structure $(\PM(\bfC), \bfE,\bfPsi^*_{\cosh,\Lio })$
  with $\bfE$ from \eqref{eq:bfE.def} and
  $\bfPsi^*_{\cosh,\Lio }$ from above.
\end{proposition}
\begin{proof} The result follows by using
  $\rmD\calE(\varrho)(\cdot) = E(\cdot)$, \ $\nabla_\bfc
  E(\bfc)=\big( \log (c_i/c^*_i)\big)_{i=1,\ldots,I}$, and 
\[
\rmD_\xi \bfPsi^*_{\cosh,\Lio }\big(\varrho,
{-}\rmD\calE(\varrho)\big)  [\eta] = \int_\bfC \sum_{r=1}^R
\kappa^r \Big(\frac{\bfc^{\bfalpha^r}}{\bfc_*^{\bfalpha^r}} \,
\frac{\bfc^{\bfbeta^r}}{\bfc_*^{\bfbeta^r}} \Big)^{1/2} 
(\sfC^*)'\big({-}\bfgamma^r\vdot \big( \log\frac{c_i}{c^*_i}\big)_{i}
\big)\big[ \bfgamma^r\vdot \nabla \eta\big] 
\rmd\varrho(\bfc),     
\]
where $\bfgamma^r = \bfalpha^r {-} \bfbeta^r$.  Using
$\sqrt{ab}\,(\sfC^*)'\big(\log(a/b)\big)=a{-}b$ and the definition of
$\bfR$ gives $\rmD_\xi \bfPsi^*_{\cosh,\Lio }\big(\varrho,
{-}\rmD\calE(\varrho)\big) [\eta]= -\int_\bfC \bfR(\bfc)\vdot
\nabla\eta \dd\varrho(\bfc)$ which is the desired right-hand side of
\eqref{eq:Liouville} when testing with $\eta$ and integrating by
parts.
\end{proof}

As in the case of quadratic gradient structure for the Liouville
equation we may consider the first-order correction to obtain a
Fokker--Planck equation. For this we insert the improved energy
$\bfE_V$ defined in \eqref{eq:whE.V} into the dissipation potential
$\bfPsi^*_{\cosh,V}$ (cf.\ \eqref{eq:Psi.V.cosh}) to obtain a
quasilinear Fokker--Planck-type equation, namely $\dot\varrho =
\rmD_\xi\bfPsi^*_{\cosh,V} \big(\varrho,{-} \rmD
 \bfE_V (\varrho)\big)$. Using the abbreviations
$a_r:=\frac{\bfc^{\bfalpha^r}}{\bfc_*^{\bfalpha^r}} $ and $b_r:=
\frac{\bfc^{\bfbeta^r}}{\bfc_*^{\bfbeta^r}}$ we find (note
$\nabla_\bfc E^V_1(\bfc)=O(1/V)$) 
\begin{align*}
&\rmD_\xi\bfPsi^*_{\cosh,V} \big(\varrho,{-} \rmD
  \bfE_V (\varrho)\big)  \ = \  
\rmD_\xi\bfPsi^*_{\cosh,V} \Big(\varrho,  -\frac1V \log
                 \rho {-} E{-} \frac1V E^V_1 \Big) \\
&= \div\bigg(\rho \sum_{r=1}^R \kappa^r \sqrt{a_rb_r}
   \Big[(\sfC^*)'\big(\log\frac{b_r}{a_r}\big)\bfgamma^r +
   (\sfC^*)''\big(\log\frac{b_r}{a_r} \big) 
  \frac{\bfgamma^r{\vdot}\nabla \rho}{V \rho} \bfgamma^r\Big]
 +O(1/V^2)  \bigg).
\end{align*} 
Using the identities $\sqrt{ab}\,(\sfC^*)'\big(\log (b/a)\big)=b-a$
and $\sqrt{ab}\,(\sfC^*)''\big(\log (b/a)\big) = (a{+}b)/2$ the
FP equation has the expansion
\[
\dot\rho(t,\bfc) = \div_{\!\bfc}\Big( \rho(t,\bfc) \bfR(\bfc)+
\frac1V \wh\bbK_\CLE(\bfc) \nabla_{\!\bfc} \rho(t,\bfc) +O(1/V^2)_{V\to
  \infty} \Big) 
\] 
where $\wh\bbK_\CLE$ is exactly the same as obtained in \eqref{eq:whK.CLE}
by a completely different approach.

\section{Hybrid models}
\label{s:Hybrid}

We show in this section how the different gradient structures for RRE,
for CME, and for the FPE can be combined to obtain hybrid models,
which are combinations of several models depending on the desired
accuracy.  The importance here is to use the proper rescalings in
terms of the volume $V$ to make the different descriptions compatible.
We do not consider a full theory, but highlight first the general
strategy of model reduction for gradient systems in Section
\ref{su:Reduction} and then illustrate this by a simple example 
in Section \ref{su:Exa.CME2RRE}.  A nontrivial case of a rigorous
coarse graining in this spirit is given in \cite{MieSte19?CGED}, where
a linear RRE with a small parameter $\eps$ is considered. The
elimination of the fast relaxations in the time scale $\eps$ leads to
a coarse-grained gradient system.

In Section \ref{su:Hyb.RRE.FP} we discuss the general coupling of
the FPE to a RRE and the similar coupling of the CME to a RRE, both
leading to so-called mean-field equations, where a linear equation for
a probability density is nonlinearly coupled to an ODE.  Finally, we
discuss the mixed discrete and continuous description, where the CME
is used for small numbers of particles and the FPE is used for larger
numbers.

\subsection{Coarse graining for gradient  systems}
\label{su:Reduction}

If a gradient system $(\sfX,\sfE_X,\Psi_X)$ is more complicated than
what is needed, one is interested in approximating the system by a
simpler model that still contains the most important features.  We
explain how this can be done while keeping the gradient structure.

We assume that the relevant states $\sfx\in \sfX$ can be described by
states $\sfy \in \sfY$ and that there is a reconstruction mapping
$\sfx= \Phi(\sfy)$, i.e.,\ $\Phi(\sfY)$ is a subset (or submanifold)
of $\sfX$. We now pull back the gradient structure
$(\sfX,\sfE_X,\Psi_X)$ to an approximative gradient structure
$(\sfY,\sfE_Y,\Psi_Y)$.  The natural approach is to restrict the
energy functional and the (primal) dissipation potential as follows:
\begin{equation}
  \label{eq:GS.restrict}
  \sfE_\sfY(\sfy)=\sfE_\sfX(\Phi(\sfy)) \quad\text{and} \quad 
\Psi_\sfY(\sfy,\dot\sfy):=\Psi_\sfX(\Phi(\sfy),\rmD \Phi(\sfy)\dot\sfy).
\end{equation}
The solutions $\sfy:[0,T]\to \sfY$ of the coarse-grained gradient
system $(\sfY,\sfE_\sfY,\Psi_\sfY)$ will provide good approximations
$\wh x: t\mapsto \Phi(\sfy(t))\in \sfX$ of the true solutions of the
full GS $(\sfX,\sfE_X,\Psi_X)$, if the set $\Phi(\sfY)$ approximates a
flow-invariant subset of $\sfX$.

In reaction systems, the primal dissipation potential $\Psi_\sfX$ is
usually not known explicitly. Hence, it is desirable to have a method
for reducing the dual dissipation potential $\Psi_\sfX^*$ directly to
$\Psi_\sfY^*$, in the case where $A = \rmD\Phi(\sfy) : \sfY \to \sfX$
is injective but its adjoint mapping $A^*:\sfX^*\to \sfY^*$ has a
large kernel. The following exact result will be the motivation for
our modeling approximations in the subsequent subsections.

\begin{proposition}
\label{pr:Psi*.restr}
Consider reflexive Banach spaces $\sfX$ and\/ $\sfY$ and a real-valued
dissipation potential $\Psi:\sfX\to {[0,\infty[}$ (i.e.\ lower
semicontinuous, convex, and $\Psi(0)=0$) that is superlinear, i.e.\
$\Psi(\mathsf v)/ \|\mathsf v\|_\sfX \to \infty$ for
$\|\mathsf v\|_\sfX\to \infty$.  Assume that  the
bounded linear operator $A : \sfY \to \sfX$ has closed range. 
Then the dissipation potential
$\wt\Psi:\sfY\to {[0,\infty[};\ \sfy \mapsto \Psi(A\sfy)$ satisfies
  \begin{equation}
    \label{eq:Psi*.reduced}
    \wt\Psi{}^*(\eta)= \inf\bigset{\Psi^*(\xi)}{ A^*\xi =\eta} \quad
    \text{ for all } \eta \in \sfY^*,
  \end{equation}
where we use the convention $\inf\emptyset = \infty$.   
\end{proposition}

\begin{proof}
  For the proof we use the saddle-point theory in
  \cite[Ch.\,VI.2]{EkeTem76CAVP}.

  Fix $\eta \in \sfY^*$ and assume first that
  $\eta \notin \Ran (A^*)$. Since $\Ran(A) \subset \sfX$ is closed,
  the Closed Range Theorem yields that $\Ran(A^*) \subset \sfY^*$ is
  closed as well, and $\Ran(A^*) = \Ker(A)^\perp$. Consequently, there
  exists $\tilde \sfy \in \Ker(A)$ such that
  $\ip{\eta, \tilde \sfy} \neq 0$, and we obtain
  \begin{align*}
	\wt\Psi{}^*(\eta)
	= \sup_{\sfy \in \sfY} 
			\big(
				\langle \eta, \sfy \rangle 
				- \Psi(A\sfy) 
			\big) 
	\geq \sup_{\lambda \in \R} 
			\big(
				\lambda \langle \eta, \tilde \sfy \rangle 
				- \Psi(\lambda A \tilde\sfy) 
			\big) 
	= \infty.			
  \end{align*}
  This yields \eqref{eq:Psi*.reduced}, since the right-hand side is
  clearly infinite as well.

  Fix now $\eta \in \Ran(A^*)$ and define the Lagrangian function
  $L:\sfX\ti \sfX^*\to {[{-}\infty,\infty[}$ via
  \[
    L(\sfx,\xi) = -\langle \xi,x\rangle + \Psi(x) - \chi^*(\xi) \quad
    \text{with }\chi^*(\xi)=\left\{ \ba{cl} 0 &\text{for }A^*\xi=\eta,\\
      \infty&\text{otherwise}. \ea\right. 
  \]
  For notational convenience we set
  \[
    h(\sfx)=\sup_{\xi\in \sfX^*} L(\sfx,\xi), \quad g(\xi)=
    \inf_{\sfx\in \sfX} L(\sfx,\xi), \quad P:=\inf_\sfX h, \quad D:=
    \sup_{\sfX^*} g.
  \]
  Classical duality theory yields the trivial inequality $P \geq D$.
  Clearly, $L(\cdot,\xi)$ is convex and lower semicontinuous, whereas
  $L(\sfx,\cdot)$ is concave and upper semicontinuous, since the
  boundedness of $A^*$ implies that
  $\bigset{\xi \in \sfX^*}{ A^*\xi =\eta}$ is closed.

  Using $\eta \in \Ran(A^*)$, we find $\xi_\eta \in \sfX^*$ with
  $A^* \xi_\eta= \eta$, so that our assumptions guarantee the
  coercivity of $\sfx \mapsto L(\sfx, \xi_\eta) \in \R$. Hence, we can
  apply \cite[Chap.\,VI, Prop.\,2.3]{EkeTem76CAVP}, which shows that
  there is no duality gap:
  \begin{align}
    \label{eq:Saddle.a}
    & P=\inf_{\sfx\in \sfX}h(\sfx) =
      \min_{\sfx\in \sfX} \Big(\sup_{\xi\in \sfX^*} L(\sfx,\xi)\Big)
      \ = \  \sup_{\xi\in \sfX^*} \Big( \inf_{\sfx\in \sfX}
      L(\sfx,\xi)\Big) = \sup_{\xi\in \sfX^*} g(\xi) = D.
  \end{align}
  We relate $P$ and $D$ with the two sides in our desired formula
  \eqref{eq:Psi*.reduced}. On the one hand,
  \begin{align*}
    h(\sfx)
    &= \sup_{\xi\in \sfX^*} L(\sfx,\xi) = \Psi(\sfx) +
      \sup_{\xi\in \sfX^*} \big( \langle \xi,{-}\sfx\rangle -
      \chi^*(\xi)\big)   \\
    &= \Psi(\sfx) -\langle \xi_\eta,\sfx\rangle +
      \mu({-}\sfx) \ \text{with } \mu(\sfx):= 
      \sup_{\zeta\in \sfX^*} \big( \langle \zeta,\sfx\rangle -
      \delta_0(A^*\zeta)\big),  
  \end{align*}
  where in the last step we have substituted $\xi= \xi_\eta+\zeta$
  with $A^*\xi_\eta= \eta$ and introduced $\delta_0(\wt\eta)=0$ for
  $\wt\eta=0$ and $\infty$ otherwise. Thus, we conclude
  \[
    h(\sfx) = \Psi(A\sfy) - \langle \eta,\sfy\rangle \ \text{ for }
    \sfx=A\sfy \quad \text{and} \quad h(\sfx)=\infty \ \text{ for }
    \sfx \not\in \Ran(A).
  \]
  Thus, taking the minimum over all of $\sfX$ is the same as taking it
  over $\Ran(A)$, namely
  \[
    P = \inf_{\sfx \in \sfX} h(\sfx) = \inf_{\sfy\in \sfY} \big(
    \Psi(A\sfy)-\langle \eta, \sfy\rangle\big) = - \wt\Psi{}^*(\eta).
  \]

  On the other hand, the definition of
  $g(\xi)=\inf_{\sfx \in \sfX} L(\sfx, \xi)$ immediately gives
  $g(\xi)=-\Psi^*(\xi)-\chi^*(\xi)$. Hence, we arrive at
  \[
    D = \sup_{\xi\in \sfX^*} g(\xi) = - \inf_{\xi\in \sfX^*}
    \big(\Psi^*(\xi)+\chi^*(\xi)\big)= - \inf\bigset{\Psi^*(\xi)}{
      A^*\xi=\eta}.
  \]
  As a result, formula \eqref{eq:Psi*.reduced} follows from $P=D$.
\end{proof}

In our applications below (as well as in most others) the explicit
minimization in \eqref{eq:Psi*.reduced} is too complicated to be
executed.  However, as the coarse-graining mapping through $\Phi$ is
usually only an approximation, it may suffice to approximate the
minimizers suitably.  In general, one has to find an approximation
$\xi=\sfM(\sfy,\eta) \in \sfX^*$ and sets
\begin{align}
\label{eq:GS.restrict.K}
&\Psi^*_\sfY(\sfy,\eta)= \Psi^*_\sfX(\Phi(\sfy),\sfM(\sfy,\eta)) \quad
\text{or} \quad  \tfrac12   
\langle \eta, \wt\sfK_\sfY(\sfy) \eta\rangle =
\langle \mathsf L(\sfy) \eta, \sfK_\sfX(\Phi(\sfy)) \mathsf L(\sfy) \eta\rangle, 
\end{align}
where $\mathsf L(\sfy):\sfY^*\to \sfX^*$ is the linear version of
$\sfM$.  Of course, when constructing $\sfM$ or $\mathsf L$ one should
keep \eqref{eq:Psi*.reduced} in mind to preserve all interesting
properties inherited by the coarse-graining process.

\subsection{A simple example: from CME to RRE}
\label{su:Exa.CME2RRE}

We apply the above idea with $(\sfX,\sfE_\sfX,\Psi^*_\sfX)$ being
$(\PM(\calN),\calE_V,\calK_V)$ and with $(\sfY,\sfE_\sfY,\Psi_\sfY)$
being $(\bfC,E,\bbK)$.  The embedding mapping
$\Phi_V : \bfC \to \PM(\calN)$ is given by the Poisson distributions
\[
\Phi_V(\bfc):= \Big( \ee^{-V|\bfc|_1}
\,\frac{(V\bfc)^\bfn}{\bfn!}\Big)_{\bfn\in \calN}  ,\ \text{ where
}|\bfc|_1=\sum_{1}^I c_i.
\]
In the simple example $\dot c=1-c$ treated in Section
\ref{su:Comp.CME.Lio.FP} the image of $\Phi_V$ defines an exactly
invariant submanifold, but this is no longer true for nonlinear
equations or systems. Nevertheless our construction provides the
surprising identity
\[
\sfE_\sfY(\bfc)= \calE_V(\Phi_V(\bfc)) = E(\bfc),
\]
with the old $E$ defined in \eqref{eq:RRE.E.K} which is independent of $V$. 

To reduce the dual dissipation potential $\Psi_\sfX$ defined via $\calK_V$ we use the derivative 
\[
\rmD \Phi_V(\bfc) \bfw = \Big(\ee^{-V|\bfc|_1}
\,\frac{(V\bfc)^\bfn}{\bfn!} \sum_{i=1}^I \big(\frac{n_i}{c_i}-V \big)w_i  \Big)_{\bfn\in \calN} .
\]
Thus, the adjoint operator $\rmD \Phi_V(\bfc)^*$ maps
$\bfmu=(\mu_\bfn)$ to $\bfzeta=(\zeta_i)_{i=1,\ldots,I}$ via 
\[
\bfmu\mapsto \bfzeta
	= \rmD \Phi_V(\bfc)^*\bfmu
	= \Big(\sum_{\bfn\in \calN}
			 \ee^{-V|\bfc|_1} \,\frac{(V\bfc)^\bfn}{\bfn!} \big(\frac{n_i}{c_i}-V \big) \mu_\bfn 
\Big)_{i=1,\ldots,I}  .
\]
In general, one is not able to solve the minimization problem \eqref{eq:Psi*.reduced} that produces $\Psi^*_\sfY$ from $\Psi_\sfX^*$, so instead we construct a linear mapping $\bfzeta\mapsto \wt\bfmu=\sfM_V(\bfc)\bfzeta$ that approximates the minimizer for $V\to \infty$ and satisfies $\bfzeta= \rmD\Phi_V(\bfc)^*\mathsf M_V(\bfc)\bfzeta$. Indeed, we search for $\wt\bfmu$ in the linear form
$\wt\bfmu^\bfa_\bfn = \bfa\vdot \bfn$ for $\bfn\in \calN$ and obtain
\begin{align*}
\rmD \Phi_V(\bfc)^*\wt\bfmu^\bfa&=\rmD \Phi_V(\bfc)^*( \bfa{\vdot}
\bfn)_{\bfn\in \calN} = \bigg(\sum_{\bfn\in \calN} \ee^{-V|\bfc|_1}
\,\frac{(V\bfc)^\bfn}{\bfn!} \big(\frac{n_i}{c_i}-V \big) 
  \sum_{j=1}^I a_j n_j \bigg)_{i=1,\ldots,I}\\
& =\bigg(\sum_{j=1}^I \sum_{\bfn\in \calN}  \ee^{-V|\bfc|_1}
       \,\frac{(V\bfc)^\bfn}{\bfn!} \Big(\frac{n_ia_jn_j}{c_i}-Va_jn_j 
      \Big)\bigg)_{i=1,\ldots,I} \\
& =\bigg(\sum_{j=1}^I \Big(V^2\frac{c_i
  a_jc_j}{c_i} +\delta_{ij}Va_i - V^2a_jc_j   \Big) \bigg)_{i=1,\ldots,I}
\ = \ V\bfa, 
\end{align*}
where we used the identities $ \sum_{\bfn\in \calN} \ee^{-V|\bfc|_1}
\,\frac{(V\bfc)^\bfn}{\bfn!} n_i = Vc_i$ and $ \sum_{\bfn\in \calN} \ee^{-V|\bfc|_1}
\,\frac{(V\bfc)^\bfn}{\bfn!} n_i n_j= V^2c_ic_j + \delta_{ij}
Vc_i$. Thus, we choose the simple operator $\sfM_V$ of the form 
\[
\bfzeta \ \mapsto \ \bfmu=\sfM_V(\bfc) \bfzeta = \big(\frac1V
\bfzeta \vdot \bfn \big)_{\bfn\in \calN}.
\]

For inserting $ \mu=\sfM_V(\bfc)\bfzeta $ and $\bfu=\Phi_V(\bfc)$ into
the full dual dissipation potential $\Psi^*_\sfX$, we use the form
\eqref{eq:Psi*Vu-xi.b} and the relations $(\sfM_V(\bfc)
\bfzeta)_{\bfn{+}\bfalpha} -(\sfM_V(\bfc) \bfzeta)_{\bfn{+}\bfbeta} =
\frac1V \bfzeta\vdot(\bfalpha{-}\bfbeta)$ and
\[
\BBneu V{\bfbeta}{\bfn} \big(\Phi_V(\bfc)\big)_{\bfn+\bfbeta} 
=\frac{V (\bfn{+}\bfbeta)!}{V^{|\bfbeta|}\bfn!}\:
  \ee^{-V|\bfc|_1}\frac{(V\bfc)^{\bfn+\bfbeta} }{(\bfn{+}\bfbeta)!} =
  V\, \ee^{-V|\bfc|_1}\frac{(V\bfc)^\bfn}{\bfn!}\;\bfc^{\bfbeta}. 
\]
With this,  we find an
approximation of the reduced dual dissipation potential $\Psi^*_\sfY$,
namely
\begin{align*}
\Psi^*_\sfY(\bfc,\bfzeta) &:= \Psi^*_\sfX\big(\Phi_V(\bfc), \mathsf
M_V(\bfc) \bfzeta \big) \\
& = \frac V2\sum_{r=1}^R \sum_{\bfn\in\calN}
   V \, \ee^{-V|\bfc|_1}\frac{(V\bfc)^\bfn}{\bfn!} \:
   \Lambda\big(k^r_\FW \bfc^{\bfalpha^r},  
    k^r_\BW \bfc^{\bfbeta^r}\big)
    \big( \frac{(\bfbeta^r{-}\bfalpha^r){\vdot}\bfzeta}V\big)^2 \\
& = \frac12 \sum_{r=1}^R  \Lambda\big(k^r_\FW \bfc^{\bfalpha^r},  
    k^r_\BW \bfc^{\bfbeta^r}\big)
    \big( (\bfbeta^r{-}\bfalpha^r){\vdot}\bfzeta\big)^2 \ = \
    \frac12 \bfzeta \vdot \bbK(\bfc) \bfzeta. 
\end{align*} 
Thus, the gradient system $(\sfY,\sfE_\sfY,\Psi^*_\sfY)$ obtained by
the abstract reduction procedure is exactly given by $(\bfC,E,\bbK)$,
which is the gradient system  for the RRE \eqref{eq:ReactKin} studied
in Theorem \ref{thm:RRE.GradSys}.

\subsection{Coupling a RRE to a Fokker--Planck equation}
\label{su:Hyb.RRE.FP}

In many applications one is interested in the microscopic description
of some variables $c_j$, while other variables $c_i$ can be described
more macroscopically. We first start from the simplified FPE
\eqref{eq:FPsimple} as the gradient system $(\PM(\bfC),
\wt\bfE_V, \bfK)$ and partition the components of $\bfc$ into
\emph{stochastic}
and \emph{macroscopic} parts, $\bfc_\rms$ and $\bfc_\rmm$ respectively, via
\begin{align*}
&\bfc=(c_1,\ldots,c_J,c_{J+1},\ldots,c_I) = (\bfc_\rms,\bfc_\rmm) \quad
\text{with }\\
&\bfc_\rms:=(c_1,\ldots,c_J)\in \bfC_\rms:={[0,\infty[}^J 
\text{ and }\bfc_\rmm:=(c_{J+1},\ldots,c_I)\in \bfC_\rmm:={[0,\infty[}^{I-J},
\end{align*}
In the notation of Section \ref{su:Reduction} we let
$\sfX=\PM(\bfC_\rms \ti \bfC_\rmm)$ and $\sfY=\PM(\bfC_\rms)\ti
\bfC_\rmm$. 

For the mapping $\Phi:\sfY\to \sfX$ we choose the product
ansatz
\begin{align*}
\Phi_V(\varrho_\rms,\wh\bfc_\rmm)({\rm d}c_1,\ldots,{\rm d}c_I):=
 \varrho_\rms({\rm d}\bfc_\rms) \prod_{j=J+1}^I \sfW(c_j;\wh c_j,V) \dd c_j,
\end{align*}
where the probability densities $\sfW(\cdot;\wh a,V)$ are given as follows:
\begin{align*}
\sfW(a;\wh a,V)&:=\frac{1}{\wh a\,Z(V\wh a)} \exp\big({-}V \wh a\: 
  \LB(a/\wh a)\big) 
 \ \text{ with } Z(v):=\int_0^\infty \!\!\exp\big({-}v\, 
 \LB(z)\big) \dd z .
\end{align*}
According to Section \ref{su:Reduction} the functional
$\wt\sfE^V_\sfY= \wt\bfE_V \circ \Phi_V$ on $\sfY$ is then given by
\begin{align*}
 \wt\sfE^V_\sfY( \varrho_s,\wh\bfc_\rmm) & =
  \int_{\bfC_\rms} \!\!
  	\Big(\frac1V \rho_s(\bfc_s)\log\rho_s(\bfc_\rms)  
  + \rho_s(\bfc_\rms) E_\rms(\bfc_\rms) \Big)\dd \bfc_\rms
  + \frac{\wt Z(V)}{V}
  +\sum_{j=J+1}^I \wh e_V(\wh c_j,c^*_j) 
  \\
\text{where }& \wh e_V(\wh a,a^*) := A(V\wh a)\wh a \log\big(\frac{\wh
  a}{a^*}\big) -\wh a + a^*- \frac{\log \big(\wh a Z(V\wh a)\big)}{V}
\end{align*}
with $E_\rms(\bfc_\rms) = \sum_{i=1}^J c_i^* \LB(c_i/c_i^*)$ and
$A(v)= \int_0^\infty z \exp\big(-v \LB(z)\big) \dd z/Z(v)$.  

It can be shown that $A(v)\geq 1$ and
$e_V(\wh a, a^*)\geq a^*\LB(\wh a/a^*)$ for all $V$, and for
$V\to \infty$ we obtain $e_V(\wh a, a^*)\to a^*\LB(\wh a/a^*)$. To
simplify the model we are therefore allowed to replace the last term
in $\wt\sfE^V_\sfY$ by the relative entropy
$E_\rmm(\wh\bfc_\rmm) =\sum_{j=J+1}^I c_j^* \LB(\wh c_j/c_j^*)$ for
the RRE.  Neglecting the irrelevant constant term
${\wt Z(V)}/{V}$, we obtain the hybrid energy again as a
relative entropy, namely
\begin{align*}
\mfE^\text{FP-RR}_V(\varrho_s,\wh\bfc_\rmm)
	 & = \int_{\bfC_\rms}\!\!\Big(
		\frac1V \rho_s(\bfc_\rms)\log\rho_s(\bfc_\rms) 
	+ \rho_s(\bfc_\rms) 
		E_\rms(\bfc_\rms)\Big) \dd \bfc_\rms 
	+ E_\rmm(\wh\bfc_\rmm) . 
\end{align*}

For the Onsager operator we also use a cruder  reduction than the
minimization advocated in Section \ref{su:Reduction}. We simply
postulate the Onsager operator $\sfK_V$ via the dual dissipation potential 
\[
\Psi^*_{V,\text{FP-RR}}(\varrho_s,\wh\bfc_\rmm; \xi,\zeta) 
	= \frac12\int_{\bfC_\rms} 
		\rho_s(\bfc_\rms)  
		\binom{\nabla_s \xi(\bfc_\rms)}{\zeta} 
			\cdot \bbK(\bfc_\rms,\wh\bfc_\rmm)
		\binom{\nabla_s \xi(\bfc_\rms)}{\zeta} 
	\dd \bfc_\rms,
\]
where $\xi\in\rmC^1(\bfC_\rms)$ and $\zeta\in \R^{I-J}$. 
Indeed, in the sense of the general reduction method explained in Section \ref{su:Reduction} we see that $\mfK^\text{FP-RR}_V$ is obtained from $\bfK$ by inserting
$\varrho = \varrho_s(\rmd \bfc_\rms){\otimes}\delta_{\wh\bfc_\rmm}$  and 
$\Xi=\sfM(\xi,\zeta):(\bfc_\rms,\bfc_\rmm)\mapsto \xi(\bfc_\rms)+
\zeta{\cdot}\bfc_\rmm $.  

Thus, the hybrid model induced by the gradient system
$(\PM(\bfC_\rms)\ti \bfC_\rmm, \mfE^\text{FP-RR}_V,\mfK^\text{FP-RR}_V)$ is given by
the coupled system for $\rho\in \PM(\bfC_\rms)$ and $\wh\bfc_\rmm\in
\bfC_\rmm$: 
\begin{align*} 
\dot \rho(\bfc_\rms) & = \div_\rms\Big(
\bbK_{\rms\rms}(\bfc_\rms,\wh\bfc_\rmm)\big(\tfrac1V\nabla_\rms \rho(\bfc_\rms) {+}
\rho(\bfc_\rms)\nabla_\rms E_\rms(\bfc_\rms)\big)  + \rho(\bfc_\rms)
\bbK_{\rms\rmm}(\bfc_\rms,\wh\bfc_\rmm)\nabla_\rmm E_\rmm(\wh\bfc_\rmm) \Big) ,\\
  \dot{\wh\bfc}_\rmm  & = - \!\int_{\bfC_\rms}\!\!\!\Big(
  \bbK_{\rms\rmm}^\tra(\bfc_\rms,\wh\bfc_\rmm) \big(\tfrac1V 
\nabla_\rms\rho(\bfc_\rms){+} \rho(\bfc_\rms)\nabla_\rms E_\rms(\bfc_\rms)\big) +
\rho(\bfc_\rms)\bbK_{\rmm \rmm}(\bfc_\rms,\wh\bfc_\rmm) \nabla_{\rmm}
E_\rmm(\wh\bfc_\rmm) \Big) \!\dd \bfc_\rms.  
\end{align*}
It is interesting to see that the last terms can be rewritten in terms
of the RRE $\dot\bfc = - \bbK(\bfc) \rmD E(\bfc)= -\bfR(\bfc)=
-(\bfR_\rms(\bfc_\rms,\wh\bfc_\rmm),\bfR_\rmm(\bfc_\rms,\wh\bfc_\rmm))$,
viz.\  
\begin{align*}
\dot \rho(\bfc_\rms) & = \div_\rms\Big(\tfrac1V
\bbK_{\rms\rms}(\bfc_\rms,\wh\bfc_\rmm)\nabla_\rms \rho(\bfc_\rms) + \rho(\bfc_\rms)
\bfR_\rms(\bfc_\rms,\wh\bfc_\rmm)  \Big)  ,\\ 
 \dot{\wh\bfc}_\rmm & = - \!\int_{\bfC_\rms} \!\!\!
  \Big(\tfrac1V  \bbK_{\rms\rmm}^\tra(\bfc_\rms,\wh\bfc_\rmm) 
     \nabla_\rms\rho(\bfc_\rms)  + \rho(\bfc_\rms) \bfR_\rmm (
     \bfc_\rms, \wh\bfc_\rmm) \Big) \dd \bfc_\rms.  
\end{align*} 
This reveals that the system is a classical mean-field model, which is
linear in the density $\rho$ for the component $\bfc_\rms$ while it is
nonlinearly coupled to the ODE for the component $\wh\bfc_\rmm$.

\subsection{Coupling a RRE to a CME}
\label{su:Hyb.RRE.CME}

In analogy to the coupling of an RRE for some macroscopic $\bfc_\rmm$
to a Fokker--Planck equation we can directly couple the CME to an RRE,
which leads to hybrid system defined on
$\PM(\N_0^J)\ti {[0,\infty[}^{I-J}$.  Instead of given the general
derivation as in Section \ref{su:Hyb.RRE.FP}, we just give an explicit
example.
    
For $\beta \in \N_0$ we consider the simple reaction
$X_1 \REVER{}{} \beta X_2$ with stoichiometric vectors
$\bfalpha=(1,0)$, $\bfbeta=(0,\beta)$, and $\bfgamma=(1,-\beta)$.  The
associated system of RREs is given by
\begin{equation}
  \label{eq:Hyb.RRS}
  \dot c_1 = c_2^\beta - c_1,\quad \dot c_2 = \beta\,(c_1{-}c_2^\beta).
\end{equation}
We have the conservation relation
$\bbQ\bfc = \beta c_1 + c_2 = q$ and the detailed-balance steady state
$\bfc_*=(1,1)^\top$. The associated CME on $\calN=\N_0^2$ takes the
form
\begin{equation}
  \label{eq:Hyb.CME}
  \dot u_\bfn
 = (n_1{+}1) u_{\bfn+(1,-\beta)} 
 - \Big( n_1 + \frac{n_2!}{V^{\beta-1} (n_2{-}\beta)!}c\Big)u_\bfn
   + \frac{(n_2{+}\beta)!}{V^{\beta-1}n_2!} \,u_{\bfn+(-1,\beta)}
   \text{ for } \bfn \in \calN.
\end{equation}
The detailed-balance steady state by $\bfw^V_\bfn = (w_{n_1}^V
w_{n_2}^V)_{\bfn \in \calN}$ with $w^V_n = \ee^{-V}V^n/n!$. 
As in \eqref{eq:Psi*Vu-xi.b} the full Onsager operator $\calK_V$ is defined via
\[
\ip{\bfmu,\calK_V(\bfu)\bfmu}
= \sum_{\bfn\in \calN}  
   \Lambda\big(\tfrac{n_1{+}1}{V}u_{\bfn{+}(1,0)}, 
	 \tfrac{(n_2{+}\beta)!}{V^\beta n_2!} u_{\bfn+(0,\beta)} 
	\big) 
  \big(V(\mu_{\bfn+(1,0)} {-} \mu_{\bfn+(0,\beta)})\big)^2.
\] 

We partition $\bfc=(c_1,c_2)=(c_\rms,c_\rmm)$, i.e.,\ we keep $c_1\in
{[0,\infty[}$ in stochastic description via the distribution
$\bfv=(v_m)_{m \in \N_0}\in \PM(\N_0)$, while $c_2\in {[0,\infty[}$
will be treated macroscopically. Thus, we define the gradient system
$(\PM(\N_0)\ti{[0,\infty[}),\mfE_V^\text{CM-RR},\mfK_V^\text{CM-RR})$
with relative entropy 
and Onsager operator defined via 
\begin{align*}
&\mfE^\text{CM-RR}_V(\bfv,c_2)= \sfE(c_2)+ \frac1V \sum_{m\in\N_0} v_m\log (v_m/w^V_m) ,
\quad \text{where } \sfE(z)=\LB(z),
\\ 
& \ip{\tbinom\bfxi\zeta, \mfK^\text{CM-RR}_V(\bfv,c_2)\tbinom\bfxi\zeta} = 
V\sum_{m\in \N_0} 
	\Lambda\big( 
			\tfrac{m{+}1}{V}v_{m+1} ,
			v_m c_2^\beta
		   \big) \big( 
\xi_m{+}\tfrac{\beta}{V} \zeta - \xi_{m+1}\big)^2,
\end{align*}
for $\bfxi : \N_0 \to \R$ and $\zeta \in \R$.  Again,
$\mfK^\text{CM-RR}_V$ is obtained from $\calK_V$ by inserting
$\bfu_{m,n_2}= v_m \delta_{\lfloor Vc_2 \rfloor}(n_2)$  and
$\Xi = \sfM (\bfxi,\zeta): (m,n_2) \mapsto \xi_m +\frac1V n_2\zeta$
and performing an approximation for large $V$. The associated
evolution equation is the hybrid system
\begin{align*}
  &\dot v_m= Vc_2^\beta v_{m-1} -\big(m+ Vc_2^\beta\big)v_m +
  (m{+}1)v_{m+1}  
\quad \text{for }m\in \N_0 \ \text{ (with }v_{-1}=0\text{)},\\
  & \dot c_2= \beta \Big( \frac1V \sum_{m\in \N} m v_m - c_2^\beta
  \Big).
\end{align*} 
Clearly, this system is consistent with the conservation law
$\bbQ\bfc = \beta c_1+c_2=\text{const.}$, 
in the sense that $c_1 := \frac{1}{V}\sum_{m\in \N} m v_m$ satisfies
$\dot c_1 = c_2^\beta - c_1 = - \dot c_2 / \beta$.

\subsection{Combining CME and Fokker--Planck descriptions}
\label{su:Hyb:DiscCont}

We consider the simplest nontrivial model, namely the scalar RRE
$\dot c = a - b c$ with $a, b > 0$,  which is induced by the
reaction $\emptyset \REVER{b}{a} X$. This corresponds to
$\alpha = 0$, $\beta = 1$, $k_\FW = a$, and $k_\BW = b$. We have the
following three derived gradient systems:

\begin{enumerate}
\item[(1)]
The RRE $\dot c = a-bc$ is generated by the gradient system $(\R_+,
\bbK,E)$ with steady state $c_*=a/b$, \ $\bbK(c) = \Lambda(a, bc)$,
and $E(c) = \frac ab \LB(bc/a)$.

\item[(2)]
The associated chemical master equation $\dot\bfu =\calB_V \bfu$ is generated by the gradient
system $(\PM(\N_0),\calE_V,\calK_V)$ and reads 
\begin{equation}
  \label{eq:Hyb.discr}
  \dot u_n = Va u_{n-1} - \big( Va {+}bn\big) u_n + b(n{+}1)u_{n+1}
  \quad
\text{for }n\in \N_0 \quad\text{(with $u_{-1}=0$)}
\end{equation} 
and has the steady state $\bfw_V=(\ee^{-Va/b} (Va/b)^n/n!)_{n\in
  \N_0}$. The entropy and Onsager operator are
\begin{align*}
& 
 \calE_V(\bfu)= \frac1V\sum_{n\in \N_0}u_n \log(u_n/w^V_n)
 \quad \text{and}
   \\&\calK_V(\bfu)= V^2a\sum_{n\in \N_0} w^V_n
\Lambda\big(\tfrac{u_n}{w^V_n},\tfrac{u_{n+1}}{w^V_{n+1} }\big)
(\bfe_{n}{-}\bfe_{n+1})\otimes (\bfe_n{-}\bfe_{n+1}). 
\end{align*}

\item[(3)]
The associated Fokker--Planck equation \eqref{eq:FP-V} takes the
form 
\begin{equation}
  \label{eq:Hyb.FP}
  \dot\rho = \pl_c\Big(\frac{\Lambda(a,bc)}V \pl_c\rho +
  \big(
 bc{-}a +\frac{\Lambda(a,bc)}{2Vc{+}1/3} \big)\rho\Big)
  \text{ for }t,c>0\text{ and } \rho(t,0)=0.
\end{equation}
This equation has the equilibrium solution $W_V:c \mapsto \sfW(c;a/b,V)$ (cf.\
\eqref{eq:sfW}) and is generated by the gradient system
$(\PM({]0,\infty[}),\bfE_V,\bfK)$ with
\begin{align*}
&\bfE_V(\rho) = \frac1V\int_0^\infty  \rho(c) \log \Big(\frac{\rho(c)}{W_V(c)}\Big)\dd
c\quad \text{and} \quad
\bfK(\rho)\xi = - \big(\rho \Lambda(a,bc) \xi'\big)'. \medskip
\end{align*}
\end{enumerate}

To combine the description via the CME and the Fokker--Planck equation
we consider the mixed state space $\mfN:=\{0,1,\ldots,N{-}1\}\cup
{[N/V,\infty[}$. Hence, $n\in \{0,\ldots,N{-}1\}$ counts the number of
atoms, while for $n\geq N$ we use the concentration $c=n/V \geq N/V$
as a continuous variable to describe the state. A typical choice could
be $1 \ll V \ll N$ to be sure to capture all small discrete effects. 

The hybrid gradient
system $(\PM(\mfN),\mfE_{V,N},\mfK_{V,N})$ is described by measures 
\[
\mfu= \sum_{n=0}^{N-1} u_n\delta_{n} + U(c) \dd c|_{[N/V,\infty[} \ \in \PM(\mfN).
\]
The idea is now to choose $\mfE_{V,N}$ and $\mfK_{V,N}$ rather than
to model the evolution equation. 

We first choose the equilibrium state in the form 
\[
\mfw^{V,N} = \sum_{n=0}^{N-1} w^V_n\delta_n + W^V(c) \dd c:=
\sum_{n=0}^{N-1} \ee^{-Va/b}\frac{(Va/b)^n}{n!} \delta_{n} + 
\frac1{Z_{V,N}} W_V(c)\,\dd c ,
\]
where $Z_{V,N} $ is uniquely determined by asking $\int_\mfN \dd\mfw^{V,N}
=1$. The entropy functional is defined via the obvious relative
entropy per volume, namely  
\begin{align*}
\mfE_{V,N}(\mfu) &= \frac1V \int_\mfN 
	 \log\Big(\frac{\rmd \mfu}{\rmd\mfw^{V,N}}\Big) 
	 	 \dd \mfu\\
& = \frac1V \sum_{n=0}^{N-1}
\LB\big(\frac{u_n}{w_n^V} \big)w_n^V + \frac1V\int_{N/V}^\infty \LB\big(
\frac{U(c)}{W_V(c)}\big) W_V(c)\dd c
 , 
\end{align*}
where $\frac{\rmd \mfu}{\rmd\mfw} $ denotes the Radon-Nikodym
derivative. 

The difficult part is the modeling of the Onsager operator $\mfK_{V,N}(\bfu)$
as it includes the crucial transfer between the discrete and the
continuous parts of the hybrid model. We define $\mfK$ in terms of its
associated quadratic form acting on smooth functions $\xi:\mfN\to \R$,
where we write $\xi_n$ for $\xi(n)$ and $W(c)$ for
$W_V(c)$: 
\begin{align*}
\langle \xi, \mfK_{V,N}(\bfu)\xi\rangle &= V^2 a
\sum_{n=1}^{N-1}  w_{n-1}\Lambda\big(\tfrac{u_{n-1}}{w_{n-1}} ,
\tfrac{u_n}{w_n}\big) (\xi_{n-1}{-}\xi_n)^2 \\
&\quad + V^2 \wh a w_{N-1}\,\Lambda\big(\tfrac{u_{N-1}}{w_{N-1}},
 \tfrac{U(N/V)}{W(N/V)} \big) \big( \xi_{N-1}{-}\xi(N/V)\big)^2
\\ & \quad
 + \int_{N/V}^\infty \Lambda(a,bc) \xi'(c)^2 U(c) \dd c .
\end{align*}
While the first and the third terms on the right-hand side give the
purely discrete and the continuous parts of the state space,
respectively, we see that the second term is the new term that couples
the discrete and the continuous parts. The parameter $\wh a$ is still
to be chosen, the natural parameter being $a$. 

The evolution equation for $\mfu$ is again a linear equation of the
form $\dot \mfu = \mfB_{V,N} \mfu$, i.e.,\ it corresponds to a continuous-time Markov
process. It consists of a discrete part, as in \eqref{eq:Hyb.discr}
but only for $n=0,\ldots,N-2$, and a continuous part, as in
\eqref{eq:Hyb.FP} but only for $c>N/V$. The new structure is the coupling
between the two subsystems which gives rise to the following
conditions:
\begin{align*}
\dot u_{N-1}&= V a\, u_{N-2} - \big(V\wh a +b\, (N{-}1) \big)
u_{N-1} +  V\wh a \,\frac{w_{N-1}}{W(N/V)}\,U(N/V),  \\
 0 &=  V\wh a\, \Big(\frac{w_{N-1}}{W(N/V)}\,U(N/V) - u_{N-1}
 \Big) + \frac1V W(N/V) \Lambda(a, b N/V)\Big( \frac UW\Big)'(N/V).  
\end{align*} 
By our definition of $\mfw^{V,N}$ we have $\frac{w_{N-1}}{W(N/V)}
\approx Nb/(aV^2)$ and see that for $\wh a =a$ these conditions take
the approximate form 
\begin{align*}
\dot u_{N-1}&\approx V a\, u_{N-2} - \big(V a +b\, (N{-}1) \big)
u_{N-1} +  b \tfrac{N}{V}\,U( \tfrac NV),  \\
 0 &\approx   \tfrac1V \Lambda(a, b \tfrac NV) U'(\tfrac NV) +  b\tfrac{N}{V} 
 U(\tfrac NV) - aV\,u_{N-1},
\end{align*}
where the second relation clearly shows the corresponding Robin
boundary condition connecting the parabolic Fokker--Planck equation to
the discrete system on $\{0,\ldots,N{-}1\}$.  Note that $u_n$ and $U$
are scaled such that $V u_{N-1}$ is comparable to $U(N/V)$.

\paragraph{Acknowledgments.} The research of A.M. was partially
supported by the Deutsche Forschungsgemeinschaft (DFG) via the
Collaborative Research Center SFB\,1114 \emph{Scaling Cascades in
  Complex Systems} through the Subproject C05 \emph{Effective models
  for materials and interfaces with multiple scales}.

J.M. gratefully acknowledges support by the European Research Council
(ERC) under the European Union's Horizon 2020 research and innovation
programme (grant agreement No.\,716117), and by the Austrian Science
Fund (FWF), Project SFB F65. 

The authors thank Christof Sch\"utte, Robert I. A. Patterson, and
Stefanie Winkelmann for helpful and stimulating discussions.

\small
\renewcommand{\baselinestretch}{0.95}
\bibliographystyle{my_alpha}
\bibliography{alex_pub,bib_alex}

\end{document}